\documentclass{gtart} 

\input gtoutput
\volumenumber{3}\papernumber{17}\volumeyear{1999}
\pagenumbers{405}{466}\published{17 December 1999}
\proposed{Joan Birman}\seconded{Walter Neumann, Cameron Gordon}
\received{7 January 1999}\revised{18 November 1999}
\accepted{6 December 1999}

\usepackage{amssymb,leqno}

\input rlepsf
\let\relabela\adjustrelabel
\let\ss\scriptstyle

\def\ints{\mathbb Z}

\def\pp{\prime\prime}
\let\endpf\endproof
\def\lan{\langle}
\def\ran{\rangle}
\def\rig{\rightarrow}
\def\mapright#1{\smash{\mathop{\longrightarrow}\limits^{#1}}}
\def\strut{\vrule width 0pt height 12pt}

\newtheorem{theorem}{Theorem}
\newtheorem{proposition}[theorem]{Proposition}
\newtheorem{lemma}[theorem]{Lemma}
\newtheorem{cor}{Corollary}[theorem]

\theoremstyle{remark}
\newtheorem{remark}{Remark}
\newtheorem{definition}{Definition}

\def\cM{{\mathcal M}}

\begin{document}

\title{An elementary approach to the mapping class\\group of a surface}

\author{Bronislaw Wajnryb}
\address{Department of Mathematics\\Technion, 32000 Haifa, Israel}
\email{wajnryb@techunix.technion.ac.il}

\begin{abstract}
We consider an oriented surface $S$ and a cellular complex $X$ of curves
on $S$, defined by Hatcher and Thurston in 1980. We prove by
elementary means, without Cerf theory, that the complex $X$ is
connected and simply connected. From this we derive an explicit
simple presentation of the mapping class group of $S$, following the
ideas of Hatcher--Thurston and Harer.
\end{abstract}
\asciiabstract{We consider an oriented surface S and a cellular
complex X of curves on S, defined by Hatcher and Thurston in
1980. We prove by elementary means, without Cerf theory, that the
complex $X$ is connected and simply connected. From this we derive an
explicit simple presentation of the mapping class group of S,
following the ideas of Hatcher-Thurston and Harer.}




\maketitlepage

\section{Introduction}
Let $S$ be a compact oriented surface of genus $g\ge 0$ with $n\ge 0$
boundary components and $k$ distinguished points. The mapping class 
group $\cM_{g,n,k}$ of $S$ is the group
of the isotopy classes of orientation preserving homeomorphisms of $S$ 
which keep the boundary of $S$ and all the distinguished points
 pointwise fixed. In this paper we study
the problem of finding a finite presentation for $\cM_{g,n}=\cM_{g,n,0}$. 
We restrict our 
attention to the case $n=1$. The case $n=0$ is easily obtained from
the case $n=1$. In principle the case of $n>1$ (and the case of $k>0$)
can be obtained from the case $n=1$ via standard exact sequences, but
this method does not produce a global formula for the case of several
 boundary components 
 and the presentation (in contrast to the ones we shall describe 
 for the case $n=0$ and $n=1$) becomes rather ugly. On the other hand 
 Gervais in \cite{Gervais} succeeded recently to produce a finite 
 presentation of $\cM_{g,n}$ starting from the results in \cite{Wajnryb}
 and using a new approach.

A presentation for $g=1$ has been known for a long time. A quite simple
presentation for $g=2$ was established in 1973 in \cite{Birman-Hilden},
but the method did not generalize to higher genus. 
In 1975 McCool proved in \cite{McCool}, by purely algebraic methods, that
 $\cM_{g,1}$ is finitely presented for any 
genus $g$. It seems that extracting an explicit finite presentation from
his proof is very difficult. In 1980 appeared the groundbreaking paper
of Hatcher and Thurston
\cite{Hatcher-Thurston} in which they gave an algorithm for 
constructing a finite
 presentation
for the group $\cM_{g,1}$ for an arbitrary $g$. In
1981 Harer applied their algorithm in \cite{Harer} to obtain a finite 
(but very unwieldy)
explicit presentation of $\cM_{g,1}$. His presentation was simplified by
Wajnryb in 1983 in \cite{Wajnryb}. A subsequent Errata \cite{Errata}
corrected small errors in the latter. The importance of the full circle of
 ideas in these papers can be jugded from a small sample of subsequent work 
 which relied on the presentation in \cite{Wajnryb}, eg \cite{Kohno},
 \cite{Lu}, \cite{Matsumoto}, \cite{Matveev-Polyak}.
The proof of Hatcher and Thurston was deeply original, and solved an 
outstanding
open problem using novel techniques. These included arguments based upon 
Morse and Cerf Theory, as presented by Cerf in \cite{Cerf}.\par

In this paper we shall give, in one place, a complete hands-on proof of
a simple presentation for the groups $\cM_{g,0}$ and $\cM_{g,1}$.
Our approach will follow the lines set in \cite{Hatcher-Thurston}, but
 we will be 
able to use elementary methods in the proof of the connectivity and
simple connectivity of the cut system complex. In particular, our work 
does not rely on Cerf theory. At the same time we will gather all of
the computational details in one place, making the result accesible for 
independent checks. Our work yields a slightly different set of generators
and relations from the ones used in \cite{Wajnryb} and in \cite{Errata}.
The new presentation makes the computations in section 4  a little simpler.
We shall give both presentations and prove that they are equivalent. 
 \par

A consequence of this paper is that, using Lu \cite{Lu} or Matveev and
Polyak \cite{Matveev-Polyak}, the fundamental theorem of the Kirby
calculus \cite{Kirby} now has a completely elementary proof (ie, one
which makes no appeal to Cerf theory or high dimensional arguments).

 This paper is organised as follows. In section 2 we give a new proof
 of the main theorem of Hatcher and Thurston. In section 3 we derive a
 presentation of the mapping class group following (and explaining)
 the procedure described in \cite{Hatcher-Thurston} and in
 \cite{Harer}.  In section 4 we reduce the presentation to the simple
 form of Theorem \ref{simple presentation} repeating the argument from
 \cite{Wajnryb} with changes required by a slightly different setup.
In section 5 we deduce the case of a closed surface and in section 6 we 
translate the presentation into the form given in [20], see Remark 1
below.

We start with a definition of a basic element of the mapping
 class group.
\begin{definition}\label{Dehn twist} {\rm A (positive) Dehn twist 
with respect to a simple closed  curve
$\alpha$  on an oriented surface $S$ is an isotopy class of
 a homeomorphism $h$ of $S$, supported in a regular neighbourhood $N$ of
 $\alpha$ (an annulus), obtained as follows: we cut $S$ open along $\alpha$,
 we rotate
one side of the cut by 360 degrees to the right 
(this makes sense on an oriented surface)
and then glue the surface back together, damping out the rotation to the
identity at the boundary of $N$.
The Dehn twist (or simply twist) with respect to $\alpha$ will 
be denoted by $T_\alpha$. If curves $\alpha$ and $\beta$ intersect only 
at one point and are transverse then $T_\alpha(\beta)$, up to an isotopy,
 is obtained from the 
union $\alpha\cup\beta$ by splitting the union at the intersection point.}
\end{definition}
\begin{figure}[ht!]
\epsfxsize1.05\hsize\cl{\hglue -0.5cm
\relabelbox\small
\epsffile{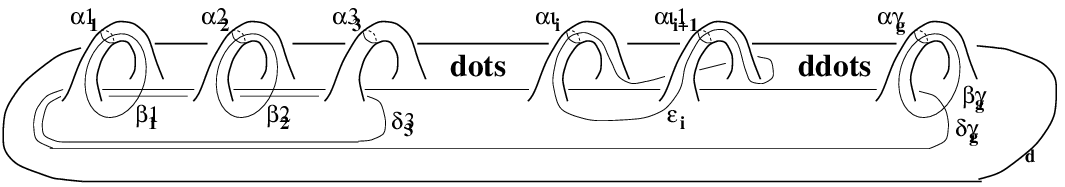}
\relabel{a1}{$\alpha_1$}
\relabel{a2}{$\alpha_2$}
\relabel{a3}{$\alpha_3$}
\relabel{ai}{$\alpha_i$}
\relabela<-1pt,1pt> {ai1}{$\alpha_{i+1}$}
\relabel{ag}{$\alpha_g$}
\relabel{b1}{$\beta_1$}
\relabel{b2}{$\beta_2$}
\relabela<-1pt,0pt> {bg}{$\beta_g$}
\relabel{d3}{$\delta_3$}
\relabela<-1pt,-3pt> {dg}{$\delta_g$}
\relabela<2pt,-2pt> {d}{$\partial$}
\relabel{e}{$\epsilon_i$}
\relabel{dots}{$\ldots$}
\relabel{ddots}{$\ldots$}
\endrelabelbox}
\caption{Surface $S_{g,1}$}
\label{general surface}
\end{figure}
\par
We shall say that two elements $a,b$ of a group are braided (or satisfy
the braid relation) if $aba=bab$.

\begin{theorem} \label{simple presentation} Let $S_{g,1}$ be a compact, 
orientable surface of genus $g\geq 3$ with one boundary component. 
Let $a_i,b_i,e_i$ denote Dehn twists
 about the curves
$\alpha_i,\beta_i,\epsilon_i$ on $S_{g,1}$ depicted on  Figure
\ref{general surface}.
The mapping class group 
$\cM_{g,1}$ of $S_{g,1}$ is generated by  elements
$b_2,b_1,a_1,e_1,a_2,e_2,
\dots,a_{g-1},e_{g-1},a_g$ and has defining relations:
\par\noindent
{\rm (M1)}\qua Elements $b_2$ and $a_2$ are braided and $b_2$ commutes with $b_1$.
Every other pair of consecutive elements (in the above order) is braided
and every other pair of non-consecutive elements commute.
\par\noindent
{\rm (M2)} \ \ \ $(b_1a_1e_1a_2)^5=b_2a_2e_1a_1b_1^2a_1e_1a_2b_2$
\par\noindent
{\rm (M3)} \ \ \ $d_3a_1a_2a_3=d_{1,2}d_{1,3}d_{2,3}$, \ where\par
$d_{1,2}=(a_2e_1a_1b_1)^{-1}b_2(a_2e_1a_1b_1)$, \ 
$d_{1,3}=t_2d_{1,2}t_2^{-1}$, \ 
$d_{2,3}=t_1d_{1,3}t_1^{-1}$, \nl
$\strut t_1=e_1a_1a_2e_1$, \  $t_2=e_2a_2a_3e_2$, \ 
$d_3 =b_2a_2e_1b_1^{-1}d_{1,3}b_1e_1^{-1}a_2^{-1}b_2^{-1}$.
\end{theorem}

\begin{theorem}\label{M2 presentation} The mapping class group 
$\cM_{2,1}$ of an orientable surface $S_{2,1}$
of genus $g=2$ with one boundary component is generated by  elements
$b_2,b_1,a_1,$ $e_1,a_2$ and has defining relations {\rm (M1)} and {\rm (M2)}.
\end{theorem}

\begin{theorem}\label{presentation closed}
The mapping class group $\cM_{g,0}$ of a compact, closed,
orientable surface of genus $g>1$ can be obtained from the above 
presentation of $\cM_{g,1}$ by adding one more relation:
\par\noindent
{\rm (M4)} \ \ \ $[b_1a_1e_1a_2\dots a_{g-1}e_{g-1}a_ga_ge_{g-1}a_{g-1}
\dots e_1a_1b_1,d_g]=1$, where
\par $t_i=e_ia_ia_{i+1}e_i$, for $i=1,2,\dots,g-1$, \
$d_2=d_{1,2}$, \nl
 $\strut d_i=(b_2a_2e_1b_1^{-1}t_2t_3\dots t_{i-1})
d_{i-1}(b_2a_2e_1b_1^{-1}t_2t_3\dots t_{i-1})^{-1}$, for $i=3,4,\dots,g$
\end{theorem} 

The presentations in Theorems 1 and 3 are equivalent to but not the same
as
those in \cite{Wajnryb} and \cite{Errata}. We now 
give alternative presentations of Theorems 1 and 3, with the goal of
correlating the work in this paper with that in \cite{Wajnryb} and
\cite{Errata}. See Remark 1, below, for a dictionary which
allows one to move between Theorems 1$^\prime$ and 3$^\prime$ and the
results in
\cite{Wajnryb} and \cite{Errata}. See Section 6 of this paper for 
a proof that the
presentations in Theorems 1 and $1^\prime$, and also Theorems 3 and
$3^\prime$,
are equivalent.
\par\medskip

\noindent
{\bf Theorem 1$^\prime$}\qua {\sl The mapping class group $\cM_{g,1}$ admits 
a presentation with generators $b_2,b_1,a_1,e_1,a_2,e_2,
\dots,a_{g-1},e_{g-1},a_g$ and with defining relations:
\par\noindent
{\rm (A)}\qua Elements $b_2$ and $a_2$ are braided and $b_2$ 
commutes with $b_1$.
Every other pair of consecutive elements (in the above order) is braided
and every other pair of non-consecutive elements commute.
\par\noindent
{\rm (B)} \ \ \ $(b_1a_1e_1)^4=
(a_2e_1a_1b_1^2a_1e_1a_2)b_2(a_2e_1a_1b_1^2a_1e_1a_2)^{-1}b_2$
\par\noindent
{\rm (C)} \ \ \ $e_2e_1b_1\tilde d_3=
\tilde t_1^{-1}\tilde t_2^{-1}b_2\tilde t_2\tilde t_1\tilde t_2^{-1}
b_2\tilde t_2b_2$
 \ where\par
$\tilde t_1=a_1b_1e_1a_1$,\ \  $\tilde t_2=a_2e_1e_2a_2$, \ \
$\tilde d_3 =(a_3e_2a_2e_1a_1u)v(a_3e_2a_2e_1a_1u)^{-1}$,\nl
$\strut u=e_2^{-1}a_3^{-1}\tilde t_2^{-1}b_2\tilde t_2a_3e_2$,\ \ 
$v=(a_2e_1a_1b_1)^{-1}b_2(a_2e_1a_1b_1)$.}

\par\medskip\noindent
{\bf Theorem 3$^\prime$}\qua {\sl The mapping class group 
 $\cM_{g,0}$ of a compact, closed,
orientable surface of genus $g>1$ can be obtained from the above 
presentation of $\cM_{g,1}$ by adding one more relation:
\par\noindent
{\rm (D)} \ \ \ $[a_ge_{g-1}a_{g-1}\dots e_1a_1b_1^2a_1e_1\dots 
a_{g-1}e_{g-1}a_g,\tilde d_g]=1$, where
\par 
$\tilde d_g=u_{g-1}u_{g-2}\dots u_1b_1(u_{g-1}u_{g-2}\dots u_1)^{-1}$,
\ \ $\strut u_1=(b_1a_1e_1a_2)^{-1}v_1a_2e_1a_1$,
\nl
$\strut u_i=(e_{i-1}a_ie_ia_{i+1})^{-1}v_ia_{i+1}e_ia_i$ \ \ 
for \ \ $i=2,\dots,g-1$, 
\nl
$\strut v_1=(a_2e_1a_1b_1^2a_1e_1a_2)b_2(a_2e_1a_1b_1^2a_1e_1a_2)^{-1}$,
\nl
$\strut v_i=\tilde t_{i-1}^{-1}\tilde t_i^{-1}
v_{i-1}\tilde t_i\tilde t_{i-1}$ \ \ for \ \ $i=2,\dots,g-1$,
\nl 
$\strut \tilde t_1=a_1e_1b_1a_1$, 
\ \ $\strut \tilde t_i=a_ie_ie_{i-1}a_i$ \ \ for \ \ $i=2,\dots ,g-1$.}
\par

\begin{remark} {\rm We now explain how to move back and forth between
the results in this paper and those in \cite{Wajnryb} and \cite{Errata}.
\begin{enumerate}
\item [(i)] The surface and the curves in \cite{Wajnryb} look different
from the surface and the curves on Figure \ref{general surface}. However
if
we compare the Dehn twist generators in Theorem $1^\prime$
 with those in
Theorem 1 of \cite{Wajnryb} and \cite{Errata} we see that corresponding
curves have the same 
intersection pattern. Thus there exists a homeomorphism of one surface
onto the 
other which takes the curves of one family onto the corresponding 
curves of the other family. The precise correspondence is given by:
\par\smallskip\noindent
\hbox{}\hspace{3.5cm}$(b_2,b_1,a_1,e_1,a_2,e_2,
\dots,a_{g-1},e_{g-1},a_g)$\nl 
\hbox{}\hspace{2.7cm}$\longleftrightarrow\ \ 
(d,a_1,b_1,a_2,b_2,a_3,\dots,b_{n-1},a_n,b_n)$,
\par\smallskip
where the top sequence refers
to Dehn twists about the curves in Figure 1 of this paper and the bottom
sequence
refers to Dehn twists about the curves in Figure 1 on page 158 of 
\cite{Wajnryb}.  
\item [(ii)] Composition of homeomorphisms 
in \cite{Wajnryb} was
performed from left to right, while in the present paper we use the 
standard composition from right to left. 
\item [(iii)] The element $d_g$ in this paper represents a Dehn twist
about the
curve $\delta_g$ in Figure \ref{general surface} of this paper. The
element
${\tilde d}_g$ in relation (D) of Theorem $3^\prime$ represents a Dehn
twist
about the curve $\beta_g$ in Figure 1.  We wrote
$d_g$ as a particular product of the generators in $\cM_{g,1}$.  
It follows from the argument in
the  last
section that any other such product representing $d_g$
 will also do. 
\item [(iv)] In the case of genus $g=2$ we should
omit relation (C) in Theorems $1^\prime$ and $3^\prime$.
\end{enumerate}
}\end{remark}
\par

Here is the plan of the proof of the theorems.
Following Hatcher and Thurston we  define a 2--dimensional cell
complex $X$ on which the mapping class group $\cM_{g,1}$ acts by cellular
  transformations
and the action is transitive on the vertices of $X$. We give a new
elementary proof of the fact that $X$ is connected and simply connected.
We then describe the stabilizer $H$ of one vertex of $X$ under the action
of $\cM_{g,1}$ and we determine an explicit presentation of $H$. Following the
algorithm of Hatcher and Thurston we get from it a presentation of $\cM_{g,1}$.
Finally we reduce the presentation to the form in Theorem
\ref{simple presentation} and as a corollary get Theorem
\ref{presentation closed}.\par\medskip
{\bf Acknowledgements}\qua I wish to thank very much the referee, who studied 
the paper very carefully and made important suggestions to improve it.

This research was partially supported by the Fund for the Promotion of
Research at the Technion and the Dent Charitable Trust --- a non
military research fund.

\section{Cut-system complex}\label{cut-system complex}

We denote by $S$ a compact, connected oriented surface of  genus $g> 0$
 with $n\ge 0$ boundary components. We denote by $\bar S$ a closed surface
 obtained from $S$ by capping each boundary component with a disk.
By a {\it curve} we shall mean a simple closed curve on $S$. We are mainly 
 interested in the isotopy classes of curves on $S$. The main goal in the 
proofs will be to decrease the number of intersection points between 
different curves. If the Euler characteristic of $S$ is negative we can put
a hyperbolic metric on $S$ for which the boundary curves are geodesics.
 Then the isotopy class of any
 non-separating curve on $S$ contains
 a unique simple closed geodesic, which is the shortest curve in its 
isotopy class. If we replace each non-separating curve
 by the unique 
geodesic isotopic to it we shall minimize the number of intersection points
between every two non-isotopic curves, by Corollary \ref{hass-scott}.
 So we can think of curves as
geodesics. In the proof we may construct new curves, which are not
geodesics and which have small intersection number with some other curves.
When we replace the new curve by the corresponding geodesic we further
 decrease the 
intersection number. If $S$ is a closed torus we can choose a flat metric
on $S$. Now geodesics are not unique but still any choice of geodesics
will minimize the intersection number of any pair of non-isotopic curves.
If two curves are isotopic on $S$ then they correspond to the same 
geodesic, but we shall call them disjoint because we can isotop one off
the other. Geodesics are never tangent.
\par      
If $\alpha$ and $\beta$ are curves then $|\alpha\cap\beta|$ denotes their
geometric intersection number, ie, the number of intersection points
 of $\alpha$ and $\beta$. If $\alpha$ is a curve we denote by
  $[\alpha]$  the homology class represented by $\alpha$ in
 $H_1(\bar S,\ints)$, up to a sign. We denote by
$i(\alpha,\beta)$  the absolute value of the algebraic intersection
  number of  $\alpha$ and $\beta$. It depends only on the classes
  $[\alpha]$ and $[\beta]$.
\par
\begin{figure}[ht!]\epsfxsize\hsize\cl{
\relabelbox\small
\epsffile{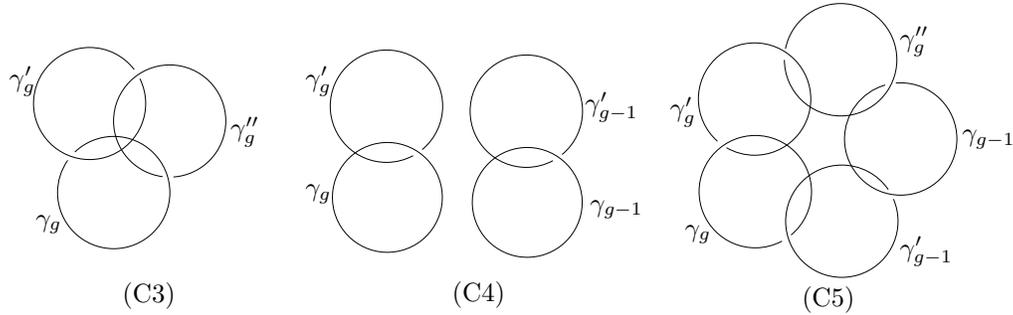}
\relabela <-1.5pt,0pt> {g31}{$\gamma_g$}
\relabela <-1.5pt,0pt> {g4}{$\gamma_g$}
\relabela <-1.5pt,0pt> {g51}{$\gamma_g$}
\relabela <-2pt,0pt> {g3}{$\gamma'_g$}
\relabela <-1.5pt,0pt> {g41}{$\gamma'_g$}
\relabela <-1.5pt,0pt> {g5}{$\gamma'_g$}
\relabela <-1.5pt,0pt> {g32}{$\gamma''_g$}
\relabela <-1.5pt,0pt> {g5m1}{$\gamma''_g$}
\relabela <-1.5pt,0pt> {g4r1}{$\gamma'_{g-1}$}
\relabela <-1.5pt,0pt> {g52}{$\gamma'_{g-1}$}
\relabela <-1.5pt,0pt> {g4r}{$\gamma_{g-1}$}
\relabela <-2.5pt,0pt> {g5m}{$\gamma_{g-1}$}
\relabela <0pt,-10pt> {C3}{(C3)}
\relabela <0pt,-4pt> {C4}{(C4)}
\relabel{C5}{(C5)}
\endrelabelbox}
\caption{Configurations of curves in paths (C3), (C4) and (C5)}
\label{faces}
\end{figure}

 We shall describe now the cut-system complex $X$ of $S$. 
 To construct $X$ we consider collections of $g$ disjoint curves
$\gamma_1,\gamma_2,\dots,\gamma_g$ in $S$ such that when we cut $S$
 open along these curves we get a connected surface
 (a sphere with $2g+n$ holes). An isotopy class of such
 a collection we call a cut system $\lan \gamma_1,\dots,\gamma_g\ran$. 
We can say that a cut system is
 a collection of geodesics. A curve is contained in a cut-system if it is
 one of the curves of the collection. If 
$\gamma_i^\prime$ is a curve in $S$, which meets $\gamma_i$ at one point
 and is
disjoint from other curves $\gamma_k$ of the cut system 
$\lan \gamma_1,\dots,\gamma_g\ran$, then
 $\lan \gamma_1,\dots,\gamma_{i-1},\gamma_i^\prime,
\gamma_{i+1},\dots,\gamma_g\ran$ forms another cut system. 
In such a situation 
the replacement  $\lan \gamma_1,\dots, \gamma_i,\dots,\gamma_g\ran\to
\lan \gamma_1,\dots,\gamma_i^\prime,\dots,\gamma_g\ran$
is called a simple move. For brewity we shall often drop the symbols for
unchanging circles and shall write $\lan
 \gamma_i\ran\to\lan \gamma_i^\prime\ran$. The cut systems on $S$ form
 the 0--skeleton (the vertices) of the complex $X$.  We join two vertices
 by an edge if and only if the corresponding cut systems are related
by a simple move.  We get the 1--skeleton $X^1$.  By a {\em path} we mean an
 edge-path in $X^1$. It consists of a 
sequence of vertices ${\bf p}=( v_1, v_2,\dots,v_k)$
where two consecutive vertices are related by a simple move. A path is 
closed if $ v_1=v_k$.
We distinguish three types of closed paths:\\
If three vertices (cut-systems) have $g-1$ curves 
$\gamma_1,\dots,\gamma_{g-1}$ in common and if the remaining three curves 
$\gamma_g,\gamma_g^\prime,\gamma_g^{\pp}$ intersect each other once,
as on Figure \ref{faces}, C3, then the vertices form a closed 
triangular path:\par
(C3)\qua a triangle\qua $\lan \gamma_g\ran\to\lan\gamma_g^\prime\ran\to\lan
 \gamma_g^{\pp}\ran\to\lan \gamma_g\ran$.\\
 If four vertices have $g-2$ curves $\gamma_1,\dots,\gamma_{g-2}$ in
 common and the remaining pairs of curves consist of 
 $(\gamma_{g-1},\gamma_g)$, $(\gamma_{g-1}^\prime,\gamma_g)$,
 $(\gamma_{g-1}^\prime,\gamma_g^\prime)$, $(\gamma_{g-1}, \gamma_g^\prime)$
 where the curves intersect as on Figure \ref{faces}, C4, then the vertices
 form a closed square path:\par
(C4)\qua a square\qua $\lan \gamma_{g-1},\gamma_g\ran\to\lan\gamma_{g-1}^\prime,
\gamma_g\ran
\to\lan\gamma_{g-1}^\prime,\gamma_g^\prime\ran
\to\lan \gamma_{g-1},\gamma_g^\prime\ran\to\lan \gamma_{g-1},\gamma_g\ran$.\\
If five vertices have $g-2$ curves $\gamma_1,\dots,\gamma_{g-2}$ in
 common and the remaining pairs of curves consist of 
 $(\gamma_{g-1},\gamma_g)$, $(\gamma_{g-1},\gamma_g^\prime)$,
 $(\gamma_{g-1}^\prime,\gamma_g^\prime)$, $(\gamma_{g-1}^\prime,
  \gamma_g^{\pp})$ and $(\gamma_g,\gamma_g^{\pp})$
 where the curves intersect as on Figure \ref{faces}, C5, then the vertices
 form a closed pentagon path:\par
(C5)\qua a pentagon\qua $\lan \gamma_{g-1},\gamma_g\ran\to\lan \gamma_{g-1},
\gamma_g^\prime\ran\to\lan\gamma_{g-1}^\prime,\gamma_g^\prime\ran
\to\lan\gamma_{g-1}^{\prime},\gamma_g^{\pp}\ran\to$\nl
$\strut\lan\gamma_g^{\pp},
\gamma_g\ran
\to\lan \gamma_{g-1},\gamma_g\ran$.\\
\par
$X$ is a 2--dimensional cell complex obtained from $X^1$ by attaching
 a 2--cell to
 every closed edge-path of type (C3), (C4) or (C5).
 The mapping class group of $S$ acts on $S$ by homeomorphisms so 
 it takes cut systems to
  cut systems. Since the edges and the faces of $X$ are determined by 
  the intersections of pairs of curves, which are clearly preserved by
  homeomorphisms, the action on $X^0$ extends to a cellular
  action on $X$.\par
 
 In this section we shall prove the main result of \cite{Hatcher-Thurston}:
   
 \begin{theorem}\label{Hatcher-Thurston} $X$ is connected and simply
   connected.\end{theorem}
 
We want to prove that every closed path $\bf p$ is null-homotopic. 
If $\bf p$ is null-homotopic we shall write ${\bf p\sim o}$.
 We start with a closed path
${\bf p}=(v_1,\dots,v_k)$ and try to simplify it. If $\bf q$ is a short-cut,
an edge path connecting a vertex $v_i$ of $\bf p$ with $v_j$, $j>i$,
 we can split
$\bf p$ into two closed edge-paths:\nl ${\bf p}_1=(v_j,v_{j+1},\dots,
v_k,v_2,\dots,v_{i-1},{\bf q})$ and ${\bf p}_2=({\bf q}^{-1},v_{i+1},v_{i+2},
\dots,v_j)$.\nl If both paths are null-homotopic in $X$ then $\bf p\sim o$.
We want to prove Theorem \ref{Hatcher-Thurston} by splitting path $\bf p$
into simpler paths according to a notion of complexity which is described 
in the next definition.

\begin{definition}\label{radius} {\rm Let ${\bf p}=(v_1,v_2,\dots,v_k)$ be a path
in $X$. Let $\alpha$ be a fixed curve of some fixed  $v_j$. We define
 {\em distance} from $\alpha$ to a vertex $v_i$ to be 
$d(\alpha,v_i)=min\{|\alpha\cap\beta|:\beta\in v_i\}$. The {\em radius of
${\bf p}$ around $\alpha$} is equal to the maximum distance from $\alpha$
to the vertices of ${\bf p}$. The path $\bf p$ is called a {\em segment}
if every vertex of $\bf p$ contains a fixed curve $\alpha$. We shall write
$\alpha$--segment if we want to stress the fact that $\alpha$ is the common
curve of the segment. If the segment has several fixed curves we can write
$(\alpha,\beta,\dots,\gamma)$--segment.}\end{definition}
\par 
 Theorem \ref{Hatcher-Thurston} will be proven by induction on the 
 genus of $S$, for the given genus it will be proven by induction on 
 the radius of a path $\bf p$ and for a given radius $m$ around a curve
 $\alpha$ we shall induct on the number of segments of $\bf p$ which have 
 a common curve $\gamma$ with $|\gamma\cap\alpha|=m$.
 The main tool in the proof of Theorem \ref{Hatcher-Thurston}
  is a reduction of the number of
intersection points of curves. 

\begin{figure}[ht!]
\relabelbox\small
\epsffile{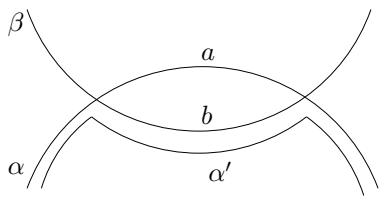}
\relabel{a}{$a$}
\relabel{b}{$b$}
\relabel{al}{$\alpha$}
\relabel{a1}{$\alpha'$}
\relabel{be}{$\beta$}
\endrelabelbox
\caption{Curves $\alpha$ and $\beta$ form a 2--gon.}
\label{2-gon}
\end{figure}

\begin{definition} \label{def 2-gon} {\rm Curves $\alpha$ and $\beta$ on $S$
have an {\em excess intersection} if there exist curves $\alpha^\prime$ and 
$\beta^\prime$, isotopic to $\alpha$ and $\beta$ respectively, and such that
$|\alpha\cap\beta|>|\alpha^\prime\cap\beta^\prime|$. Curves 
$\alpha$ and $\beta$ form a 2--{\em gon} if there are arcs $a$ of $\alpha$
and $b$ of $\beta$, which meet only at their common end points and do not 
meet other points of $\alpha$ or $\beta$ and such that $a\cup b$ bound a 
disk, possibly with holes, on $S$. The disk is called a 2--gon. We can 
{\em cut off} the 2--gon from $\alpha$ by replacing the arc $a$ of $\alpha$
by the arc $b$ of $\beta$. We get a new curve $\alpha^\prime$ 
(see Figure \ref{2-gon}).}\end{definition}
\begin{lemma}[Hass--Scott, see \cite{Hass-Scott}, Lemma 3.1]
 \label{hass}
If $\alpha,\beta$ are two curves on a surface $S$ having an excess
 intersection
 then they form a 2--gon (without holes) on $S$.
\end{lemma}
\begin{cor}\label{hass-scott} Two simple closed geodesics on $S$ have no
 excess intersection. In particular if we replace two curves by 
 geodesics in their isotopy class then the number of intersection
 points between the curves does not increase.\end{cor}

\proof If there is a 2--gon we can shorten one geodesic in its homotopy
class, by first
cutting off the 2--gon and then smoothing corners.\endpf 

\begin{lemma}\label{cut off 2-gon} Consider a finite collection of 
simple closed geodesics
on $S$. Suppose that curves $\alpha$ and $\beta$ of the collection form
a minimal 2--gon (which does not contain another 2--gon). Let $\alpha^\prime$
be the curve obtained by cutting off the minimal 2--gon from $\alpha$ and 
passing to the isotopic geodesic. Then $|\alpha\cap\alpha^\prime|=0$,
$|\beta\cap\alpha^\prime|<|\beta\cap\alpha|$ and $|\gamma\cap\alpha^\prime|
\leq |\gamma\cap\alpha|$ for any other curve $\gamma$ of the collection.
In particular if $ |\gamma\cap\alpha|=1$ then $ |\gamma\cap\alpha^\prime|=1$.
Also $[\alpha]=[\alpha^\prime]$.\end{lemma}
\proof Since the 2--gon formed by arcs $a$ and $b$ of 
$\alpha$ and $\beta$ is minimal every other curve $\gamma$ of the 
collection intersects the 2--gon along arcs which meet $a$ and $b$ once.
Thus cutting off the 2--gon will not change $|\alpha\cap\gamma|$ 
and it may only decrease after passing to the isotopic
geodesic. Clearly $\alpha$ and $\alpha^\prime$ are disjoint and
homologous on $\bar S$ and by passing to $\alpha^\prime$ we remove at
 least two intersection points of $\alpha$ with $\beta$.\endpf

\subsection{The case of genus 1}

In this section we shall assume that $S$ is a surface of genus one,
 possibly with boundary.  By $\bar S$ we shall denote the closed 
 torus obtained by glueing a disk
to each boundary component of $S$. We want to prove: 
\begin{proposition} \label{genus1} If $S$ has genus one then the
 cut system complex $X$ of $S$ is connected and simply connected.
\end{proposition} 

On a closed torus $\bar S$ the homology class
 $[\alpha]$ of a curve $\alpha$ is defined by a pair
of relatively prime integers, up to a sign, after a fixed choice of
 a basis of $H_1(\bar S,\ints)$. If $[\alpha]=(a_1,a_2)$ and 
 $[\beta]=(b_1,b_2)$ then the absolute value of their algebraic
 intersection number is equal $i(\alpha,\beta)=|a_1b_2-a_2b_1|$. If
 $\alpha$ and $\beta$ are geodesics on $\bar S$ then $|\alpha\cap\beta|=
 i(\alpha,\beta)$, therefore it is also true for curves on $S$ which form
 no 2--gons, because then they have no excess intersection on $\bar S$.
 
\begin{figure}[ht!]
\relabelbox\small
\noindent\mbox{}\hglue.7in\epsfxsize 3in\epsffile{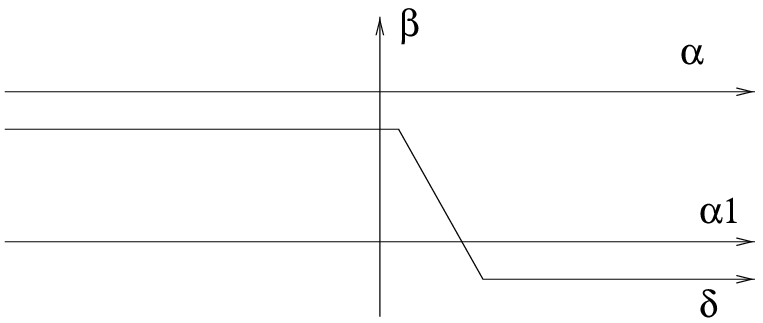}
\relabel{d}{$\delta$}
\relabel{a}{$\alpha$}
\relabel{a1}{$\alpha$}
\relabel{b}{$\beta$}
\endrelabelbox
\caption{ Curve $\delta$ has smaller intersection with $\alpha$ and $\beta$.}
\label{equal signs}
\end{figure}

\begin{lemma} Let $\alpha$, $\beta$ be nonseparating curves on $S$ and 
suppose that $|\alpha\cap\beta|=k\neq 1$. Then there exists a nonseparating 
curve $\delta$ such that if $k=0$ then $|\delta\cap\alpha|=|\delta\cap
\beta|=1$ and if $k>1$ then $|\delta\cap\alpha|<k$ and
$|\delta\cap\beta|<k$.\end{lemma}
\proof If $|\alpha\cap\beta|=0$ then the curves are isotopic on $\bar S$
and they split $S$ into two connected components $S_1$ and $S_2$. We
can choose points $P$ and $Q$ on $\alpha$ and $\beta$ respectively and
connect them by simple arcs in $S_1$ and in $S_2$. The union of the arcs 
forms the required curve $\delta$. If $k>1$ and if the curves have
an excess intersection on $\bar S$ then they form a 2--gon on $S$. We can
cut off the 2--gon decreasing the intersection, by Lemma \ref{cut off 2-gon}.
If there are no 2--gons then the algebraic and the geometric intersection 
numbers are equal. In particular all intersections have the same sign.
Consider two intersection points consecutive along $\beta$.
 Choose $\delta$ as on  Figure \ref{equal signs}. Then $|\delta\cap\alpha|=1$
 so it is nonseparating and $|\delta\cap\beta|<k$.\endpf
 
\begin{cor} If surface $S$ has genus 1 then the cut system complex of $S$
is connected.\end{cor}
\proof A cut system on $S$ is an isotopy class of a single curve. If two
curves intersect once they are connected by an edge. It follows 
from the last lemma by induction that any two curves can be connected by an 
edge-path in $X$.\endpf
We now pass to a proof that closed paths are null-homotopic.
\begin{lemma}\label{square} A closed path
 ${\bf p}=(\delta_1,\delta_2,\delta_3,\delta_4,\delta_1)$ such that
 $|\delta_2\cap \delta_4|=0$
is null-homotopic in $X$. \end{lemma}

\proof Let $\beta=T_{\delta_2}^{\pm 1}(\delta_3)$ be the image of $\delta_3$
by the Dehn Twist along $\delta_2$. Recall that as a set
 $\beta=\delta_2\cup \delta_3$ outside a small neighbourhood of
 $\delta_2\cap\delta_3$.
From this we get that $|\beta\cap \delta_2|=1$, $|\beta\cap \delta_3|=1$
 and $|\beta\cap \delta_4|=1$.
 Thus
 ${\bf p^\prime}=(\delta_1,\delta_2,\beta,\delta_4,\delta_1)$ is a closed
  path 
which is homotopic to ${\bf p}$ because ${\bf p-p}^\prime$ splits into 
a sum of two triangles  ${\bf t_1}=(\beta,\delta_2,\delta_3,\beta)$ 
and ${\bf t_2}=(\beta,\delta_3,\delta_4,\beta)$. We also have
$i(\beta,\delta_1)=|i(\delta_3,\delta_1)\pm i(\delta_2,\delta_1)|$, so for
 a suitable
choice of the sign of the twist we have
 $i(\delta_3,\delta_1)>i(\beta,\delta_1)$
 unless $i(\delta_3,\delta_1)=0$.
We may assume by induction 
that $i(\delta_3,\delta_1)=0$. If $|\delta_1\cap\delta_3|>0$
then $\delta_1$ and $\delta_3$  have an excess intersection on $\bar S$
and form a 2--gon on $S$. We can cut off the 2--gon from $\delta_3$ getting
a new curve $\beta$ such that
 $[\beta]=[\delta_3]$,
 $|\beta\cap \delta_3|=0$, $|\beta\cap \delta_1|<|\delta_3\cap \delta_1|$ 
 and $|\beta\cap\delta_i|=|\delta_3\cap\delta_i|$ for $i\neq 1$, by
 Lemma \ref{cut off 2-gon}.
 We get a new closed path
  ${\bf p^\prime}=(\delta_1,\delta_2,\beta,\delta_4,\delta_1)$
 and the difference 
between it and the old path is a closed path
 ${\bf q}=(\beta,\delta_2,\delta_3,\delta_4,\beta)$ 
with 
$|\beta\cap \delta_3|=0$ and $|\delta_2\cap \delta_4|=0$.
So by induction it suffices to
assume that $|\delta_1\cap \delta_3|=0$. If we now let
 $\beta=T_{\delta_2}(\delta_3)$ then $\beta$
intersects each of the four curves once so our path is a sum of four
 triangles and thus is null-homotopic in $X$.
\endpf

\begin{lemma} \label{no 2-gons} If
 ${\bf p}=(\alpha_1,\dots,\alpha_k)$ is a closed path then  there
exists a closed path
 ${\bf p^\prime}=(\alpha_1^\prime,\dots,\alpha_k^\prime)$
 such that  ${\bf p}^\prime$ is homotopic 
 to $\bf p$ in $X$, $[\alpha_i]=[\alpha_i^\prime]$ for all $i$ and 
 the collection of curves $\alpha_1^\prime,\dots,\alpha_{k-1}^\prime$
 forms no 2--gons.
\end{lemma}
\proof Suppose that there exists a 2--gon
bounded by arcs of two curves in  ${\bf p}$. Then there also exists 
a minimal 2--gon bounded by arcs of curves $\alpha_i$
and $\alpha_j$. If we cut off the 2--gon from $\alpha_i$ we get a curve
$\alpha_i^\prime$ such that $[\alpha_i]=[\alpha_i^\prime]$,
  $\alpha_i^\prime$ is disjoint from  $\alpha_i$, 
   $|\alpha_j\cap\alpha_i|>|\alpha_j\cap\alpha_i^\prime|$ and 
   $|\alpha_m\cap\alpha_i^\prime|\leq|\alpha_m\cap\alpha_i|$
   for any other curve $\alpha_m$, by Lemma \ref{cut off 2-gon}.
It follows that if we replace $\alpha_i$ by $\alpha_i^\prime$ 
 we get a new closed path ${\bf p^\prime}$ with
   a smaller number of
   intersection points between its curves. The difference of the two
   paths is a closed path ${\bf q}=(\alpha_{i-1},\alpha_i,\alpha_{i+1}, 
   \alpha_i^\prime,\alpha_{i-1})$ which is null-homotopic by Lemma
    \ref{square}. Thus
   ${\bf p^\prime}$ is homotopic to ${\bf p}$ in $X$. Lemma \ref{no 2-gons}
follows by induction on the total number of intersection points
between pairs of curves of ${\bf p}$.\endpf
\noindent
{\bf Proof of Proposition \ref{genus1}}\qua
  Let ${\bf p}=(\alpha_1,\dots,\alpha_k,\alpha_1)$ 
 be a closed path in $X^1$.  We may assume that the path has no 2--gons.
 By Lemma \ref{hass} there is no excess intersection on the closed torus
 $\bar S$. It means that the geometric intersection number
 of two curves of ${\bf p}$ is equal to their algebraic intersection
  number. 
 Let $m=max\{i(\alpha_1,\alpha_j)|j=1,\dots,k\}$ be the {\em radius} of
  ${\bf p}$ around $\alpha_1$. Suppose first that $m=1$.
 Two disjoint curves on $\bar S$ have the same homology class, and two
 curves representing the same class have algebraic intersection equal to 0.
 It follows that two consecutive curves in a path cannot be both disjoint
  from $\alpha_1$. If $k>4$ then either $|\alpha_1\cap\alpha_3|=1$ or
  $|\alpha_1\cap\alpha_4|=1$. We get an edge which splits ${\bf p}$ into
  two shorter closed paths with radius $1$. If $k=4$ and
  $|\alpha_1\cap\alpha_3|=0$ then ${\bf p\sim o}$ by Lemma \ref{square}.
  If $k=3$ then ${\bf p}$ is a triangle. Suppose now that $m>1$. We may
   assume by induction that every path of radius less than $m$ is
null-homotopic and every path of radius $m$ which has less
 curves $\alpha_j$ with $|\alpha_1\cap\alpha_j|=m$ is also null-homotopic.
Choose the smallest $i$ such that $i(\alpha_1,\alpha_i)=m$. Then
 $i(\alpha_1,\alpha_{i-1})< m$ and
$i(\alpha_1,\alpha_{i+1})\leq m$. Choose a basis of the homology group
 $H_1(\bar S)$ which
contains the curve $\alpha_1$. A homology class of a curve is then
 represented by 
a pair of integers $(a,b)$. We consider the homology classes and their
intersection numbers up to a sign. We have $[\alpha_1]=(1,0)$, 
$[\alpha_{i-1}]=(a,b)$,
$[\alpha_i]=(p,m)$ and $[\alpha_{i+1}]=(c,d)$. 
The intersection form is defined by $i((a,b),(c,d))=|ad-bc|$. Thus
$am-bp=\pm 1,cm-dp=\pm 1,
|b|< m, |d|\leq m$. 
We get $m(ad-bc)=(\pm d\pm b)$.
 Since 
$2m>|b|+|d|$ we must have $|ad-bc|=1$ or $ad-bc=0$. In the first case 
 $i(\alpha_{i-1},\alpha_{i+1})=|\alpha_{i-1}\cap \alpha_{i+1}|=1$. We
 can
\lq\lq cut off" the triangle
 ${\bf q}=(\alpha_{i-1},\alpha_i,\alpha_{i+1},\alpha_{i-1})$ getting
 a path which is null-homotopic by the induction hypothesis.
If $ad-bc=0$ then
 $i(\alpha_{i-1},\alpha_{i+1})=|\alpha_{i-1}\cap \alpha_{i+1}|=0$.
  Let $\beta=T_{\alpha_{i-1}}^{\pm 1}(\alpha_i)$.
Then $|\beta\cap \alpha_{i-1}|=1$
 and
$|\beta\cap \alpha_{i+1}|=1$.
We can replace $\alpha_i$ by $\beta$ getting a new closed path. 
Their difference is
the closed path
 ${\bf q}=(\alpha_{i-1},\alpha_i,\alpha_{i+1},\beta,\alpha_{i-1})$ which 
is null homotopic by
 lemma
 \ref{square}. Thus the new path is homotopic to the old path. For
 a suitable 
choice of the sign of the Dehn twist we have $i(\beta,\alpha_1)<m$.
It may happen that $|\beta\cap\alpha_1|\geq m$. We can get rid of 2--gons 
 by Lemma \ref{no 2-gons} and thus get rid of the excess intersection. We
get a homotopic path which is null-homotopic by the induction hypothesis.
 This concludes the proof of Proposition \ref{genus1} \endpf

\subsection{Paths of radius $0$}

From now until the end of section \ref{cut-system complex}, we assume
that $S$ is a surface of genus $g>1$ with a finite number of boundary
components. We denote by $X$ the cut-system complex of $S$. We assume:

\medskip\noindent
{\bf Induction Hypothesis 1}\qua The cut-system complex of a surface of genus
less than $g$ is connected and simply-connected.\par\medskip

We want to prove that every closed path in $X$ is null-homotopic.
We shall start with paths of radius zero. The simplest paths of radius 
zero are closed segments.

\begin{lemma} A closed segment is null-homotopic in $X$.\end{lemma}

\proof When we cut $S$ open along the common curve $\alpha$ the remaining
curves of each vertex form a vertex of a closed
path in the cut-system complex of a surface of a smaller genus.  By Induction
Hypothesis 1 
it is a sum of paths of type (C3), (C4) and (C5) there. When we
adjoin $\alpha$ to every vertex we get a splitting of the original
paths into null-homotopic paths.\endpf
In a similar way we prove:

\begin{lemma}\label{common curve} If two vertices of $X$ have
 one or more curves in common we can connect them by a path all of whose 
 vertices contain the common curves.\end{lemma}
 
 \proof If we cut
  $S$ open along the common curves the
 remaining collection of curves form two vertices of the cut-system 
 complex on the new surface of smaller genus. They can be connected by a path.
 If we adjoin all the common curves to each vertex of this path we get a path in $X$
 with the required properties.\endpf

 \begin{remark}\label{non-separating pair} {\rm 
 If $\alpha$ and $\beta$ are two
 disjoint non-separating curves on $S$ then $\alpha\cup\beta$ does not
 separate $S$ if and only if $[\alpha]\neq[\beta]$. In this case the pair
  $\alpha$,
 $\beta$ can be completed to a cut-system on $S$.}
\end{remark}
   
We shall now construct two simple types of null-homotopic paths
in $X$.

\begin{lemma}\label{hexagon relation} Let $\alpha_1$, $\alpha_2$, $\alpha_3$ 
be disjoint
curves such that the union of any two of them does not separate $S$ but
 the union of all three
separates $S$. Then there exist disjoint curves
 $\beta_1$, $\beta_2$, $\beta_3$ and a closed path\\
{\rm (C6)}\qquad\quad $\lan\alpha_1,\alpha_2\ran\to\lan\alpha_1,\beta_2\ran\to\lan\alpha_1,
\alpha_3\ran\to\lan\beta_1,\alpha_3\ran\to\lan\alpha_2,\alpha_3
\ran$\nl\noindent\hbox{}\hglue 9cm$\to\lan\alpha_2,\beta_3\ran\to\lan\alpha_2,\alpha_1\ran$,\\
which is null-homotopic in $X$.\end{lemma}

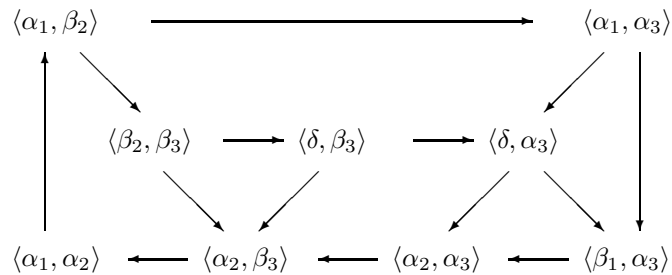
\begin{figure}[ht!]\unitlength0.9pt\small
\cl{\begin{picture}(330,140)
\put(30,20){$\lan\alpha_1,\alpha_2\ran$}
\put(105,23){\vector(-1,0) {25}}
\put(110,20){$\lan\alpha_2,\beta_3\ran$}
\put(185,23){\vector(-1,0) {25}}
\put(190,20){$\lan\alpha_2,\alpha_3\ran$}
\put(265,23){\vector(-1,0) {25}}
\put(270,20){$\lan\beta_1,\alpha_3\ran$}
\put(70,70){$\lan\beta_2,\beta_3\ran$}
\put(120,73){\vector(1,0) {25}}
\put(150,70){$\lan\delta,\beta_3\ran$}
\put(200,73){\vector(1,0) {25}}
\put(230,70){$\lan\delta,\alpha_3\ran$}
\put(30,120){$\lan\alpha_1,\beta_2\ran$}
\put(90,123){\vector(1,0) {160}}
\put(270,120){$\lan\alpha_1,\alpha_3\ran$}
\put(60,110){\vector(1,-1) {25}}
\put(280,110){\vector(-1,-1) {25}}
\put(295,110){\vector(0,-1) {75}}
\put(45,35){\vector(0,1) {75}}
\put(255,60){\vector(1,-1) {25}}
\put(240,60){\vector(-1,-1) {25}}
\put(160,60){\vector(-1,-1) {25}}
\put(95,60){\vector(1,-1) {25}}
\end{picture}} 
\caption{Hexagonal path}
\label{hexagon}
\end{figure}

\proof  Let $\gamma_3$,\dots, $\gamma_g$ be a cut system 
on a surface  $S-(\alpha_1\cup\alpha_2\cup\alpha_3)$ (not connected),
ie, a collection of curves which does not separate the surface 
any further. Then $\lan\alpha_1,\alpha_2,\gamma_3,\dots,\gamma_g\ran$
is a cut system on $S$. Let $S_1$ and $S_2$ be the components of a surface
 obtained by cutting $S$
open along all $\alpha_i$'s and $\gamma_j$'s. An arc connecting different
components of the boundary does not separate the surface so we can find
(consecutively) disjoint arcs  $b_1$ connecting $\alpha_1$ with $\alpha_2$,
$b_2$ connecting $\alpha_2$ with $\alpha_3$ and $b_3$ 
connecting $\alpha_3$ with $\alpha_1$ in $S_1$, and similar arcs 
$b_1^\prime$,
$b_2^\prime$ and $b_3^\prime$ in $S_2$ with the corresponding ends
coinciding in $S$. The pairs of corresponding arcs form the
required curves $\beta_1$, $\beta_2$ and $\beta_3$ and we get a closed
path ${\bf p}$ described in  (C6). Moreover the curves $\beta_2$ and 
$\beta_3$ are disjoint and $[\beta_2]\neq[\beta_3]$. To prove that 
the path is null-homotopic in $X$ we choose a curve
$\delta=T_{\alpha_2}(\beta_1)$.
 Then $|\delta\cap\alpha_1|=|\delta\cap\alpha_2|=|\delta\cap\beta_1|=
 |\delta\cap\beta_2|=1$ and $|\delta\cap\alpha_3|=|\delta\cap\beta_3|=0$.
 Figure \ref{hexagon} shows how ${\bf p}$ splits into a sum of triangles
 (C3), squares
 (C4) and pentagons (C5) and therefore is null-homotopic in $X$.\endpf

\begin{proposition} If $\bf p$ is a path of radius $0$
around a curve $\alpha$ then $\bf p\sim o$.
\end{proposition}
\proof  Let $v_0$ be a vertex of ${\bf p}$ containing $\alpha$.
 We shall prove the
  proposition by induction on the number of segments of $\bf p$ having
 a fixed curve disjoint from $\alpha$.
 Consider the maximal $\alpha$--segment of ${\bf p}$ which
contains the vertex $v_0$. We shall call
it the first segment. Let $v_1$ be the last vertex of the first segment.
The next vertex contains a curve $\beta$ disjoint from $\alpha$ such that
 $\beta$ is  the common curve of the next segment of ${\bf p}$.
 Since $|\alpha\cap\beta|=0$ the simple move from $v_1$ to the next vertex
 does not involve $\beta$ hence $v_1$ also contains $\beta$. Let $v_2$ be
the last vertex of the second segment.
If there are
only two segments then $v_2$  also
contains both $\alpha$ and $\beta$. By Lemma \ref{common curve} there is
an $(\alpha,\beta)$--segment connecting $v_1$ and $v_2$. Then ${\bf p}$ is
a sum of a closed $\alpha$--segment and a closed $\beta$--segment.
 So we may assume that there
is a third segment. The vertex $\tilde v$ of $\bf p$ following $v_2$ 
contains a curve $\gamma$ disjoint from $\alpha$ and $\gamma$ is the
common curve of the third segment. Let $v_3$ be the last vertex of the
third segment. We shall reduce the number of segments.
There are three cases.\par\medskip\noindent
{\bf Case 1}\qua Vertex  $v_2$ does not contain $\gamma$.\qua Since
 $\tilde v$ contains $\gamma$ and does not contain $\beta$ we have
  $|\beta\cap\gamma|=1$. Let $S_1$ be a surface of genus $g-1$ obtained 
  by cutting $S$ open 
along $\beta\cup\gamma$. Vertices $v_2$ and $\tilde v$ have $g-1$ curves 
in common 
and the common curves form a cut system $u$ on $S_1$. The union
$\beta\cup\gamma$ cannot separate $S-\alpha$ hence $\alpha$ does not
 separate $S_1$ and it belongs to a cut-system $u^\prime$ on $S_1$. 
 Vertices
 $u$ and $u^\prime$ can be connected by a path $\bf q$ in the cut-system 
 complex
 of $S_1$. If we adjoin $\beta$ (respectively $\gamma$) to each vertex of 
 $\bf q$ we get a path $\bf q_2$ (respectively $\bf q_1$). Path $\bf q_2$
 connects  $v_2$ to a vertex $u_2$ containing $\alpha$ and $\beta$ and
 path $\bf q_1$ connects $\tilde v$ to a vertex $u_1$ 
 containing $\alpha$ and 
 $\gamma$. The corresponding
  vertices of $\bf q_1$ and $\bf q_2$ are connected by an edge so
  the middle rectangle on Figure \ref{radius 0, case 1}
   splits into a sum of squares of type (C4) and is null-homotopic. 
   We can connect 
   $v_1$ to $u_2$ by an $(\alpha,\beta)$--segment so
 the triangle on Figure \ref{radius 0, case 1} is a closed 
 $\beta$--segment and is also 
 null-homotopic. The part of $\bf p$ between $v_1$ and $\tilde v$ 
 can be replaced
 by the lower path on Figure \ref{radius 0, case 1}. We get a new path
 ${\bf p}^\prime$, which has a smaller number of segments (no
 $\beta$--segment) and is homotopic to $\bf p$ in $X$.  
   
\begin{figure} [ht!]\unitlength0.9pt\small
\cl{\begin{picture}(310,105)
\put(30,85){$v_0$}
\put(50,88){\vector (1,0){30}}
\put(60,95){$\alpha$}
\put(90,85){$v_1$}
\put(110,88){\vector (1,0){50}}
\put(130,95){$\beta$}
\put(165,85){$v_2$}
\put(185,88){\vector (1,0){10}}
\put(200,85){$\tilde v$}
\put(220,88){\vector (1,0){50}}
\put(240,95){$\gamma$}
\put(275,85){$v_3$}
\put(105,75){\vector (1,-1){50}}
\put(98,45){${(\alpha,\beta)}$}
\put(170,35){\vector (0,1){40}}
\put(160,50){$\beta$}
\put(205,35){\vector (0,1){40}}
\put(210,50){$\gamma$}
\put(165,20){$u_2$}
\put(185,23){\vector (1,0){10}}
\put(200,20){$u_1$}
\end{picture}}
\caption{Path of radius 0, Case 1}
\label{radius 0, case 1}
\end{figure}
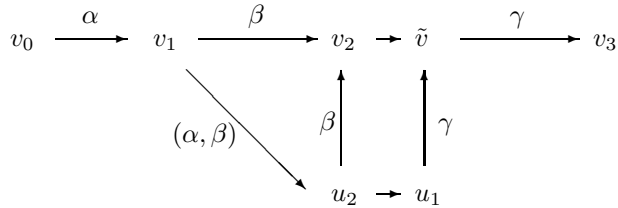    
\par\medskip\noindent
{\bf Case 2}\qua Vertex $v_2$ contains $\gamma$ and $\alpha\cup\gamma$ does not
 separate $S$.\qua If there exists a vertex $v$ which contains $\alpha$ and
 $\beta$ and $\gamma$ we can connect it to $v_1$ and $v_2$ 
 using Lemma \ref{common curve}. We get a closed segment and the
 remaining path has one segment less (Figure \ref{radius 0, case 2}).
  Otherwise
 $\alpha\cup\beta\cup\gamma$ separate $S$ and we can apply Lemma 
 \ref{hexagon relation}. There exist vertices $w_1$ containing $\alpha$
  and $\beta$,
 $w_2$ containing $\beta$ and $\gamma$ and $w_3$ containing $\alpha$ and
 $\gamma$ and a $\beta$--segment from $w_1$ to $w_2$, a $\gamma$--segment
 from $w_2$ to $w_3$ and an $\alpha$--segment from $w_3$ to $w_1$. The sum
 of the segments is null-homotopic. We now connect $v_1$ to $w_1$ by an
 $(\alpha,\beta)$--segment and $v_2$ to $w_2$ by a $(\beta,\gamma)$--segment.
 Thus the second segment of $\bf p$ can be replaced by a sum of an
 $\alpha$--segment and a $\gamma$--segment, and the difference is a closed 
 $\beta$--segment plus a null-homotopic hexagonal path of Lemma
  \ref{hexagon relation}
   (see the right side of Figure \ref{radius 0, case 2}).\\
\begin{figure} [ht!]\unitlength0.9pt\small
\cl{\begin{picture}(370,125)
\put(10,85){$v_1$}
\put(30,88){\vector (1,0){50}}
\put(50,95){$\beta$}
\put(85,85){$v_2$}
\put(105,88){\vector (1,0){50}}
\put(125,95){$\gamma$}
\put(160,85){$v_3$}
\put(25,75){\vector (1,-1){50}}
\put(18,45){${(\alpha,\beta)}$}
\put(90,35){\vector (0,1){40}}
\put(95,50){${(\beta,\gamma)}$}
\put(85,20){$v$}
\put(220,13){\vector (1,0){50}}
\put(280,25){\vector (0,1){25}}
\put(280,70){\vector (0,1){25}}
\put(205,100){\vector (0,-1){75}}
\put(220,108){\vector (1,0){50}}
\put(295,108){\vector (1,0){50}}
\put(210,20){\vector (2,1){60}}
\put(240,20){$\alpha$}
\put(285,35){$\gamma$}
\put(285,80){$(\beta,\gamma)$}
\put(210,65){$(\alpha,\beta)$}
\put(240,115){$\beta$}
\put(315,115){$\gamma$}
\put(245,45){$\beta$}
\put(200,10){$w_1$}
\put(275,10){$w_3$}
\put(275,55){$w_2$}
\put(275,105){$v_2$}
\put(200,105){$v_1$}
\put(350,105){$v_3$}
\end{picture}}
\caption{Path of radius 0, Case 2}
\label{radius 0, case 2}
\end{figure}
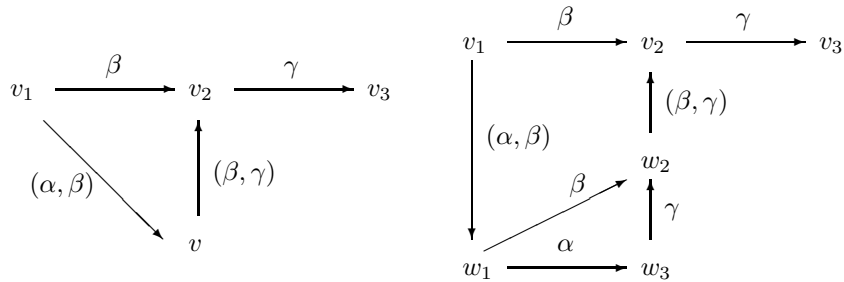 
\par\medskip\noindent 

{\bf Case 3}\qua Vertex $v_2$ contains $\gamma$ and $\alpha\cup\gamma$ 
 separates $S$ into two surfaces $S_1$ and $S_2$.\qua If there were only 
 three segments then, as at the vertex
  $v_1$, the first vertex of the first segment would contain both $\alpha$
  and $\gamma$ and their union would not separate $S$. This contradicts
  our assumptions. It follows in particular that every closed path of 
  radius zero with
  at most three segments (where the common curve of each segment is
  disjoint from a fixed curve of the first segment) is null-homotopic. 
  We may assume that
  the path $\bf p$ has
  a fourth segment with a fixed curve $\delta$ disjoint from $\alpha$.
  Since $[\gamma]=[\alpha]$ we cannot have $|\gamma\cap\delta|=1$.
  Therefore $\delta$ is not involved in the simple move from $v_3$ to
  the next vertex and $v_3$ contains $\delta$. In particular
   $[\gamma]\neq[\delta]$ and $[\alpha]\neq[\delta]$. We may assume that
 $\beta$ lies in $S_1$. If $\delta$ lies in $S_2$ then there is a 
 vertex $w$ 
 which contains $\alpha$ and $\beta$ and $\delta$. We can connect $w$ to 
 $v_1$ by an $(\alpha,\beta)$--segment, to $v_2$ by a $\beta$--segment
  and to $v_3$
 by a $\delta$--segment. We get a new path, homotopic to $\bf p$, 
 which does not contain $\beta$--segment nor $\gamma$--segment (see Figure 
 \ref{radius 0, case 3}, left part.)
 
 Suppose now that $\delta$ lies in $S_1$. Consider the cut-system complex 
 $X_1$ of $S_1$ and choose a vertex $s$ of $X_1$ which contains $\delta$ and 
 a vertex $s^\prime$ of $X_1$ which contains $\beta$. Let $\bf q$ be
 a path in $X_1$ which connects $s$ to $s^\prime$. Let $t$ be a fixed
 vertex of the cut-system complex $X_2$ of $S_2$ (if $X_2$ is not empty.) 
 We add $\alpha$ and all curves of
$t$  to each vertex of the path $\bf q$ and get an $\alpha$--segment in $X$ 
connecting a vertex $w_2$, containing $\delta$, to a vertex $w_2^\prime$,
containing $\beta$. Then we add $\gamma$ and  all curves of
$t$  to each vertex of the path $\bf q$ and get a $\gamma$--segment in $X$ 
connecting a vertex $w_3$, containing $\delta$, to a vertex $w_3^\prime$,
containing $\beta$. We now connect $v_1$ to $w_2$ by an $\alpha$--segment, 
$v_2$ to $w_2^\prime$ by a $\beta$--segment, $v_3$ to $w_3^\prime$ by a 
$\gamma$--segment and $v_4$ to $w_3$ by a $\delta$--segment
(see Figure \ref{radius 0, case 3}, the right side.)
 Corresponding vertices of the two vertical segments
on Figure \ref{radius 0, case 3}, the right side, have a common curve
$\delta_i$, a curve of a vertex of the path $\bf q$ disjoint from 
$\alpha$ and from $\gamma$, and can be connected by a $\delta_i$--segment.
We get a \lq\lq ladder\rq\rq such that each 
 small rectangle in this ladder has radius zero around $\gamma$ and 
 consists of only three segments.  Therefore it is null-homotopic.
  Every other closed path
 on Figure \ref{radius 0, case 3}, the right side, has a similar property.
  We get a new path, homotopic to $\bf p$, 
 which does not contain $\beta$--segment nor $\gamma$--segment.\endpf
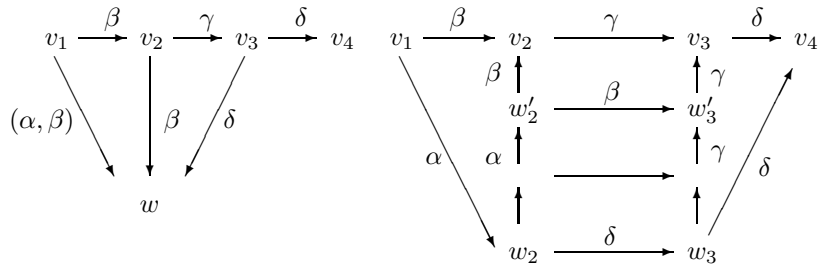
\begin{figure}[ht!]\unitlength0.9pt\small
\cl{\begin{picture}(370,125)
\put(30,110){$v_1$}
\put(45,113){\vector (1,0){20}}
\put(55,118){$\beta$}
\put(70,110){$v_2$}
\put(85,113){\vector (1,0){20}}
\put(95,118){$\gamma$}
\put(110,110){$v_3$}
\put(125,113){\vector (1,0){20}}
\put(135,118){$\delta$}
\put(150,110){$v_4$}
\put(75,105){\vector (0,-1){50}}
\put(80,75){$\beta$}
\put(70,40){$w$}
\put(35,105){\vector (1,-2){25}}
\put(15,75){${(\alpha,\beta)}$}
\put(115,105){\vector (-1,-2) {25}}
\put(105,75){$\delta$}
\put(175,110){$v_1$}
\put(190,113){\vector (1,0){30}}
\put(200,118){$\beta$}
\put(225,110){$v_2$}
\put(245,113){\vector (1,0){50}}
\put(265,118){$\gamma$}
\put(300,110){$v_3$}
\put(320,113){\vector (1,0){20}}
\put(327,118){$\delta$}
\put(345,110){$v_4$}
\put(230,35){\vector (0,1){15}}
\put(230,60){\vector (0,1){15}}
\put(230,90){\vector (0,1){15}}
\put(225,80){$w_2^\prime$}
\put(215,60){$\alpha$}
\put(215,95){$\beta$}
\put(225,20){$w_2$}
\put(300,80){$w_3^\prime$}
\put(305,90){\vector (0,1){15}}
\put(305,60){\vector (0,1){15}}
\put(305,35){\vector (0,1){15}}
\put(310,63){$\gamma$}
\put(265,26){$\delta$}
\put(265,87){$\beta$}
\put(310,95){$\gamma$}
\put(330,55){$\delta$}
\put(180,105){\vector (1,-2){40}}
\put(245,83){\vector (1,0){50}}
\put(190,60){$\alpha$}
\put(245,55){\vector (1,0){50}}
\put(245,23){\vector (1,0){50}}
\put(300,20){$w_3$}
\put(310,30){\vector (1,2){35}}
\end{picture}}
\caption{Path of radius 0, Case 3}
\label{radius 0, case 3}
\end{figure}

\subsection{The general case}
We now pass to the general case and we want to prove it by induction 
on the radius of a closed path.\par\medskip

\noindent {\bf Induction hypothesis 2}\qua A closed path of radius less 
than $m$ is
null-hom\-o\-topic.\par\medskip
We want to prove that a closed path $\bf p$ of radius $m$ around a curve
$\alpha$ is null-homotopic.
The general idea is to construct a  {\em short-cut},
an edge-path which splits $\bf p$ and is close to $\alpha$ and 
to a fixed $\beta$--segment. The first step is to construct one
intermediate curve.
\begin{lemma}\label{curve with n big} Let $\gamma_1$ and $\gamma_2$
 be non-separating curves on $S$  such that  $|\gamma_1\cap\gamma_2|=n>1$.
 Then there exists
 a non-separating
   curve $\delta$ such that  $|\gamma_1\cap\delta|<n$ and
  $|\gamma_2\cap\delta|<n$. Suppose that we are also given 
 non-separating curves  $\alpha$, $\beta$ and an integer $m>0$ such that 
  $|\alpha\cap\beta|\leq m$, $|\gamma_1\cap\alpha|<m$, 
  $|\gamma_2\cap\alpha|\leq m$,
  and  $|\beta\cap\gamma_1|=|\beta\cap\gamma_2|=0$. Then we can find 
  a curve $\delta$ as above which also satisfies $|\delta\cap\alpha|<m$
  and $|\delta\cap\beta|=0$.\end{lemma}
  \proof We orient the curves
 $\gamma_1$ and $\gamma_2$ and split the union $\gamma_1\cup\gamma_2$ into a 
 different union of oriented simple closed curves as follows. We start
 near an intersection point, say $P_1$, on the side of $\gamma_2$ after 
 $\gamma_1$
 crosses it and on the side of $\gamma_1$ before $\gamma_2$ crosses it.
Now we move parallel to $\gamma_1$ to the
  next
  intersection point with $\gamma_2$, say $P_2$. We do not cross $\gamma_2$
  at $P_2$ and move parallel
   to $\gamma_2$, in the positive
   direction, back to $P_1$.
    We get a curve $\delta_1$.
     Now we start near
 $P_2$ and move parallel to $\gamma_1$ until we meet an intersection point,
  say $P_3$,
 which is either equal to $P_1$ or was not met before.
  We do not cross $\gamma_2$ at $P_3$ and move parallel 
to $\gamma_2$, in the positive direction,
 back to $P_2$.  We get a curve
$\delta_2$. And so on.
Curve $\delta_i$ meets $\gamma_1$ near some points of $\gamma_1\cap\gamma_2$,
 but not
near $P_i$ and it meets $\gamma_2$ near some points of  
$\gamma_1\cap\gamma_2$,
 but not near $P_{i+1}$. So $\delta_i$ meets both curves less than $n$ times.
 Let $\bar\gamma$ denote the (oriented) homology class represented by
 an oriented curve $\gamma$ in 
  $H_1(\bar S,\ints)$. We have
$\bar\gamma_1+\bar\gamma_2=\bar\delta_1+\dots+\bar\delta_k$. Now we  repeat
 a similar construction
for the opposite orientation of 
$\gamma_2$ starting near the same point $P_1$.
We get new curves $\epsilon_1$,\dots,$\epsilon_r$ and 
$\bar\gamma_1-\bar\gamma_2=\bar\epsilon_1+\dots+\bar\epsilon_r$.
 Also $\bar\delta_1-\bar\epsilon_1=\bar\gamma_2$.
Combining these equalities in $H_1(\bar S,\ints)$ we get
$\bar\epsilon_1+\sum_{i\neq 1}\bar\delta_i=\bar\gamma_1$,
$\bar\delta_1+\sum_{i\neq 1}\bar\epsilon_i=\bar\gamma_1$,
$\sum_{i\neq 1}\bar\delta_i-\sum_{i\neq 1}\bar\epsilon_i=\bar\gamma_2$.
 A simple closed curve separates $S$ if and only if it represents $0$ in
$H_1(\bar S,\ints)$. Since $\gamma_1$ and $\gamma_2$ are non-separating 
 it follows that either $\delta_1$
 and some $\delta_i$, $i\neq 1$, are not separating or $\epsilon_1$ and 
 some $\epsilon_i$, $i\neq 1$, are not separating. And each of them 
 meets $\gamma_1$ and $\gamma_2$ less than $n$ times, so it can be 
 chosen for $\delta$. If we are also given curves $\alpha$ and $\beta$
 and integer $m>0$ which satisfy the assumptions of the Lemma then 
 $|\gamma_1\cap\alpha|+|\gamma_2\cap\alpha|=\Sigma|\delta_i\cap\alpha|=
 \Sigma|\epsilon_i\cap\alpha|\leq 2m-1$ therefore one of the constructed
 nonseparating curves intersects $\alpha$ less than $m$ times and is
 disjoint from $\beta$.\endpf  

\begin{lemma}\label{curve disjoint} Let $\gamma_1$ and $\gamma_2$
 be disjoint non-separating curves on $S$  such that 
  $\gamma_1\cup\gamma_2$ separates $S$. Then there exists
 a non-separating curve $\delta$ such that  $|\gamma_1\cap\delta|=1$ and
  $|\gamma_2\cap\delta|=1$.  Suppose that we are also given 
 non-separating curves  $\alpha$, $\beta$ and an integer $m>0$ such that 
  $|\alpha\cap\beta|\leq m$,  $|\alpha\cap\beta|=1$ if $m=1$,
  $|\gamma_1\cap\alpha|<m$, $|\gamma_2\cap\alpha|\leq m$,
  and  $|\beta\cap\gamma_1|=|\beta\cap\gamma_2|=0$. Then we can find 
  a curve $\delta$ as above which also satisfies $|\delta\cap\alpha|<m$
  and $|\delta\cap\beta|<m$.\end{lemma}
  \proof
 By our 
assumptions  $\gamma_1\cup\gamma_2$ separates $S$ into two components
 $S_1$ and $S_2$. We can choose a simple arc $d_1$ in $S_1$ which connects 
 $\gamma_1$ to $\gamma_2$ and a simple arc $d_2$ in $S_2$ which connects 
 $\gamma_1$ to $\gamma_2$. Then we can slide the end-points of $d_1$ and 
 $d_2$ along $\gamma_1$ and $\gamma_2$ to make the end-points meet. We
 get a nonseparating curve $\delta$ which intersects $\gamma_1$ and 
 $\gamma_2$ once. Suppose that we are also given curves $\alpha$ and $\beta$
 and an integer $m>0$. We need to alter the arcs $d_1$ and $d_2$, if
 necessary, in order to decrease the intersection of $\delta$ with $\alpha$
 and $\beta$.
We may assume that $\beta$ lies in $S_1$. We have to consider several cases.
\par\medskip\noindent
{\bf Case 1}\qua Curve $\alpha$ lies in $S_1$.\qua Then  $d_2$  is disjoint 
from $\alpha$. If $m=1$ then $|\alpha\cap\beta|=1$ so the union 
$\alpha\cup\beta$ does not separate its regular neighbourhood and does not
separate $S_1$. We can choose $d_1$ disjoint from $\alpha$ and $\beta$ and
then $\delta$ is also disjoint from $\alpha$ and $\beta$.\\
 Suppose that 
$m>1$. If $\alpha$ separates $S_1$ (but does not separate $S$)
then it separates $\gamma_1$ from $\gamma_2$ in $S_1$.
 There exists an arc $d$ in $S_1$ which connects
 $\gamma_1$ with $\gamma_2$ and is disjoint from $\alpha$, if $\alpha$ does
 not separate $S_1$, or meets $\alpha$ once, if $\alpha$ separates
 $S_1$. We choose such an arc $d$ which has minimal number of intersections
 with $\beta$.
 If $|\beta\cap d|>m$ then there exist two points $P$ and $Q$ 
 of $\beta\cap d$, consecutive along $\beta$, and not separated by a point of 
 $\beta\cap\alpha$. We can move along $d$ to $P$ then along $\beta$, without
 crossing $\beta$, to $Q$, and then continue along $d$ to its end. This
 produces an arc which meets $\alpha$ at most once and has smaller number
 of intersections with $\beta$. So we may assume that $|\beta\cap d|= m$
 and that every pair
of points of $\beta\cap d$ consecutive along $\beta$ is separated 
by a point of $\beta\cap\alpha$. We now alter $d$ as follows. Consider 
the intersection $d\cap (\alpha\cup\beta)$. If the first or the last point 
along $d$ of this intersection belongs to $\alpha$ we start from this end 
of $d$. Otherwise we start from any end. We move along $d$ to the first
point, say $P$, of intersection with $\alpha\cup\beta$. If $P\in \alpha$
we continue along $\alpha$, without crossing it, to the next point of 
$\alpha\cap\beta$. Then along $\beta$, without crossing it, to the last 
point, say $Q$, of $\beta\cap d$ on $d$, and then along $d$ to its end. 
The new arc crosses $\beta $ at most once, near $Q$, and crosses $\alpha$
less than $m$ times. If $P\in \beta$ we continue along $\beta$, without
crossing $\beta$, to $Q$ and then along $d$ to its end, which produces 
a similar result. We can choose such an arc for $d_1$ and then the curve 
$\delta$ satisfies the Lemma.
\par\medskip\noindent
{\bf Case 2}\qua Curve $\alpha$ meets $\gamma_1$ or $\gamma_2$.\qua Then $m>1$,
because if $m=1$ and $|\gamma_1\cap\alpha|=0$ and $\alpha$ crosses 
$\gamma_2$ into $S_1$ then it must cross it again in order to exit
$S_1$, and this contradicts $|\gamma_2\cap\alpha|\leq m$. The arcs of 
$\alpha$ split $S_1$ (and $S_2$) into connected components. One of the
components must meet both $\gamma_1$ and $\gamma_2$ (Otherwise the union
of all components meeting $\gamma_1$ has $\alpha$ for a boundary component
and then $\alpha$ is disjoint from $\gamma_1$ and $\gamma_2$.) Choosing
$d_1$ (respectively $d_2$) in such a component we can make them disjoint
from $\alpha$. Now we want to modify $d_1$ in such a way that
$|d_1\cap\alpha|=0$ and
  $| d_1\cap\beta|<m$. There are three subcases.
  \par\medskip\noindent
{\bf Case 2a}\qua There exists an arc $a_1$ of $\alpha$ in $S_1$ which
 connects $\gamma_1$ and $\gamma_2$.\qua Choose $d_1$ parallel
to this arc. It may happen that $d_1$ meets $\beta$ $m$ times.
Then $a_1$ is the only arc of $\alpha$ which meets $\beta$.
We then modify $d_1$ as follows. We move from $\gamma_1$ along
 $d_1$
until it meets $\beta$. Then we turn along $\beta$, away from $a_1$, to
the next point of $a_1$. We turn before crossing $a_1$ and move parallel
 to $a_1$ to $\gamma_2$.
 The new
 arc does not meet $\alpha$
and meets $\beta$ less than $m$ times.
\par\medskip\noindent
{\bf Case 2b}\qua There exists an arc of $\alpha$ in $S_1$ which connects
$\gamma_1$ and $\beta$ and there exists an arc of $\alpha$ which connects 
$\gamma_2$ and $\beta$.\qua Then there exist points $P$ and $Q$ of 
$\alpha\cap\beta$, consecutive along $\beta$, and arcs $a_1$ and $a_2$
of $\alpha$ such that $a_1$ connects $\gamma_1$ to $P$ and $a_2$ connects
$Q$ to $\gamma_2$. We move along $a_1$ to $P$ then along $\beta$, without
crossing $\beta$, to $Q$, and then along $a_2$ to $\gamma_2$. The new arc 
does not meet $\alpha$ and meets $\beta$ less than $m$ times.
\par\medskip\noindent
{\bf Case 2c}\qua If an  arc of $\alpha$ in $S_1$ meets $\beta$ then it  meets 
only $\gamma_1$.\qua (The case of $\gamma_2$ is similar.)
 We consider an arc $d$ in $S_1$ which is disjoint from $\alpha$
and connects $\gamma_1$ and $\gamma_2$. We start at $\gamma_2$ and move 
along $d$ to the first point of intersection with $\beta$. Then we move 
along $\beta$, without crossing it, to the first point of intersection with 
$\alpha$. Then we move along $\alpha$, away from $\beta$, to $\gamma_1$.
 The new
arc does not meet $\alpha$ and meets $\beta$ less than $m$ times. If $\beta$ 
is disjoint from $\alpha$ then $\beta$ is either disjoint from a component
of $S_1-\alpha$ which connects $\gamma_1$ to $\gamma_2$ or 
is contained in it. We can find an arc in the component 
(disjoint from $\alpha$) which connects $\gamma_1$ with $\gamma_2$ and meets
$\beta$ at most once.

 So in each case we have an arc $d_1$ which is disjoint from $\alpha$
 and meets $\beta$ less than $m$ times. We now slide the end-points of $d_1$
 along $\gamma_1$ and $\gamma_2$ to meet the end-points of $d_2$. Each slide 
 can be done along one of two arcs of $\gamma_i$. Choosing suitably we may 
 assume that $d_1$ meets at most $m/2$ points of $\alpha$ sliding along 
 $\gamma_2$ and at most $(m-1)/2$ points of $\alpha$ sliding along $\gamma_1$.
 The curve $\delta$ obtained from $d_1$ and $d_2$ meets $\alpha$ and $\beta$
 less than $m$ times.
 \par\medskip\noindent
{\bf Case 3}\qua The curve $\alpha$ lies in $S_2$.\qua Then $|\alpha\cap\beta|=0$
so we must have $m>1$. We can choose $d_2$ which is disjoint from $\beta$ and 
meets $\alpha$ at most once and we can choose $d_1$ which is disjoint from 
$\alpha$ and meets $\beta$ at most once. 
The curve $\delta$ obtained from $d_1$ and $d_2$ meets $\alpha$ and $\beta$
 less than $m$ times.
\endpf
\par

\begin{lemma}\label{connecting path} Let $\delta_1$ and $\delta_2$
 be non-separating curves on $S$ and let $w_1$ be a vertex of $X$ containing 
$\delta_1$ and let $w_2$ be a vertex of $X$ containing $\delta_2$. Then  there
exists an edge-path ${\bf q}=(w_1=z_1,z_2,\dots,z_k=w_2)$ connecting $w_1$
and $w_2$. Suppose that we are also given 
 non-separating curves  $\alpha$, $\beta$ and an integer $m>0$ such that 
  $|\alpha\cap\beta|\leq m$, $|\alpha\cap\beta|=1$ if $m=1$,
  $|\delta_1\cap\alpha|<m$, $|\delta_2\cap\alpha|\leq m$,
  and  $|\beta\cap\delta_1|=|\beta\cap\delta_2|=0$. Then there exists a path
  $\bf q$ as above and an integer  $j$, $1\leq j < k$,
 such that $d(z_i,\beta)<m$ for all $i$,
$d(z_i,\alpha)<m$ for $1\leq i\leq j<k$ and $z_i$ contains $\delta_2$
for $j<i\leq k$.\end{lemma} 
\proof We shall prove the lemma by induction on $|\delta_1\cap\delta_2|=n$. \\
If $\delta_1=\delta_2$ we can connect $w_1$ and $w_2$ by a 
$\delta_1$--segment, by Lemma \ref{common curve}.\\
If $n=1$ there exist vertices $u_1$, $u_2$ in $X$ which are connected by an
 edge and such that $\delta_1\in u_1$, $\delta_2\in u_2$. Now we can 
 connect $u_1$ to $w_1$ and $w_2$ to $u_2$ as in the previous case.\\
If $n=0$ and $\delta_2\cup\delta_1$ does not separate $S$ then there
exists a vertex $v$ containing both curves $\delta_2$ and $\delta_1$. We can
connect $v$ to $w_1$ and $w_2$ as in the first case.\\
Suppose now that $n=0$ and that $\delta_2\cup\delta_1$ separates $S$.
Then, by Lemma \ref{curve disjoint}, there exists
 a curve $\delta$ such that  $|\delta_2\cap\delta|=|\delta_1\cap\delta|=1$.
We can find a vertex $v$ containing $\delta$ and we can connected $v$ to
$w_1$ and $w_2$ as in the second case. If we are also given curves $\alpha$,
and $\beta$ and an integer $m$ we can choose $\delta$ which also satisfies
$|\alpha\cap\delta|<m$ and  $|\beta\cap\delta|<m$. Then the path obtained
by connecting $v$ to $w_1$ and $w_2$ have all vertices in a distance 
less than $m$ from $\beta$ and in a distance less than $m$ from $\alpha$,
except for the final $\delta_2$--segment which ends at $w_2$ (Curve
$\delta_2$ may have distance $m$ from $\alpha$.)\\
If $n>1$ then by Lemma \ref{curve with n big} there exists a curve
$\delta$ such that $|\delta_1\cap\delta|<n$ and $|\delta_2\cap\delta|<n$.
We choose a vertex $v$ containing $\delta$. By induction on $n$ we
can connect $v$ to $w_1$ and $w_2$. If we are also given curves $\alpha$
and $\beta$ and an integer $m$ then we can find $\delta$ which also 
satisfies
$|\delta\cap\alpha|<m$ and  $|\delta\cap\beta|=0$. By induction on $n$ we 
can connect $w_1$ to $v$ and $v$ to $ w_2$ by a path  
the vertices of which are closer to $\beta$ than $m$, and 
closer to $\alpha$ than $m$ except for a final
 $\delta_2$--segment which ends at $w_2$.\endpf 

As an immediate corollary we get:

\begin{cor}\label{X is connected} Complex $X$ is connected.\end{cor} 

 We need one more lemma before we prove that every closed path is 
null-homotopic in $X$.

 \begin{lemma} \label{better curve} Let $\alpha$, $\beta$, $\gamma$ be
  non-separating curves
  on $S$ such that   $|\alpha\cap\beta|=m$, $|\alpha\cap\gamma|\leq m$,
$|\beta\cap\gamma|=1$. There exists a non-separating curve $\delta$
such that $|\delta\cap\alpha|<m$, $|\delta\cap\beta|=0$ and
$|\delta\cap\gamma|\leq 1$. If $m=1$ then $|\delta\cap\gamma|=0$ and
$[\delta]$ is different from $[\alpha]$, $[\beta]$ and $[\gamma]$.
\end{lemma}

\begin{figure}[ht!]
\relabelbox\small
\epsffile{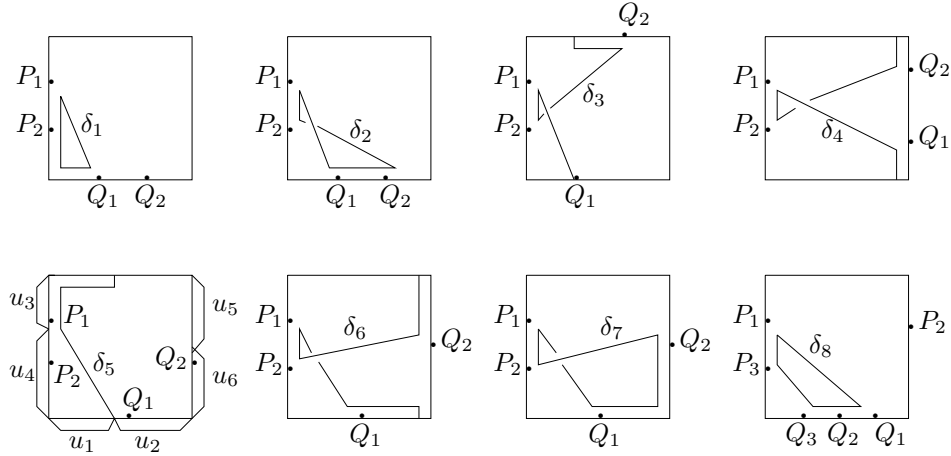}
\relabela <-0.8pt,0pt> {s1}{\llap{$P_2$}\large\bf .}
\relabela <-0.8pt,0pt> {s2}{\large\bf .}
\relabela <-0.8pt,0pt> {s3}{\large\bf .}
\relabela <-0.8pt,0pt> {s4}{\llap{$P_1$}\large\bf .}
\relabela <-0.8pt,0pt> {s5}{\llap{$P_2$}\large\bf .}
\relabela <-0.8pt,0pt> {s6}{\large\bf .}
\relabela <-0.8pt,0pt> {s7}{\large\bf .}
\relabela <-0.8pt,0pt> {s8}{\llap{$P_1$}\large\bf .}
\relabela <-0.8pt,0pt> {s9}{\llap{$P_1$}\large\bf .}
\relabela <-0.8pt,0pt> {s10}{\llap{$P_2$}\large\bf .}
\relabela <-0.8pt,0pt> {s11}{\large\bf .}
\relabela <-0.8pt,0pt> {s12}{{\large\bf .}$Q_2$}
\relabela <-0.8pt,0pt> {s13}{\llap{$P_1$}\large\bf .}
\relabela <-0.8pt,0pt> {s14}{\llap{$P_2$}\large\bf .}
\relabela <-0.8pt,0pt> {s15}{\llap{$P_1$}\large\bf .}
\relabela <-0.8pt,0pt> {s16}{\llap{$P_2$}\large\bf .}
\relabela <-0.8pt,0pt> {s17}{\large\bf .}
\relabela <-0.8pt,0pt> {s18}{\large\bf .}
\relabela <-0.8pt,0pt> {s19}{\llap{$P_1$}\large\bf .}
\relabela <-0.8pt,0pt> {s20}{\llap{$P_2$}\large\bf .}
\relabela <-0.8pt,0pt> {s21}{{\large\bf .}$Q_1$}
\relabela <-0.8pt,0pt> {s22}{{\large\bf .}$Q_2$}
\relabela <-0.8pt,0pt> {s23}{\llap{$P_1$}\large\bf .}
\relabela <-0.8pt,0pt> {s24}{\large\bf .}
\relabela <-0.8pt,0pt> {s25}{\large\bf .}
\relabela <-0.8pt,0pt> {s26}{\large\bf .}
\relabela <-0.8pt,0pt> {s27}{{\large\bf .}$\,P_1$}
\relabela <-0.8pt,0pt> {s28}{{\large\bf .}\kern-2pt\lower 7pt\hbox{$P_2$}}
\relabela <-0.8pt,0pt> {s29}{\large\bf .}
\relabela <-0.8pt,0pt> {s30}{\llap{$Q_2$}\large\bf .}
\relabela <-0.8pt,0pt> {s31}{{\large\bf .}$Q_2$}
\relabela <-0.8pt,0pt> {s32}{\large\bf .}
\relabela <-0.8pt,0pt> {s33}{\llap{$P_3$}\large\bf .}
\relabela <-0.8pt,0pt> {s34}{{\large\bf .}$P_2$}
\relabel {d1}{$\delta_1$}
\relabel {d2}{$\delta_2$}
\relabel {d3}{$\delta_3$}
\relabel {d4}{$\delta_4$}
\relabel {d5}{$\delta_5$}
\relabel {d6}{$\delta_6$}
\relabel {d7}{$\delta_7$}
\relabel {d8}{$\delta_8$}
\relabel {u1}{$u_1$}
\relabel {u2}{$u_2$}
\relabel {u3}{$u_3$}
\relabel {u4}{$u_4$}
\relabel {u5}{$u_5$}
\relabel {u6}{$u_6$}
\relabela <0pt,-2pt> {Q1}{$Q_1$}
\relabela <0pt,-2pt> {Q2}{$Q_2$}
\relabela <0pt,-2pt> {Q3}{$Q_2$}
\relabela <0pt,-2pt> {Q4}{$Q_1$}
\relabela <0pt,-2pt> {Q5}{$Q_1$}
\relabela <0pt,-2pt> {Q8}{$Q_3$}
\relabela <0pt,-2pt> {Q9}{$Q_1$}
\relabela <0pt,-2pt> {Q10}{$Q_2$}
\relabela <0pt,-2pt> {Q12}{$Q_1$}
\relabel {Q13}{$Q_1$}
\relabel {Q15}{$Q_1$}
\relabela <0pt,-1pt> {Q17}{$Q_2$}
\endrelabelbox
\caption{Constructing curve $\delta$}
\label{constructing better curve}
\end{figure}

\proof When we split $S$ along $\gamma\cup\beta$ we get a surface $S_1$ with
a \lq\lq rectangular" boundary component $\partial$ consisting of two 
$\beta$--edges (vertical) and two $\gamma$--edges
(horizontal on pictures of Figure \ref{constructing better curve}).
 We can think of $S_1$ as a rectangle
 with holes and with some handles attached to it. Curve $\alpha$ intersects
 $S_1$ along some arcs $a_i$ with end-points $P_i$ and $Q_i$ on $\partial$.
If, for some $i$, points $P_i$ and $Q_i$ lie on the same $\beta$--edge then
$m>1$ and we can construct a curve $\delta$ consisting of an arc parallel to
$a_i$ and an arc parallel to the
arc of $\beta$ which connects $P_i$ and $Q_i$ passing through the point
 $\gamma\cap\beta$. Then $|\delta\cap\beta|=0$,
$|\delta\cap\alpha|<m$ and $|\delta\cap\gamma|=1$. Recall that if two curves 
intersect exactly at one point then they are both non-separating on $S$.
Therefore $\delta$ satisfies the conditions of the Lemma. If for some $i$
points $P_i$ and $Q_i$ lie on different $\gamma$--edges of $\partial$ then
we can modify the arc $a_i$ sliding its end-point $P_i$ along the $\gamma$
edge to the point corresponding to $Q_i$. We get a closed curve $\delta$
satisfying the conditions of the Lemma. So we may assume that there are 
no arcs $a_i$ of the above types.

Suppose that for every pair $i,j$ the pairs of end-points $P_i,Q_i$ and
$P_j,Q_j$ do not separate each other on $\partial$. Then we can connect 
the corresponding end points by nonintersecting intervals inside a rectangle.
In the other words a regular neighbourhood
of $\alpha\cup\partial$ \ in $S_1$ is a planar surface homeomorphic to 
a rectangle with holes. Since $S_1$ has positive genus there exists a 
subsurface of $S_1$ of a positive genus attached to one hole or 
a subsurface of $S_1$ which connects two holes of the rectangle.
Such a subsurface contains a  curve $\delta$ which is non-separating on $S$
and is 
disjoint from $\alpha$, $\beta$ and $\gamma$ and the homology class
$[\delta]$ is different from  $[\alpha]$, $[\beta]$ and $[\gamma]$.
 This happens in particular 
when $m=1$ because then there is at most one point on every 
edge and the pairs of end points of
arcs do not separate each other.

So we may assume that $m>1$ and that there exists a pair of arcs, 
say $a_1$ and $a_2$, 
such that the pair $P_1,Q_1$ separates the pair $P_2,Q_2$ in $\partial$.
Since an arc $a_i$ does not connect different $\gamma$--edges we must have
two points, say $P_1$ and $P_2$, on the same edge. Suppose that they lie on
a $\beta$--edge, say the left edge. Choosing an intermediate point, 
if there is one, we
may assume that $P_1$ and $P_2$ are consecutive points of $\alpha$ along 
$\beta$. We have different possible configurations of pairs
 of points.
For each of them we construct curves $\delta_i$, as on Figure 
\ref{constructing better curve}. Each $\delta_i$ is disjoint from $\beta$
and intersects $\gamma$ at most once, and if it is disjoint from $\gamma$
it intersects some other curve once. So $\delta_i$ is not-separating.
We shall prove that we can always choose a suitable $\delta_i$ with
$|\delta_i\cap\alpha|<m$. Observe that $\delta_i$ may meet $\alpha$
only along the boundary $\partial$, and not along the arc connecting
$P_1$ to $P_2$.\\
{\bf Case 1}\qua Points $Q_1$ and $Q_2$ lie on the same $\gamma$--edge, say 
lower edge.\qua If there is no point of $\alpha$ on $\gamma$ to the left of 
$Q_1$ then $|\delta_1\cap\alpha|<m$. If there is a point of 
$\alpha$ on $\gamma$ to the left of
$Q_1$ then $|\delta_2\cap\alpha|<m$.\\
{\bf Case 2}\qua  Points $Q_1$ and $Q_2$ lie on different $\gamma$--edges.\qua Then
$|\delta_3\cap\alpha|<m$.\\
{\bf Case 3}\qua  Points $Q_1$ and $Q_2$ lie on the right edge.\qua Then
$|\delta_4\cap\alpha|<m$.\\
{\bf Case 4}\qua One of the points $Q_i$, say $Q_1$, lies on a $\gamma$--edge 
and the other lies on a $\beta$--edge.\qua Let $u_i$, $i=1,\dots,6$ denote
the number of intersection points of $\alpha$ with the corresponding 
piece of $\partial$ on Figure \ref{constructing better curve}. Then
 $u_3+u_4=u_5+u_6=|\alpha\cap\beta|= m$ and $u_1+u_2\leq m$. Also
$|\delta_1\cap\alpha|=u_1+u_4$, $|\delta_5\cap\alpha|=u_1+u_3$,
$|\delta_6\cap\alpha|=u_2+u_5$ and $|\delta_7\cap\alpha|=u_2+u_6$.
Moreover, since $P_2$ and $Q_2$ are connected by an arc of $\alpha$, they
represent different points on $S$ (otherwise it would be the
 only arc of $\alpha$) and $u_4\neq u_6$. It follows that 
 $|\delta_i\cap\alpha|<m$ for $i=1,5,6$ or $7$.\par
We may assume now that for every pair of arcs $(i,j)$
 whose end-points separate each other no two end-points lie on the same 
$\beta$--edge. If $P_i$ and $Q_i$ lie on a $\gamma$--edge and
$P_j$ lie in between then $a_i$ together with the interval of $\gamma$
between $P_i$ and $Q_i$ form a nonseparating curve which meets $\alpha$
less than $m$ times. So we may assume that $P_1$ 
 and $P_2$ lie on 
different $\beta$--edges, say $P_1$ on the left edge and
$P_2$ on the right edge, and $Q_1$ and $Q_2$ lie on the lower edge. 
Replacing $Q_1$ or $Q_2$ by an 
intermediate point, if necessary,  we may also assume that for every 
point $Q_i$ between $Q_1$ and $Q_2$ the corresponding point $P_i$ also
lies between $Q_1$ and $Q_2$.
Now if there is no point of $\alpha$ on the left edge below $P_1$ then 
$|\delta_1\cap\alpha|<m$. If there is such point of $\alpha$ consider the one
closest to $P_1$ and call it $P_3$. Then, by our assumptions, 
 point $Q_3$ lies on the
 lower edge to the left of $Q_2$ and 
$|\delta_8\cap\alpha|<m$.\\
This concludes the proof of the Lemma.\endpf

\begin{proposition}\label{radius m} A path $\bf p$ of radius $m$ around
 $\alpha$ is null-homotopic.
\end{proposition}

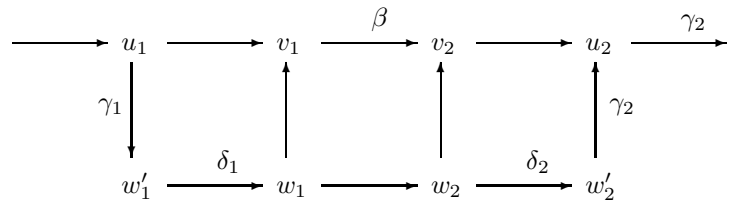
\begin{figure}[ht!]\unitlength0.9pt\small
\cl{\begin{picture}(320,105)
\put(15,88){\vector (1,0){40}}
\put(80,88){\vector (1,0){40}}
\put(145,88){\vector (1,0){40}}
\put(210,88){\vector (1,0){40}}
\put(275,88){\vector (1,0){40}}
\put(65,80){\vector (0,-1){40}}
\put(50,60){$\gamma_1$}
\put(130,40){\vector (0,1){40}}
\put(195,40){\vector (0,1){40}}
\put(260,40){\vector (0,1){40}}
\put(80,28){\vector (1,0){40}}
\put(165,95){$\beta$}
\put(295,95){${\gamma_2}$}
\put(100,35){${\delta_1}$}
\put(230,35){${\delta_2}$}
\put(265,60){$\gamma_2$}
\put(60,85){$u_1$}
\put(125,85){$v_1$}
\put(190,85){$v_2$}
\put(255,85){$u_2$}
\put(60,25){$w_1^\prime$}
\put(125,25){$w_1$}
\put(145,28){\vector (1,0){40}}
\put(190,25){$w_2$}
\put(210,28){\vector (1,0){40}}
\put(255,25){$w_2^\prime$}
\end{picture}}
\caption{Reducing a path of radius m}
\label{reducing a path of radius m}
\end{figure}

\proof Let $v_0$ be a vertex of $\bf p$ containing $\alpha$. We say that
$\bf p$ begins at $v_0$. Let $v_1$ be the first vertex of $\bf p$ which
has distance $m$ from $\alpha$. Let $\bf q$ be the maximal segment of
$\bf p$, which starts at $v_1$ and  contains some fixed curve $\beta$ 
satisfying
$|\beta\cap\alpha|=m$ and such that no vertex of $\bf q$ contains a curve
$\beta^\prime$ satisfying $|\beta^\prime\cap\alpha|<m$. Let $v_2$ be the
last vertex of $\bf q$. Let $u_1$ be the vertex of $\bf p$ preceding $v_1$
and let $u_2$ be the vertex of $\bf p$ following $v_2$. Vertex $u_1$
contains a curve $\gamma_1$ such that $|\gamma_1\cap\alpha|<m$. Vertex
$u_2$ is the first vertex of the second segment which has a fixed curve
$\gamma_2$ such that $|\gamma_2\cap\alpha|\leq m$.
If $u_1$ contains $\beta$ then $|\gamma_1\cap\beta|=0$. Otherwise, since
 $v_1$ does not
contain $\gamma_1$, the move from $u_1$ to $v_1$ involves $\gamma_1$ and 
$\beta$, so $|\gamma_1\cap\beta|=1$. If $v_2$ contains $\gamma_2$ then
$|\gamma_2\cap\beta|=0$. It may also happen that $|\gamma_2\cap\alpha|<m$
and that $\beta\in u_2$. Then also $|\gamma_2\cap\beta|=0$.
 Otherwise $|\gamma_2\cap\beta|=1$. We want to
construct vertices $w_1$, $w_1^\prime$, $w_2$ and $w_2^\prime$ and edge
paths connecting them, as on Figure \ref{reducing a path of radius m}, so
that the rectangles are null homotopic. Then we can replace the part of 
$\bf p$ between $u_1$ and $u_2$ by the path connecting consecutively
$u_1$ to $w_1^\prime$, $w_1^\prime$ to $w_1$, $w_1$ to $w_2$, $w_2$ to 
$w_2^\prime$ and $w_2^\prime $ to $u_2$. We denote the new path by 
${\bf p}^\prime$. \\
In our construction vertex $w_i$ contains a nonseparating curve $\delta_i$
disjoint from $\beta$. 
If $|\gamma_i\cap\beta|=0$ we let $\delta_i=\gamma_i$, $w_i=w_i^\prime=u_i$
and the corresponding rectangle degenerates to an edge.
 If $|\gamma_i\cap\beta|=1$ we proceed as follows. 
 By Lemma \ref{better curve} there exists
 a nonseparating
curve $\delta_i$ such that $|\delta_i\cap\beta|=0$, $|\delta_i\cap\alpha|<m$,
and $|\delta_i\cap\gamma_i|\leq 1$. If 
$|\delta_i\cap\gamma_i|=0$ (this is always the case if $m=1$) then 
$[\delta_i]\neq [\gamma_i]$ and $[\delta_i]\neq[\beta]$ because they 
have different intersections. There exists a vertex $w_i^\prime$ containing
$\delta_i$ and $\gamma_i$ and a vertex $w_i$ containing $\delta_i$ 
and $\beta$. We can connect $u_i$ to $w_i^\prime$ by a $\gamma_i$--segment,
$w_i^\prime$ to $w_i$ by a $\delta_i$--segment and $w_i$ to $v_i$ by a
$\beta$--segment. The corresponding rectangle has radius zero around 
$\delta_i$, so it is null-homotopic by the Induction Hypothesis 2.\\
 If $|\delta_i\cap\gamma_i|=1$
there exist vertices $w_i^\prime$ and $w_i$ which are connected by an edge
and contain $\gamma_i$ and $\delta_i$ respectively. We connect $w_i^\prime$
to $u_i$ by a $\gamma_i$--segment. We now apply Lemma \ref{connecting path}
to vertices $w_i$ and $v_i$ with $\delta_1,\delta_2,\alpha,\beta$
replaced by $\delta_i,\beta,\gamma_i,\beta$ respectively and $m>1$.
There exists a path connecting $w_i$ to $v_i$ such that all vertices of
the path have distance less than $m$ from $\gamma_i$ and $\beta$. The 
corresponding rectangle has radius less than $m$ around $\gamma_i$ so it
is null-homotopic, by the Induction Hypothesis 2.

We now apply Lemma \ref{connecting path} to vertices $w_1$ and $w_2$.
There exists a path ${\bf q}=(w_1=z_1,z_2,\dots,z_k=w_2)$ connecting $w_1$
and $w_2$ such that  $d(z_i,\beta)<m$ for all $i$,
$d(z_i,\alpha)<m$ for $1\leq i\leq j<k$ and $z_i$ contains $\delta_2$
for $j<i\leq k$. In particular the middle rectangle on Figure 
\ref{reducing a path of radius m} has radius less than $m$ around $\beta$
so it is null-homotopic by the Induction Hypothesis 2. All vertices 
of the new part of path ${\bf p}^\prime$ have distance less than $m$ from
$\alpha$ except for the final $\gamma_2$--segment from $w_2^\prime$ to
$u_2$, if $|\gamma_2\cap\beta|=1$, or final $\delta_2=\gamma_2$--segment
of $\bf q$, if $|\gamma_2\cap\beta|=0$ and the right rectangle degenerates.
Thus ${\bf p}^\prime$ has smaller number of segments at the distance 
$m$ from $\alpha$, it has no $\beta$--segment, and it is null-homotopic
by induction on the number of segments of $\bf p$ at the distance $m$
from some curve $\alpha$.\endpf

This concludes the proof of Theorem \ref{Hatcher-Thurston}.

\section{A presentation of $M_{g,1}$}

In this section we shall consider a surface $S$ of genus $g>1$ with one 
boundary component $\partial$. Let ${\mathcal M}_{g,1}$ be the mapping class group of $S$.
Let $X$ be the cut-system complex of $S$ described in the previous section. 
We shall establish a presentation of ${\mathcal M}_{g,1}$ via its action on $X$. 
The action of ${\mathcal M}_{g,1}$ on $X$ is defined by its action 
on vertices of $X$.
 If $v=\lan C_1,\dots,C_g\ran$ is a vertex of $X$
and $g\in{\mathcal M}_{g,1} $ then $g(v)=\lan g(C_1),\dots,g(C_g)\ran$.

We start with some properties of homeomorphisms of a surface. Then
we describe stabilizers of vertices and edges of the action of ${\mathcal M}_{g,1}$ on
$X$. Finally we consider the orbits of faces of $X$ and determine
a presentation of $X$.

In order to shorten some long formulas we shall adopt the following
notation for conjugation: $a*b=aba^{-1}$. As usually
$[a,b]=aba^{-1}b^{-1}$.\par
\begin{remark}\label{proofs}
 {\rm Some proofs of relations between homeomorphisms 
of a surface are left to the reader. The general idea of a proof 
is as follows. We split the surface into a union of disks by a finite
number of curves (and arcs with the end-points on the boundary if the 
surface has a boundary). We prove that the given product of homeomorphisms
takes each curve (respectively arc) onto an isotopic curve (arc), preserving
some fixed orientation of the curve (arc). Then the product is isotopic to 
a homeomorphism pointwise fixed on each curve and arc. But a homeomorphism 
of a disk fixed on its boundary is isotopic to the identity homeomorphism,
relative to the boundary, by Lemma of Alexander. Thus the given product of
homeomorphisms is isotopic to the identity.}\end{remark}

Dehn proved in \cite{Dehn} that every homeomorphism of $S$ is
isotopic to a product of twists. We start with some properties of twists.

\begin{lemma} \label{conjugation}
Let $\alpha$ be a curve on $S$, let $h$ be a homeomorphism and  let  
$\alpha^\prime = h(\alpha)$.
 Then $T_{\alpha^\prime}=hT_\alpha h^{-1}$.\end{lemma}

\proof Since $h$ maps $\alpha$ to $\alpha^{\prime}$ we may assume that 
(up to
isotopy) it also maps a neighborhood
$N$ of $\alpha$ to a neighborhood $N^\prime$ of $\alpha^\prime$. The
 homeomorphism
$h^{-1}$ takes $N^\prime$ to $N$, then $T_\alpha$ maps $N$ to $N$, 
twisting along
$\alpha$, and $h$  takes $N$ back
 to $N^\prime$. Since $T_\alpha$ is supported in $N$, the composite map is
supported in $N^\prime$ and is a Dehn twist about $\alpha^\prime$. \endpf

\begin{lemma}\label{basic relations}
Let $\gamma_1,\gamma_2,\dots,\gamma_k$ be a chain of
curves, ie, the consecutive curves intersect once and non-consecutive 
curves are disjoint. Let $N$
denote the regular neighbourhood of the union of these curves.
 Let $c_i$ denote the twist along $\gamma_i$. Then
the following relations hold:
 \begin{enumerate}
\item[\rm(i)] The ``commutativity relation":
$c_ic_j=c_jc_i$ if $|i-j|>1$.
\item[\rm(ii)] The ``braid relation":
$c_ic_j(\gamma_i)= \gamma_j$, \ and \ \ $c_ic_jc_i=c_jc_ic_j$ 
\ if \ $|i-j|=1$.
\item[\rm(iii)] The \lq\lq chain relation":
 If $k$ is odd then $N$ has
two boundary components, $\partial_1$ and $\partial_2$, and 
$(c_1c_2\dots c_k)^{k+1}=T_{\partial_1}
 T_{\partial_2}$.
If $k$ is even then $N$ has one boundary component $\partial_1$ and 
$(c_1c_2\dots c_k)^{2k+2}=T_{\partial_1}$.
\item[\rm(iv)]
$(c_2c_1c_3c_2)(c_4c_3c_5c_4)(c_2c_1c_3c_2)=
(c_4c_5c_3c_4)(c_2c_1c_3c_2)(c_4c_3c_5c_4)$.
\item[\rm(v)]
$(c_1c_2\dots c_k)^{k+1}=(c_1c_2\dots c_{k-1})^k(c_kc_{k-1}\dots
c_2c_1^2c_2\dots c_{k-1}c_k)=$\par
\  $(c_kc_{k-1}\dots
c_2c_1^2c_2\dots c_{k-1}c_k)(c_1c_2\dots c_{k-1})^k$.
\end{enumerate}\end{lemma}

\proof Relation (i) is obvious. It follows immediately from the definition
of Dehn twist that $c_2(\gamma_1)=c_1^{-1}(\gamma_2)$.
Both statements of (ii) follow from this and from Lemma \ref{conjugation}.
Relation (iii) is a little more complicated. 	It can be proven by the
method explained in Remark \ref{proofs}. Relations (iv) and (v) follow
from the braid relations (i) and (ii) by purely algebraic operations.
 We shall prove (iv).
\par
$(c_2c_1c_3c_2)(c_4c_3c_5c_4)(c_2c_1c_3c_2)=
c_2c_3c_1c_2c_4c_5c_3c_4c_2c_3c_1c_2=$\nl
$\strut c_2c_3c_4c_1c_5c_2c_3c_2c_4c_3c_1c_2=
c_2c_3c_4c_3c_1c_5c_2c_3c_4c_3c_1c_2=$\nl
$\strut c_4c_2c_3c_4c_1c_5c_2c_4c_3c_1c_2c_4=
c_4c_2c_3c_1c_4c_5c_4c_2c_3c_1c_2c_4=$\nl
$\strut c_4c_5c_2c_3c_1c_4c_2c_3c_1c_2c_5c_4=
c_4c_5c_2c_3c_1c_2c_1c_4c_3c_2c_5c_4=$\nl
$\strut c_4c_5c_3c_2c_3c_4c_1c_2c_3c_2c_5c_4=
c_4c_5c_3c_2c_3c_4c_3c_1c_2c_3c_5c_4=$\nl
$\strut c_4c_5c_3c_4c_2c_3c_4c_1c_2c_3c_5c_4=
(c_4c_5c_3c_4)(c_2c_1c_3c_2)(c_4c_3c_5c_4)$.\par\noindent
We now prove (v). We prove by induction that for $s\leq k$ we have
$(c_1c_2\dots c_k)^s=(c_1c_2\dots c_{k-1})^s(c_kc_{k-1}\dots c_{k-s+1})$.
Using (i) and (ii) one checks easily that for $i>1$ we have
$c_i(c_1c_2\dots c_k)=(c_1c_2\dots c_k)c_{i-1}$. Now \par
$(c_1c_2\dots c_k)^{s+1}=(c_1c_2\dots c_{k-1})^s(c_kc_{k-1}\dots c_{k-s+1})
(c_1c_2\dots c_k)=$\nl
 $\strut (c_1c_2\dots c_{k-1})^sc_1c_2\dots c_kc_{k-1}\dots c_{k-s}=
 (c_1c_2\dots c_{k-1})^{s+1}c_kc_{k-1}\dots c_{k-s}$. \\
 For $s=k+1$ we get
 $(c_1c_2\dots c_k)^{k+1}=(c_1c_2\dots c_{k-1})^k(c_kc_{k-1}\dots c_1
 c_1c_2\dots c_k)$. 
 
 This proves the first equality in (v). By (i) and (ii) 
  $(c_kc_{k-1}\dots
c_2c_1^2c_2\dots c_{k-1}c_k)$ commutes with $c_i$ for $i<k$, which implies the second 
equality.
 \endpf

\begin{figure}[ht!]
\relabelbox\small
\epsffile{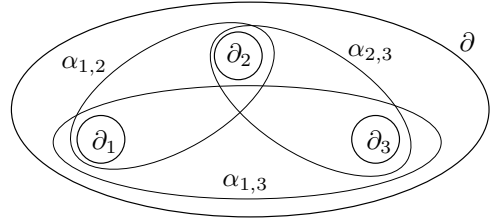}
\relabel {a12}{$\alpha_{1,2}$}
\relabel {a23}{$\alpha_{2,3}$}
\relabela <0pt,-2pt> {a13}{$\alpha_{1,3}$}
\relabela <0pt,-5pt> {d}{$\partial$}
\relabela <0pt,-7pt> {d1}{$\partial_1$}
\relabela <0pt,-7pt> {d2}{$\partial_2$}
\relabela <0pt,-7pt> {d3}{$\partial_3$}
\endrelabelbox
\caption{Lantern relation}
\label{fig lantern}
\end{figure}

\begin{figure}[ht!]
\relabelbox\small
\epsffile{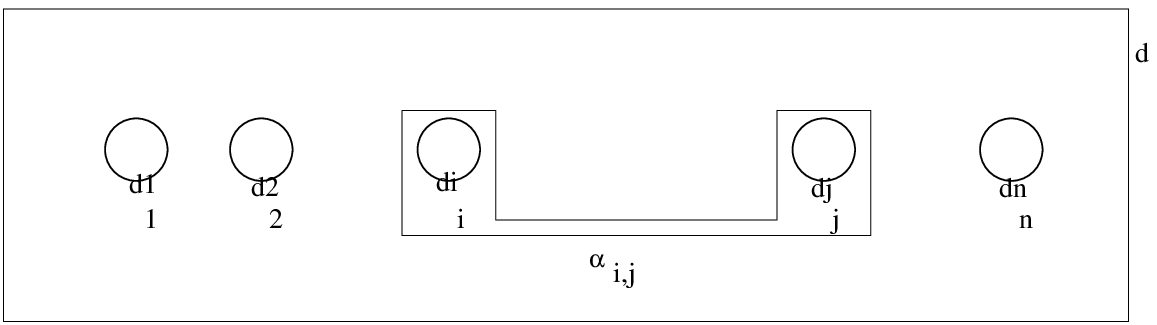}
\relabel {a}{$\alpha_{i,j}$}
\relabela <0pt,-7pt> {d}{$\partial$}
\relabela <0pt,-7pt> {d1}{$\partial_1$}
\relabela <0pt,-7pt> {d2}{$\partial_2$}
\relabela <0pt,-7pt> {di}{$\partial_i$}
\relabela <0pt,-7pt> {dj}{$\partial_j$}
\relabela <0pt,-7pt> {dn}{$\partial_n$}
\endrelabelbox
\caption{Generators of $\cM_{0,n}$}
\label{disk with holes}
\end{figure}

The next lemma was observed by Dennis 
Johnson in \cite{Johnson} and was called a {\em lantern relation}.
 \begin{lemma}\label{lantern} Let $U$ be a disk with the outer 
 boundary $\partial$ and with $3$ inner holes
 bounded by curves $\partial_1,\partial_2,\partial_3$ which form
 vertices of a triangle in the  clockwise order. 
 For $1\leq i<j\leq 3$ let $\alpha_{i,j}$
 be the simple closed curve in $U$ which bounds a neighbourhood of
 the \lq\lq edge" \ $(\partial_i,\partial_j)$ of the triangle
 (see Figure \ref{fig lantern}). Let
$d$ be the twist along $\partial$, $d_i$ the twist along $\partial_i$
and $a_{i,j}$ the twist along $\alpha_{i,j}$.\\ Then
$dd_1d_2d_3=a_{1,2}a_{1,3}a_{2,3}$.
\end{lemma}
We now describe a presentation of the mapping class group of a disk 
with holes.
\begin{lemma}\label{disk with holes presentation} 
Let $U$ be a disk with the outer 
 boundary $\partial$ and with $n$ inner holes
 bounded by curves $\partial_1,\partial_2,\dots,\partial_n$.
For $1\leq i<j\leq n$ let $\alpha_{i,j}$
 be the simple closed curve in $U$ shown on Figure \ref{disk with holes},
  separating two holes $\partial_i$ and 
 $\partial_j$ from the other holes. 
  Let
$d$ be the twist along $\partial$, $d_i$ the twist along $\partial_i$
and $a_{i,j}$ the twist along $\alpha_{i,j}$. Then the mapping class 
group of $U$ has a presentation with generators $d_i$ and $a_{i,j}$ and with
relations \\
{\rm(Q1)}\qua    $[d_i,d_j]=1$ \ and  \ $[d_i,a_{j,k}]=1$  
           for all $i,j,k$.\\
{\rm(Q2)}\qua   pure braid relations\par
       {\rm(a)}\qua $a_{r,s}^{-1}*a_{i,j}=a_{i,j}$
            if \ $ r<s<i<j$ \ or \ $i<r<s<j$,\par
       {\rm(b)}\qua $a_{r,s} ^{-1}*a_{s,j}=a_{r,j}*a_{s,j}$
            if \ $r<s<j$,\par
       {\rm(c)}\qua $a_{r,j}^{-1}*a_{r,s}=a_{s,j}*a_{r,s}$
            if \ $r<s<j$,\par
       {\rm(d)}\qua $[a_{i,j},a_{r,j}^{-1}*a_{r,s}]=1$
           if \ $r<i<s<j$.\\

\end{lemma}

\proof  Relations (Q2) come in place of standard relations in 
the pure braid group on $n$ strings and we shall first prove the 
equivalence of (Q2) to the standard presentation of the pure braid group.
The standard presentation has generators $a_{i,j}$ and relations
(this is a corrected version of the relations in \cite{Birman}):

  (i)\qua $a_{r,s}^{-1}*a_{i,j}=a_{i,j}$
            if \ $ r<s<i<j$ \ or \ $i<r<s<j$,\par
       (ii)\qua      $a_{r,s} ^{-1}*a_{s,j}=a_{r,j}*a_{s,j}$
            if \ $r<s<j$,\par
    (iii)\qua $[a_{r,j},a_{s,j}]=[a_{r,s}^{-1},a_{r,j}^{-1}]$
            if \ $r<s<j$,\par
       (iv)\qua $a_{r,s}^{-1}*a_{i,j}=[a_{r,j},\ a_{s,j}]*a_{i,j}$
           if \ $r<i<s<j$.\\
So relations (a) and (b) are the same as (i) and (ii) respectively.
We can substitute relation (iii) in (iv) and get (d), after cancellation
of $a_{r,s}$. When we substitute (ii) for the first three terms of (iii)
we get (c), after cancellation of $a_{r,s}^{-1}$.\\
We now consider the disk $U$ with holes.
When we glue a disk with a distinguished center 
to each curve $\partial_i$ we get a disk with $n$ distinquished points.
Its mapping class group is isomorphic to the pure braid group $P_{n}$ with 
generators $a_{i,j}$ and with relations (Q2). In the passage from the 
mapping class group of $U$ to the mapping class group of the punctured disk
we kill exactly the twists $d_i$, which commute with everything. One can 
check that the removal of the disks does not affect the relations (Q2) so
the mapping class group of $U$ has a presentation with relations (Q1) and
(Q2).\endpf

We now consider the surface $S$.
We shall fix some curves on $S$. The surface $S$ consists of a disk
with $g$ handles attached to it. For $i=1,\dots,g$
and $j=1,\dots,g-1$ we 
 fix curves $\alpha_i$, $\beta_i$, $\epsilon_j$
  (see Figure \ref{general surface}). Curve $\alpha_i$ is
a meridian curve across the $i$-th handle, $\beta_i$ is a curve along 
the $i$-th handle and $\epsilon_i$ \ runs along $i$-th handle and $(1+i)$-th
handle. Curves $\alpha_1,\alpha_2,\dots,\alpha_g$ form a cut-system.
\par
\begin{figure}[ht!]
\relabelbox\small
\epsffile{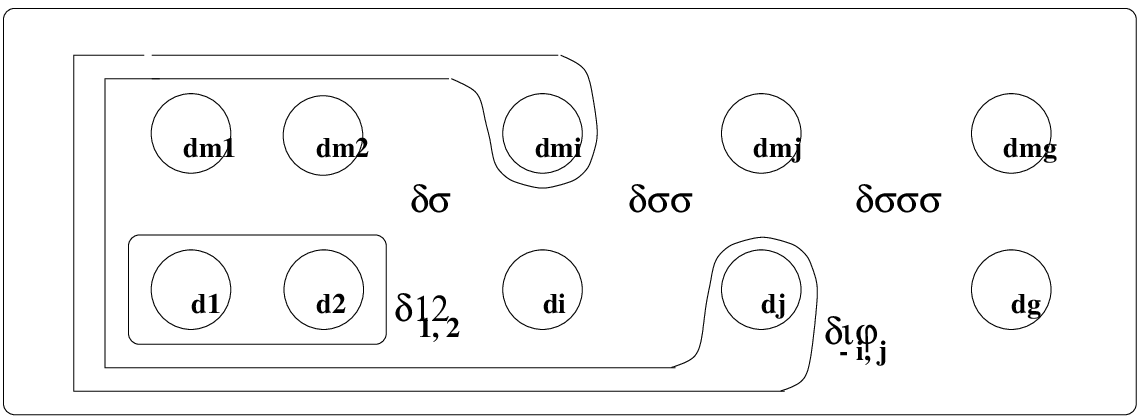}
\relabela <-5pt,4pt> {dm1}{$\partial_{-1}$}
\relabela <-5pt,4pt> {dm2}{$\partial_{-2}$}
\relabela <-5pt,4pt> {dmi}{$\partial_{-i}$}
\relabela <-5pt,4pt> {dmj}{$\partial_{-j}$}
\relabela <-5pt,4pt> {dmg}{$\partial_{-g}$}
\relabela <-5pt,4pt> {d1}{$\partial_1$}
\relabela <-5pt,4pt> {d2}{$\partial_2$}
\relabela <-5pt,4pt> {di}{$\partial_i$}
\relabela <-5pt,4pt> {dj}{$\partial_j$}
\relabela <-5pt,4pt> {dg}{$\partial_g$}
\relabela <0pt,0pt> {d12}{$\delta_{1,2}$}
\relabela <0pt,0pt> {dij}{$\delta_{-i,j}$}
\relabela <0pt,0pt> {ds}{$\ldots$}
\relabela <0pt,0pt> {dss}{$\ldots$}
\relabela <0pt,0pt> {dsss}{$\ldots$}
\extralabel <-15pt, 100pt> {$\partial$}
\endrelabelbox
\caption{Curves $\delta_{i,j}$ on $S_0$}
\label{delta curves}
\end{figure}

We denote by $I_0$ the set of indices $I_0=\{-g,1-g,2-g,\dots,-1,1,2,
\dots,g-1,g\}$.
 When we cut $S$ open along the curves
$\alpha_1,\dots,\alpha_g$ we get a disk $S_0$ with $2g$ holes bounded 
by curves $\partial_i$, $i\in I_0$, where curves $\partial_i$ and 
$\partial_{-i}$ correspond to 
the same curve $\alpha_i$ on $S$ (see Figure \ref{delta curves}). The 
glueing back map identifies $\partial_i$ with $\partial_{-i}$ according 
to the reflection with respect to the $x$--axis. 
Curves on $S$ can be represented on $S_0$. If a curve on $S$ meets 
some curves 
$\alpha_i$ then it is represented on $S_0$ by a disjoint union of 
several
arcs. In particular $\epsilon_i$ is represented by two arcs joining 
$\partial_{-i}$ to $\partial_{-i-1}$ and $\partial_i$ to $\partial_{i+1}$.
We denote by
$\delta_{i,j}$, $i<j\in I_0$, curves on $S$ represented on 
 Figure \ref{delta curves}. Curve $\delta_{i,j}$
separates holes $\partial_i$ and $\partial_j$ from the other holes on $S_0$.
We now fix some elements of ${\mathcal M}_{g,1}$.

\begin{definition}\label{generators}{\rm We denote by $a_i,b_i,e_j$ the
Dehn twists along the curves $\alpha_i,\beta_i,\epsilon_j$ respectively.
 We fix the following elements of ${\mathcal M}_{g,1}$:

$s=b_1a_1a_1b_1$.\nl
$\strut t_i=e_ia_ia_{i+1}e_i$ for $i=1,\dots,g-1$.\nl
$\strut d_{1,2}=(b_1^{-1}a_1^{-1}e_1^{-1}a_2^{-1})*b_2$.

For $i<j\in I_0$ we let\par
 $d_{i,j}=(t_{i-1}t_{i-2}\dots t_1t_{j-1}t_{j-2}\dots t_2)*d_{1,2}$
    if $i>0$,\nl
 $\strut d_{i,j}=(t_{-i-1}^{-1}t_{-i-2}^{-1}\dots t_1^{-1}s^{-1}t_{j-1}t_{j-2}
 \dots t_2)*d_{1,2}$  if $i<0$ and $i+j>0$,\nl
  $\strut d_{i,j}=(t_{-i-1}^{-1}t_{-i-2}^{-1}\dots t_1^{-1}s^{-1}t_{j}t_{j-1}
 \dots t_2)*d_{1,2}$  if $i<0$, $j>0$ and $i+j<0$,\nl
  $\strut d_{i,j}=(t_{-j-1}^{-1}t_{-j-2}^{-1}\dots t_1^{-1}t_{-i-1}^{-1}
  t_{-i-2}^{-1}
 \dots t_2^{-1}s^{-1}t_1^{-1}s^{-1})*d_{1,2}$  if $j<0$,\nl
 $\strut d_{i,j}=(t_{j-1}^{-1}d_{j-1,j}t_{j-2}^{-1}d_{j-2,j-1}\dots 
 t_1^{-1}d_{1,2})*
 (s^2a_1^4)$ if $i+j=0$.}
\end{definition}

 The products described in the above definition represent very 
simple elements of ${\mathcal M}_{g,1}$ and we shall explain their meaning now.
We shall first define special homeomorphisms of $S_0$.

\begin{definition}\label{half-twist} {\rm A {\em half-twist} $h_{i,j}$
along a curve $\delta_{i,j}$ is an isotopy class (on $S_0$ relative to 
its boundary) of 
 a homeomorphism of $S_0$ which is fixed
outside $\delta_{i,j}$ and is equal to a counterclockwise
\lq\lq rotation" by $180$ degrees
inside $\delta_{i,j}$. In particular $h_{i,j}$ switches the two holes 
$\partial_i$ and $\partial_j$ inside $\delta_{i,j}$ so it is not fixed on 
the boundary
of $S_0$, but $h_{i,j}^2$ is fixed on the boundary of $S_0$ and is
isotopic to the full Dehn twist along $\delta_{i,j}$.} \end{definition}

\begin{lemma}\label{action of tk} The result of the action of $t_k$
(respectively $s$) on a curve $\delta_{i,j}$ is the same as the result of the
action of the product of half-twists $h_{k,k+1}h_{-k-1,-k}$ (respectively
the result of the action of $h_{-1,1}$) on $\delta_{i,j}$.
 So  $t_k$  rotates $\delta_{i,j}$ around  $\delta_{k,k+1}$ counterclockwise
and around $\delta_{-k-1,-k}$ conterclockwise and switches the corresponding 
holes. If the pair $(i,j)$ is disjoint from $\{k,k+1,-k,-k-1\}$ than $t_k$
leaves $\delta_{i,j}$ fixed. It also leaves curves $\delta_{k,k+1}$ and
$\delta_{-k-1,-k}$ fixed. In a similar way $s$ rotates $\delta_{i,j}$
counterclockwise around $\delta_{-1,1}$ and switches the holes. If
$(i,j)$ is disjoint from $(-1,1)$ than $s$ leaves $\delta_{i,j}$  fixed.
Also $\delta_{-1,1}$ is fixed by $s$. \end{lemma}
\proof The result of the action can be checked directly.\endpf
\begin{lemma}\label{dij is a twist}
 The element $d_{i,j}$ of ${\mathcal M}_{g,1}$ is equal to the twist along
the curve $\delta_{i,j}$ for all $i,j$.\end{lemma}
\proof We start with the curve $\beta_2$
and apply to it the product of twists $b_1^{-1}a_1^{-1}e_1^{-1}a_2^{-1}$. 
We get the curve 
$\delta_{1,2}$. Therefore, by Lemma \ref{conjugation}, $d_{1,2}$ is 
equal to the twist along $\delta_{1,2}$.
For $i\neq -j$ we start with
$\delta_{1,2}$ and apply consecutive factors $t_i$ and $s$, one at a time.
We check that the result is $\delta_{i,j}$ so $d_{i,j}$ is the twist along
$\delta_{i,j}$, by Lemma \ref{conjugation}.
For $d_{-1,1}$ we use (iii) of Lemma \ref{basic relations}.
The curve $\delta_{-1,1}$
is the boundary of a regular neighbourhood of $\alpha_1\cup\beta_1$ 
 and $(a_1b_1)^6=s^2a_1^4$ by
(i) and (ii) of Lemma \ref{basic relations},
therefore $d_{-1,1}=(a_1b_1)^6$ is the twist along $\delta_{-1,1}$
by (iii) of Lemma \ref{basic relations}. Now we apply a suitable
product of $t_i$'s and $d_{i,i+1}$'s to $\delta_{-1,1}$ and get 
$\delta_{-j,j}$. Therefore $d_{-j,j}$ is equal to the twist along 
$\delta_{-j,j}$ by Lemma \ref{conjugation}. \endpf

We can explain now the relations in Theorem \ref{simple presentation}.
\begin{lemma} The relations {\rm (M1), (M2)} and {\rm (M3)} from Theorem
\ref{simple presentation} are satisfied in ${\mathcal M}_{g,1}$.\end{lemma}
\proof Relations (M1) follow from Lemma \ref{basic relations} (i) and (ii).
Curves $\beta_1,\alpha_1,\epsilon_1$ form a chain. One boundary component
of a regular neighbourhood of $\alpha_1\cup\beta_1\cup\epsilon_1$ is
equal to $\beta_2$. It is easy to check that 
$a_2e_1a_1b_1^2a_1e_1a_2(\beta_2)$ is equal to the other boundary component.
By Lemma \ref{basic relations} (iii) and by Lemma \ref{conjugation} we have 
$(b_1a_1e_1)^4=b_2a_2e_1a_1b_1^2a_1e_1a_2b_2(a_2e_1a_1b_1^2a_1e_1a_2)^{-1}$.
This is equivalent to relation (M2) by
Lemma \ref{basic relations} (v). Consider now relation (M3).
Applying consecutive twists one can check that 
$(b_2a_2e_1b_1^{-1})(\delta_{1,3})=\delta_3$, where
$\delta_3$ is the curve on Figure \ref{general surface}. Thus $d_3$ 
represents the twist along $\delta_3$.
When we cut $S$ along
curves $\alpha_1$, $\alpha_2$, $\alpha_3$ and $\delta_3$ we split off 
a sphere with four holes from surface $S$. Since elements $d_{i,j}$
represent twists along curves $\delta_{i,j}$, relation (M3) follows
from lantern relation, Lemma \ref{lantern}.\endpf
Our first big task is to establish a presentation of a stabilizer of one
vertex of $X$. Let $v_0$
 be a fixed vertex of $X$ corresponding 
to the cut system $\lan\alpha_1,\alpha_2,\dots,\alpha_g\ran$.
 Let $H$ be the stabilizer of $v_0$ in ${\mathcal M}_{g,1}$.

\begin{proposition} \label{H presentation} The stabilizer $H$ of vertex
 $v_0$
admits the following presentation:\par
The set of generators consists of $a_1$, $a_2, \dots,a_g$, 
$s$, $t_1$, $t_2$, $\dots$, $t_{g-1}$ and the $d_{i,j}$'s for
$i<j, \ \ i,j\in I_0$.\par

The set of defining relations consists of:

\noindent
{\rm(P1)}\qua    $[a_i,a_j]=1$ \ and  \ $[a_i,d_{j,k}]=1$  
           for all $i,j,k\in I_0$.\\
{\rm(P2)}\qua   pure braid relations\par
       {\rm(a)}\qua $d_{r,s}^{-1}*d_{i,j}=d_{i,j}$
            if \ $ r<s<i<j$ \ or \ $i<r<s<j$,\par
       {\rm(b)}\qua $d_{r,s} ^{-1}*d_{s,j}=d_{r,j}*d_{s,j}$
            if \ $r<s<j$,\par
       {\rm(c)}\qua $d_{r,j}^{-1}*d_{r,s}=d_{s,j}*d_{r,s}$
            if \ $r<s<j$,\par
       {\rm(d)}\qua $[d_{i,j},d_{r,j}^{-1}*d_{r,s}]=1$
           if \ $r<i<s<j$.\\
{\rm(P3)}\qua  $t_it_{i+1}t_i=t_{i+1}t_it_{i+1}$ for \ $i=1,\dots,g-2$\ 
   and \ $[t_i,t_j]=1$ if \ $1\leq i<j-1<g-1$.\\
{\rm(P4)}\qua   $s^2=d_{-1,1}a_1^{-4}$ and \  
     $t_i^2=d_{i,i+1}d_{-i-1,-i}a_i^{-2}a_{i+1}^{-2}$ \  
      for  $i=1,\dots,g-1$.\\
{\rm(P5)}\qua  $[t_i,s]=1$   for \ $i=2,\dots,g-1$.\\ 
{\rm(P6)}\qua   $st_1st_1=t_1st_1s$.\\
{\rm(P7)}\qua   $[s,a_i]=1$ for $1\leq i\leq g$,
\ \ $t_i*a_i=a_{i+1}$ for  $1\leq i\leq g-1$,\nl
$\strut[a_i,t_j]=1$ \  for $1\leq i \leq g$,  $j\neq i,i-1$.\\ 
{\rm(P8)}\qua $s*d_{i,j}=d_{i,j}$ \ if \ $i\neq \pm1$ \ and \ $j\neq\pm1$
       \  or \ if \ $i=-1$ \ and \ $j=1$,\nl
        $\strut s*d_{-1,j}=d_{1,j}$ \ for \ $2\leq j\leq g$, 
        \ \ $s*d_{i,-1}=d_{i,1}$ \ for \ $-g\leq i\leq -2$,\nl
        $\strut t_k*d_{i,j}=d_{i,j}$ \ if \ $1\leq k\leq g-1$ \ and
          \ $(\,j=i+1=k+1$, \  or \ $j=i+1=-k$ \ \
        or  $i,j\notin\{\pm k,\pm (k+1)\}\,)$, 
        \nl
        $\strut t_k*d_{k,j}=d_{k+1,j}$ \ for \ $1\leq k\leq g-1$ \ and \
         $k+2\leq j\leq g$, \nl
        $\strut t_k*d_{i,-k-1}=d_{i,-k}$ \ for \ $1\leq k\leq g-1$ \ and \
        $-g\leq i\leq -k-2$,\nl
 $\strut t_k*d_{-k-1,k}=d_{-k,k+1}$, \
  $t_k*d_{-k-1,k+1}=d_{k,k+1}*d_{-k,k}$, \ for \ $1\leq k\leq g-1$ \nl
    $\strut t_k*d_{-k-1,j}=d_{-k,j}$ \ for \ $1\leq k\leq g-1$ \ and \
    $j> -k$, \    $j\neq k,k+1$,\nl
    $\strut t_k*d_{i,k}=d_{i,k+1}$ \ for \ $1\leq k\leq g-1$ \ and \
    \ $i<k$, \  $i\neq -k,-k-1$.

\end{proposition}
\proof  An element of $H$ leaves the cut-system
$v_0$ invariant but it may permute the curves $\alpha_i$ and may reverse 
their orientation. Clearly $a_i$ belongs to $H$. One can easily check that 
$t_i(\alpha_i)=\alpha_{i+1}$, $t_i(\alpha_{i+1})=\alpha_i$,
$t_i(\alpha_k)=\alpha_k$ for $k\neq i,i+1$. $s(\alpha_1)$
is equal to $\alpha_1$ with the opposite orientation and $s$ is fixed
on other $\alpha_i$'s. Thus $t_i$'s and $s$ belong to $H$. 
By Lemma \ref{dij is a twist} we know
that $d_{i,j}$ is a twist along the curve $\delta_{i,j}$ so it also belongs 
to $H$. We shall prove in the next section that the relations (P1) -- (P8)
follow from the relations (M1) -- (M3), so they are satisfied
 in ${\mathcal M}_{g,1}$ and thus
also in $H$.
The group $H$ can be defined by two exact sequences.
\begin{equation}\label{sigma sequence}
1\rig \ints_2^g\rig \pm\Sigma_g\rig \Sigma_g\rig 1.
\end{equation}
\begin{equation}\label{H sequence}
1\rig H_0\rig H\rig \pm \Sigma_g\rig 1.
\end{equation}
Before defining the objects and the homomorphisms in these sequences
we shall recall the following fact from group theory.
\begin{lemma}\label{general presentation}
Let 
$1\rig A\rig B\rig  C\rig 1$ be an exact sequence of groups
with known presentations $A=\lan a_i|Q_j\ran$ and $C=\lan c_i| W_j\ran$.
A presentation of $B$ can be obtained as follows:
Let $b_i$ be a lifting of $c_i$ to $B$. Let $R_j$ be a word obtained from
$W_j$ by substitution of $b_i$ for each $c_i$. Then $R_j$ represents 
an element $d_j$ of $A$ which we write as a product of generators 
$a_i$ of $A$.
Finally for every $a_i$ and $b_j$  the conjugate $b_j*a_i$
represents an element 
$a_{i,j}$ of $A$, which we write as a product of the generators $a_i$.\\
Then $B=\lan a_i,b_j|Q_j, R_j=d_j, b_j*a_i=a_{i,j}\ran$.\end{lemma}

We now describe the sequence (\ref{sigma sequence})
and the group $\pm\Sigma_g$. This is the group of permutations of the set
$I_0=\{-g,1-g,\dots,-1,1,2,\dots,g\}$ such that $\sigma(-i)=-\sigma(i)$.
 The homomorphism
$\pm\Sigma_g\to\Sigma_g$ forgets the signs. A generator of the kernel
changes the sign of one letter. The sequence splits, $\Sigma_g$ can be
considered as the subgroup of the permutations which take positive numbers
to positive numbers.
 Let $\tau_i=(i,i+1)$ be a transposition 
in $\Sigma_g$ for $i=1,2,\dots,g-1$. Then \\
(S1)\qua $ [\tau_i,\tau_j]=1$ for $|i-j|>1$,\\
(S2)\qua $\tau_i*\tau_{i+1}=\tau_{i+1}^{-1}*\tau_i$ for $i=1,\dots,g-2$,\\
(S3)\qua $\tau_i^2=1$ for $i=1,\dots,g-1$.

 This defines a presentation of $\Sigma_g$. Further let $\sigma_i$ for 
$i=1,\dots,g$ denote the change of sign of the
$i$-th letter in a signed permutation. Then $\sigma_i^2=1$ and
 $[\sigma_i,\sigma_j]=1$ for all $i,j$. Finally 
the conjugation gives $\tau_i*\sigma_i=\sigma_{i+1}$, 
$\tau_i*\sigma_{i+1}=\sigma_i$
 and $[\sigma_j,\tau_i]=1$ for $j\neq i$ and 
$j\neq i+1$. In fact it suffices to use one generator
 $\sigma=\sigma_1$ and the relations
$\sigma_i=(\tau_{i-1}\tau_{i-2}\dots \tau_1)*\sigma$. We get the relations

\noindent
(S4)\qua $\sigma^2=1$,\\
(S5)\qua $[\sigma,\tau_j]=1$ \ and \
 $[(\tau_i\tau_{i-1}\dots \tau_1)*\sigma,\tau_j]=1$ \ for \ $1\leq i\leq g-1$
 \ and \ $j\neq i$ \ and \ $j\neq i+1$,\\
(S6)\qua $[(\tau_i\tau_{i-1}\dots \tau_1)*\sigma,\sigma]=1$ \ and \
$[(\tau_i\tau_{i-1}\dots \tau_1)*\sigma,
(\tau_j\tau_{j-1}\dots \tau_1)*\sigma]=1$ \ for \  $1\leq i,j\leq g-1$.

Group $\pm\Sigma_g$ has a presentation with generators
 $\sigma,\tau_1,\dots,\tau_{g-1}$
and with defining relations  (S1) -- (S6).

We shall describe now the sequence (\ref{H sequence}). 
 A homeomorphism in $H$
may permute the curves $\alpha_i$ and may change their orientations. 
 We fix an orientation of each curve $\alpha_i$ and 
 define a homomorphism $\phi_1\co H\to
 \pm\Sigma_g$ as follows: a homeomorphism $h$ is mapped onto
  a permutation $i\mapsto \pm j$ if 
$h(\alpha_i)=\alpha_j$ and the sign is \lq\lq$+$" if the
 orientations of $h(\alpha_i)$ and of $\alpha_j$ agree, and
 \lq\lq$-$" otherwise. If $h$ preserves 
 the isotopy class of $\alpha_i$
 and preserves its orientation then it is isotopic to a homeomorphism 
 fixed on $\alpha_i$. 
  The kernel of $\phi_1$ is the subgroup $H_0$ of 
the elements of $H$ represented by the homeomorphisms
which keep the curves $\alpha_1,\alpha_2,\dots,
\alpha_g$ pointwise fixed. We want to find a presentation of $H$ from
the sequence (\ref{H sequence}). We start with a presentation of $H_0$.
An element of $H_0$ induces a homeomorphism of $S_0$.
When we glue back the corresponding pairs of boundary components of $S_0$
we get the surface $S$. This glueing map 
 induces
a homomorphism from the mapping class group of $S_0$
 onto $H_0$.
 
\begin{figure}[ht!]
\relabelbox\small
\epsffile{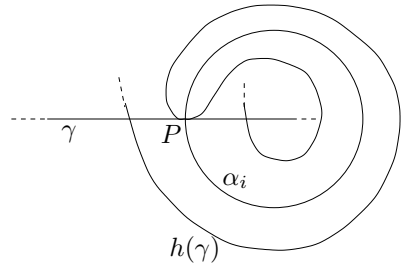}
\relabela <-3pt,0pt> {P}{$P$}
\relabela <0pt,0pt> {a}{$\alpha_i$}
\relabela <0pt,0pt> {g}{$\gamma$}
\relabela <-2pt,0pt> {h}{$h(\gamma)$}
\endrelabelbox
\caption{Self-touching 2--gon}
\label{fake 2-gon}
\end{figure}

  We shall prove that the kernel of this homomorphism is 
  generated by products $d_id_{-i}^{-1}$, where $d_i$ is the
 twist along curve $\partial_i$,  so both twists are
  identified with $a_i$ in $H_0$. It suffices to assume, by induction, 
  that we glue only one pair of boundaries $\partial_i$ and $\partial_{-i}$
on $S_0$ and get a nonseparating curve $\alpha_i$ on $S$. Homeomorphism 
$d_id_{-i}^{-1}$ induces a spin map $s$ of $S$ along $\alpha_i$
(see \cite{Birman}, Theorem 4.3 and Fig. 14). Let $\gamma$ 
be a curve on $S$ which meets $\alpha_i$ at one point $P$. Let $h_0$ be 
a homeomorphism of $S_0$ which induces a homeomorphism $h$ of $S$
isotopic to the identity and fixed on $\alpha_i$. We shall prove that
for a suitable power $k$ the map $s^kh$ is 
fixed on $\gamma$ (after an isotopy of $S$ liftable to $S_0$). 
If $h(\gamma)$
forms a 2--gon with $\gamma$ we can get rid of the 2--gon by an isotopy 
liftable to $S_0$ (fixed on $\alpha_i$). If $h(\gamma)$ and $\gamma$ 
form no 2--gons then they are tangent at $P$. Let us move $h(\gamma)$ off 
$\gamma$, near the point $P$, to a curve $\gamma^\prime$.
 If $\gamma^\prime$ and
$\gamma$ intersect then they form a 2--gon. But then $h(\gamma)$ 
and $\gamma$ form a 
self-touching 2--gon (see Figure \ref{fake 2-gon}). Clearly the spin map
$s$ or $s^{-1}$ removes the 2--gon.
 If $\gamma^\prime$ and $\gamma$ are disjoint
then they bound an annulus. If the annulus contains the short arc of 
$\alpha_i$ between  $\gamma^\prime$ and $\gamma$ then
 $h(\gamma)$ is isotopic to
$\gamma$ by an isotopy liftable to $S_0$. If 
 the annulus contains the long arc of 
$\alpha_i$ between  $\gamma^\prime$ and $\gamma$ then the situation is 
similar to Figure \ref{fake 2-gon}, but $h(\gamma)$ does not meet $\gamma$
outside $P$. The spin map $s$ or $s^{-1}$ takes $h(\gamma)$ onto a curve
isotopic to
$\gamma$ by an isotopy liftable to $S_0$.
 So we may assume, by induction on 
$|h(\gamma)\cap\gamma|$, that $h(\gamma)=\gamma$.
 We proceed as in 
Remark \ref{proofs}. We spilt $S$ into disks by curves $\alpha_i$, $\gamma$,
and other curves disjoint from $\alpha_i$. Homeomorphism $h$ takes each 
curve to an isotopic curve and the 2--gons which appear do not contain 
$\alpha_i$. Thus all isotopies are liftable to $S_0$.

It follows that the kernel is generated by   $d_id_{-i}^{-1}$'s.
   By Lemma \ref{disk with holes
presentation} the mapping class group of $S_0$ has a presentation with
generators $d_k$ and $d_{i,j}$ and with relations (Q1) and (Q2),
where $a_{i,j}=d_{i,j}$, $i,j\in I_0$. Therefore
 $H_0$ has a presentation with generators $a_1,\dots, a_g$ and
$d_{i,j}$ for $i<j\in I_0$ and with relations (P1), (P2). In these 
relations $d_{i,j}$ is represented by a Dehn twist along $\delta_{i,j}$.\par
We come back to sequence (\ref{H sequence}).
We see from the action of $t_i$ and $s$ on $\alpha_j$ that
 we can lift $\tau_i$ to $t_i$ and $\sigma$ to $s$.
Relations (S1) and (S2) lift to (P3). Relations (S3) and (S4) lift to (P4).
Relation (S5) lifts to
 $[(t_it_{i-1}\dots t_1)*s,t_j]=1$ for $j\neq i$ and $j\neq i+1$
and it follows from (P3) and (P5). We shall deal with (S6) a little later.
We now pass to the conjugation of the generators of $H_0$ by $s$ and $t_k$.
 Since $s^2\in H_0$ and
$t_k^2\in H_0$, by (P4),
it suffices to know the result of the conjugation of each generator of $H_0$
by either $s$ or by $s^{-1}$, and by either $t_k$ or by $t_k^{-1}$,
the other follows. The result of the conjugation is described in relations
(P7) and (P8).
 Finally we lift the relation (S6). We start with the case
 $[\tau_1*\sigma,\sigma]$. It lifts to
$t_1st_1^{-1}st_1s^{-1}t_1^{-1}s^{-1}=$ \ (by (P6)) 
$t_1st_1^{-1}ss^{-1}t_1^{-1}s^{-1}t_1=t_1st_1^{-2}s^{-1}t_1^{-1}t_1^2$. We have
$t_1^2\in H_0$, by (P4), and the conjugation of an element of $H_0$
by $s$ and $t_i$ is already determined by (P7) and (P8). So we know how 
to lift $[\tau_1*\sigma,\sigma]=1$. In the general case we have a commutator 
$[(t_it_{i-1}\dots t_1)*s,(t_jt_{j-1}\dots t_1)*s]$. If $i>j$ then, by
(P3) and (P5), this commutator is equal to the conjugate of 
$[t_1*s,s]$ by $t_jt_{j-1}\dots t_1t_it_{i-1}\dots t_2$. This is 
a conjugation of an element of $H_0$ by $t_k$'s, so the result is
determined by (P7) and (P8).\par
This concludes the proof of Proposition \ref{H presentation}.\endpf

We shall prove now that ${\mathcal M}_{g,1}$ acts 
transitively on vertices of $X$, so
 there is only one vertex orbit, and that $H$ acts transitively on edges
incident to $v_0$, so there is only one edge orbit.
When we cut $S$ open along the curves
$\alpha_1,\dots,\alpha_g$ we get a disk with $2g$ holes.
 It is homeomorphic to any other such disk and we can prescribe the 
 homeomorphism on the boundary components, hence every two cut systems 
 can be transformed into each other by a homeomorphism and 
 ${\mathcal M}_{g,1}$ acts 
 transitively on $X^0$.  Let $w$ be a vertex connected to $v_0$ by an edge.
Then $w$ contains a curve $\beta$ which intersects some curve $\alpha_i$
at one point and is disjoint from the other curves of $v_0$. We cut $S$
open along $\alpha_1,\alpha_2,\dots,\alpha_g,\beta$ and get a disk
with $2g-1$ holes, including one \lq\lq big" hole bounded by arcs of 
both $\alpha_i$ and $\beta$. If we cut $S$ along curves belonging to another 
edge incident to $v_0$ we get a similar situation. There exists 
a homeomorphism of one disk with $2g-1$ holes onto the other which 
preserves the identification of curves of the boundary and takes the set
$\alpha_1,\alpha_2,\dots,\alpha_g$ onto itself. It induces a homeomorphism of
$S$ which leaves $v_0$ invariant and takes one edge onto the other. Thus
$H$ acts transitively on the edges incident to $v_0$.

We now fix one such edge and describe its stabilizer. If we replace 
curve $\alpha_1$ of cut-system $v_0$ by curve $\beta_1$ we get a new
cut-system connected to $v_0$ by an edge. We denote the cut-system
 by $v_0^\prime$ and the edge by ${\bf e}_0$. Let $H^\prime$ be
 the stabilizer of ${\bf e}_0$ in $H$. 
\begin{lemma} \label{H prime generators} The stabilizer $H^\prime$ of 
the edge ${\bf e}_0$ is generated by $a_1^2s$, $t_1st_1$, $a_2$,
 $d_{2,3}$, $d_{-2,2}$, $d_{-1,1}d_{-1,2}d_{1,2}a_1^{-2}a_2^{-1}$,
and $t_2,\dots,t_{g-1}$.
\end{lemma} 
 \proof  An element $h$ of $H^\prime$ may
permute curves $\alpha_2,\dots,\alpha_g$ and may reverse their orientations.
It may also reverse simultaneously orientations of both curves $\alpha_1$ 
and $\beta_1$
(it preserves the orientation of $S$ so it must preserve algebraic 
intersection of curves). We check that $a_1^2s$ reverses orientations of 
$\beta_1$ and $\alpha_1$, $t_1st_1$ reverses orientation of $\alpha_2$
and leaves $\beta_1$ and $\alpha_1$ invariant. Elements $t_i$ permute
the curves $\alpha_2,\dots,\alpha_g$. Modulo these homeomorphisms we 
may assume that $h$ is fixed on all curves $\alpha_i$ and $\beta_1$.
It induces a homeomorphism of $S$ cut open along all these curves. We 
get a disk with holes and, by Lemma \ref{disk with holes presentation},
its mapping class group is generated by twists around holes and twists
along suitable curves surrounding two holes at a time. Element $d_{-1,1}=
(a_1^2s)^2$ represents twist around the big hole (the cut along 
$\alpha_1\cup\beta_1$). Conjugates of $a_2$ by $t_j$'s produce twists 
around other holes. The conjugate of $d_{2,3}$ by $(t_1st_1)^{-1}$ is equal
to $d_{-2,3}$ and the conjugate of $d_{-2,3}$ by $(t_2t_1st_1)^{-1}$ 
is equal
to $d_{-3,-2}$. Every $d_{k,k+1}$ can be obtained from these by conjugation
by $t_i$'s, $i>1$.
 It is now clear that every curve $\delta_{i,j}$ with
$i,j\neq\pm1$ can be obtained from $d_{2,3}$ and $d_{-2,2}$ by conjugation
by the elements chosen in the Lemma. So the corresponding twists are 
products of the above generators. We also need twists along curves which 
surround the \lq\lq big" hole and another hole. 
Consider such a curve on $S_0$. It must contain inside both holes 
$\partial_{-1}$ and $\partial_1$, including the arc between them 
corresponding to $\beta_1$, and one other hole. One such curve, 
call it $\gamma$, contains 
$\partial_{-1},\partial_1,\partial_2$. By Lemma \ref{lantern} the twist
along $\gamma$ is equal to $d_{-1,1}d_{-1,2}d_{1,2}a_1^{-2}a_2^{-1}$. Any
other curve which contains $\partial_{-1},\partial_1,\partial_i$ with 
$i>1$ is obtained from $\gamma$ by application of $t_i$'s. The curve
which contains $\partial_{-1},\partial_1,\partial_{-2}$ is obtained from
$\gamma$ by application of $(t_1st_1)^{-1}$ and a curve which contains
$\partial_{-1},\partial_1,\partial_i$ with $i<-2$ is obtained from 
the last curve by application of $t_i^{-1}$'s. Therefore all generators
of $H^\prime$ are products of the generators in the Lemma.\endpf

We now distinguish one more element of 
${\mathcal M}_{g,1}$, which does not belong to $H$.
Let $r=a_1b_1a_1$. The element $r$ is a \lq\lq quarter-twist". It 
switches curves $\alpha_1$ and $\beta_1$, $r^2=sa_1^2=h_{-1,1}$ is
a half-twist around $\delta_{-1,1}$ and $r^4=d_{-1,1}$. Also
$r$ leaves the other curves $\alpha_i$ fixed, so it switches the vertices
of the edge ${\bf e}_0$, $r(v_0)=v_0^\prime$ and $r(v_0^\prime)=v_0$.

We now describe precisely a construction from \cite{Laudenbach} and 
\cite{Hatcher-Thurston} which will let us determine a presentation of 
${\mathcal M}_{g,1}$. This construction was very clearly explained 
in \cite{Heusner}.

Let us consider a free product $(H*\ints)$ where the group $\ints$ is
generated by $r$. An $h$-{\em product} is an element of $(H*\ints)$ with positive 
powers of $r$, so it has a form $h_1rh_2r\dots h_krh_{k+1}$. We have 
an obvious map $\eta\co  (H*\ints)\to {\mathcal M}_{g,1}$ through which the 
$h$--products 
act on $X$. We shall prove that $\eta$ is onto and we shall find
the $h$ products which normally genarate the kernel of $\eta$.

To every edge-path ${\bf p}=(v_0,v_1,\dots, v_k)$ which begins at $v_0$
we assign  an $h$--product
$g=h_1rh_2r\dots h_krh_{k+1}$ such that $h_1r\dots h_mr(v_0)=v_m$
 for $m=1,\dots,k$. We construct it as follows.
There exists $h_1\in H$
such that $h_1(v_0)=v_0$ and $h_1(v_0^\prime)=v_1$. Then 
$h_1r(v_0)=v_1$. Next we transport the second edge to $v_0$.
$(h_1r)^{-1}(v_1)=v_0$ and $(h_1r)^{-1}(v_2)=v_1^\prime$. There exists 
$h_2\in H$ such that $h_2r(v_0)=v_1^\prime$ and $h_1rh_2r(v_0)=v_2$.
And so on. Observe that the elements $h_i$ in the 
$h$--product corresponding 
to an edge-path $\bf p$ are not uniquely determined. In particular 
$h_{k+1}$ is an arbitrary element of $H$.
The construction implies:

\begin{lemma} The generators of $H$ together with the element $r$
generate the group ${\mathcal M}_{g,1}$.
\end{lemma}
\proof Let $g$ be an element of ${\mathcal M}_{g,1}$. Then $g(v_0)$ 
is a vertex of $X$ and 
can be connected to $v_0$ by an edge-path $ {\bf p}=(v_0,v_1,\dots, v_k
=g(v_0))$. Let $g_1=h_1rh_2r\dots h_kr$ be an $h$--product corresponding to
$\bf p$. Then $g_1(v_0)=v_k=g(v_0)$,
  therefore $\eta(g_1^{-1})g=h_{k+1}$
leaves $v_0$ fixed and belongs to the stabilizer $H$ of $v_0$. It follows
that $g=\eta(h_1rh_2r\dots h_krh_{k+1}) $ in ${\mathcal M}_{g,1}$.\endpf

By the inverse process we define an edge-path induced by the $h$--product
$g=h_1rh_2r\dots h_krh_{k+1}$. The edge-path starts at $v_0$ then
 $v_1=h_1r(v_0)$,
$v_2=h_1rh_2r(v_0)$ and so on. The last vertex of the path is 
$v_k=g(v_0)$.
\begin{remark} {\rm An $h$--product $g$ represents an element 
in $H$ if and only
 if 
$v_k=v_0$. This happens if and only if the corresponding edge-path is 
closed.
 We can multiply such $g$
by a suitable element of $H$ on the right and get an $h$--product which 
represents
the identity in ${\mathcal M}_{g,1}$ 
and induces the same edge-path as $g$.}\end{remark}

We now describe a presentation of ${\mathcal M}_{g,1}$.

\begin{theorem}\label{G presentation} The mapping class group 
${\mathcal M}_{g,1}$
admits the following presentation:\par
The set of generators consists of $a_1$, $a_2, \dots,a_g$, 
$s$, $t_1$, $t_2$, $\dots$, $t_{g-1}$, $r$ and the $d_{i,j}$'s for
$i<j, \ \ i,j\in I_0$.\par
The set of defining relations consists of relations {\rm (P1) -- (P8)} of
Proposition \ref{H presentation} and of the following relations:\\
{\rm (P9)}\qua $r$ commutes with  $a_1^2s$, $t_1st_1$, $a_2$,
 $d_{2,3}$, $d_{-2,2}$, $d_{-1,1}d_{-1,2}d_{1,2}a_1^{-2}a_2^{-1}$,
and $t_2,\dots,t_{g-1}$.\\
{\rm (P10)}\qua $r^2=sa_1^2$.\\
{\rm (P11)}\qua $(k_ir)^3=(k_isa_1)^2$ for $i=1,2,3,4$, \ \  where\nl
$\strut k_1=a_1$, $k_2=d_{1,2}$, $k_3=a_1^{-1}a_2^{-2}d_{1,2}
d_{-2,1}d_{-2,2}$, $k_4=a_1^{-1}a_2^{-1}a_3^{-1}d_{1,2}
d_{1,3}d_{2,3}$,\nl
$\strut (rk_5rk_5^{-1})^2=sa_1^2k_5sa_1^2k_5^{-1}$, \ where \ $\strut k_5=a_2t_1d_{1,2}^{-1}$, 
   $(ra_1t_1)^5=(sa_1^2t_1)^4$.\end{theorem}

Relations (P9) --  (P11) 
say that some elements of $H*\ints$ belong to $ker(\eta)$,
 so Theorem \ref{G presentation} claims
that $(H*\ints)$ modulo relations (P9) --  (P11) is isomorphic to
${\mathcal M}_{g,1}$.

Relations  (P9) tell us that $r$ commutes with the generators of
the stabilizer $H^\prime$ of the edge ${\bf e}_0$. We shall prove this
claim now. An element $h$ of $H^\prime$ leaves the edge ${\bf e}_0$ fixed
but it may reverse the orientations of $\alpha_1$ and $\beta_1$. The
element $r^2=sa_1^2$ does exactly this and commutes with $r$. Modulo
this element we may assume that $h$ leaves $\alpha_1$ and $\beta_1$ 
pointwise fixed. But then we may also assume that it leaves some 
neighbourhood of $\alpha_1$ and $\beta_1$ 
pointwise fixed. On the other hand $r$ is equal to the identity
 outside a neighbourhood
of $\alpha_1$ and $\beta_1$ so it commutes with $h$.\par

From relations (P9) and (P10) we get information about $h$--products.

\noindent {\bf Claim 1}\qua If two $h$--products represent the same element in
 ${\mathcal M}_{g,1}$ and induce the same edge-path then
they are equal in $(H*\ints)/(P9)$.

\proof If two $h$--products $g_1=h_1r\dots rh_{k+1}$ and
 $g_2=f_1r\dots rf_{k+1}$
induce the same edge-path ${\bf p}=(v_0,v_1,\dots,v_k)$
then $h_1r(v_0)=f_1r(v_0)$. Therefore
$h_1^{-1}f_1r(v_0)=r(v_0)$ and $h_1^{-1}f_1\in H^\prime$ commutes with $r$
in $(H*\ints)/(P9)$.
Now
 $f_1rf_2=h_1h_1^{-1}f_1rf_2=h_1rf_2^\prime$ in $(H*\ints)/(P9)$.
Therefore $g_2$ and a new $h$--product
$h_1rf_2^\prime rf_3r\dots rf_{k+1}$ are equal in $(H*\ints)/(P9)$ and
induce the same edge-path $\bf p$.
If we apply $r^{-1}h_1^{-1}$ to the vertices $(v_1,v_2,\dots,v_k)$
of $\bf p$ we get a shorter edge-path which starts at $v_0$ and is
induced by both shorter $h$--products $h_2r\dots rh_{k+1}$ and
$f_2^\prime r\dots rf_{k+1}$. 
Claim 1 follows by induction on $k$.\endpf

 Two different edge-paths may be homotopic in the $1$--skeleton $X^1$.
 This means that there is a backtracking $v_i,v_{i+1},v_i$ along
 the edge-path.\par\medskip\noindent
 {\bf Claim 2}\qua If two $h$--products represent the same element in
 ${\mathcal M}_{g,1}$ and induce edge-paths which are equal
 modulo back-tracking then the $h$--products are equal in
 $(H*\ints)/((P9),(P10))$.
 \proof  Consider an $h$--product $g=g_ih_{i+1}rh_{i+2}r$, where
 $g_i$ is an $h$--product inducing a shorter edge-path $\bf p$
  and the edge-path induced by $g$ has a back-tracking
  at the end: \  $g_i(v_0)=v_i$,
  $g_ih_{i+1}r(v_0)=v_{i+1}$, and 
  $g_ih_{i+1}rh_{i+2}r(v_0)=v_i$.
  Clearly the $h$--product $g_ih_{i+1}rr$ induces the same
   edge-path.
  In particular $g_ih_{i+1}r^2(v_0)=v_i$ hence there exists 
  $h^\prime\in H$
  such that $\eta(g_ih_{i+1}r^2h^\prime)=\eta(g_ih_{i+1}rh_{i+2}r)$. 
  Now by Claim 1
  the $h$--products are equal in $(H*\ints)/(P9)$. But 
   $g_ih_{i+1}r^2h^\prime$ is equal in $(H*\ints)/(P10)$ to a shorter 
   $h$--product
  which induces the edge-path $\bf p$. Claim 2 follows by induction
  on the number of back-trackings. \endpf

Relations (P11) correspond to edge-paths of type (C3),
(C4) and (C5). Six relations correspond to six 
particular edge-paths. 

As in (P11) we let $k_1=a_1$, \ $k_2=d_{1,2}$, \ $k_3=d_{1,2}d_{-2,1}
d_{-2,2}a_2^{-2}a_1^{-1}$, \ 
 \nl
 $\strut k_4=a_1^{-1}a_{2}^{-1}a_3^{-1}d_{1,2}d_{1,3}d_{2,3}=d_3$,
\ $k_5=a_2d_{1,2}^{-1}t_1$, \ $k_6=a_1t_1$.\par
We now choose six $h$--products.
For $i=1,2,3,4$ let $g_i=(k_ir)^3$,
$g_5=(rk_5rk_5^{-1})^2$ and $g_6=(rk_6)^5$. Product $g_i$ 
appears in the corresponding relation in (P11).

\begin{figure}[ht!]
\cl{\relabelbox\small
\epsfxsize.8\hsize\epsffile{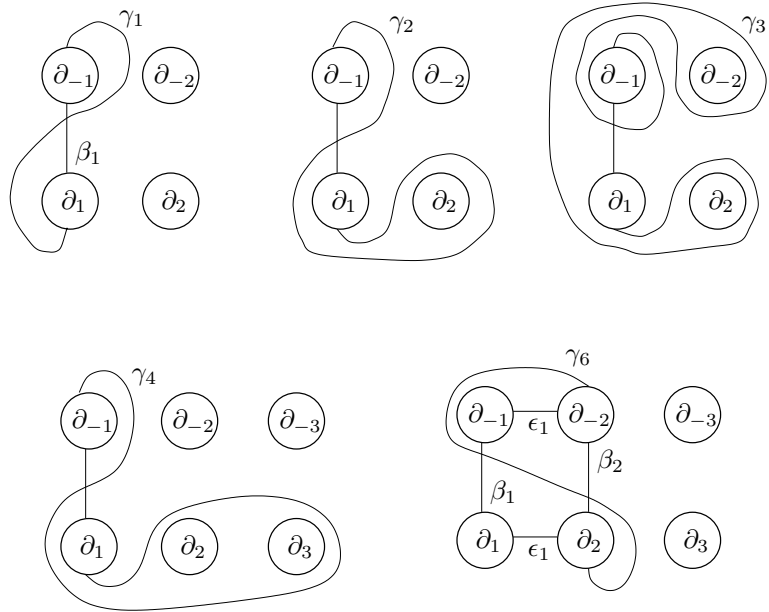}
\relabela <0pt,0pt> {b}{$\beta_1$}
\relabela <0pt,0pt> {b1}{$\beta_1$}
\relabela <0pt,0pt> {b2}{$\beta_2$}
\relabela <0pt,0pt> {g1}{$\gamma_1$}
\relabela <0pt,0pt> {g2}{$\gamma_2$}
\relabela <0pt,0pt> {g3}{$\gamma_3$}
\relabela <0pt,0pt> {g4}{$\gamma_4$}
\relabela <0pt,0pt> {g6}{$\gamma_6$}
\relabela <0pt,0pt> {e}{$\epsilon_1$}
\relabela <0pt,0pt> {e1}{$\epsilon_1$}
\relabela <-3pt,3pt> {d1}{$\partial_1$}
\relabela <-3pt,3pt> {d2}{$\partial_2$}
\relabela <-6pt,3pt> {dm1}{$\partial_{-1}$}
\relabela <-6pt,3pt> {dm2}{$\partial_{-2}$}
\relabela <-3pt,3pt> {2d1}{$\partial_1$}
\relabela <-3pt,3pt> {2d2}{$\partial_2$}
\relabela <-6pt,3pt> {2dm1}{$\partial_{-1}$}
\relabela <-6pt,3pt> {2dm2}{$\partial_{-2}$}
\relabela <-3pt,3pt> {3d1}{$\partial_1$}
\relabela <-3pt,3pt> {3d2}{$\partial_2$}
\relabela <-6pt,3pt> {3dm1}{$\partial_{-1}$}
\relabela <-6pt,3pt> {3dm2}{$\partial_{-2}$}
\relabela <-3pt,3pt> {4d1}{$\partial_1$}
\relabela <-3pt,3pt> {4d2}{$\partial_2$}
\relabela <-6pt,3pt> {4dm1}{$\partial_{-1}$}
\relabela <-6pt,3pt> {4dm2}{$\partial_{-2}$}
\relabela <-3pt,3pt> {5d1}{$\partial_1$}
\relabela <-3pt,3pt> {5d2}{$\partial_2$}
\relabela <-6pt,3pt> {5dm1}{$\partial_{-1}$}
\relabela <-6pt,3pt> {5dm2}{$\partial_{-2}$}
\relabela <-3pt,3pt> {d3}{$\partial_3$}
\relabela <-3pt,3pt> {d31}{$\partial_3$}
\relabela <-6pt,3pt> {dm3}{$\partial_{-3}$}
\relabela <-6pt,3pt> {dm31}{$\partial_{-3}$}
\endrelabelbox}
\caption{Curves $\gamma_i$}
\label{gamma curves}
\end{figure}

For $i=1,\dots,6$ let $\gamma_i$ be the corresponding curve on Figure
\ref{gamma curves}, represented on surface $S_0$ 
($\gamma_5=\beta_2$).

For $i=1,2,3,4$ homeomorphism $k_i$ 
leaves all curves $\alpha_i$ invariant. Also 
$k_ir(\alpha_1)=\gamma_i$, $k_ir(\gamma_i)=\beta_1$,
$k_ir(\beta_1)=\alpha_1$.
 It follows that the $h$--product $g_i$
represents the edge path 
${\bf p}_i=(\lan\alpha_1\ran\to\lan\gamma_i\ran\to\lan\beta_1
\ran\to\lan\alpha_1\ran)$.\\
It is not hard to check that for $i=5,6$ the $h$--product $g_i$ represents
 the edge-path ${\bf p}_i$, where\\
${\bf p}_5=(\lan\alpha_1,\alpha_2\ran\to\lan\beta_1,
\alpha_2
\ran\to\lan\beta_1,\gamma_5
\ran\to\lan\alpha_1,\gamma_5\ran\to
\lan\alpha_1,\alpha_2\ran)$.\\
${\bf p}_6=(\lan\alpha_1,\alpha_2\ran\to\lan\beta_1,\alpha_2
\ran\to\lan\beta_1,\epsilon_1
\ran\to\lan\gamma_6,\epsilon_1\ran\to\lan\gamma_6,\alpha_1\ran\to
\lan\alpha_2,\alpha_1\ran)$.

Since $g_i$ represents a closed edge-path it is equal in ${\mathcal M}_{g,1}$
 to some
element $h_i\in H$. Then $V_i=g_ih_i^{-1}$ also represents ${\bf p}_i$
and is equal to the identity in ${\mathcal M}_{g,1}$. We shall prove in 
the next section 
that $h_i$ is equal in $H$ to the right hand side of  the corresponding 
relation in (P11) and thus $V_i=1$ in $(H*\ints)/((P9),(P10),(P11))$.
\par
We already know, by Theorem \ref{Hatcher-Thurston}, that every closed 
path in $X$ is a sum of paths of type (C3), (C4) and (C5). Some of these 
paths are represented by the paths ${\bf p}_1$ -- ${\bf p}_6$. We say that a 
closed path is {\em conjugate} to a path ${\bf p}_i$ if it has a form 
${\bf q}_1{\bf q}_2{\bf q}_1^{-1}$, where ${\bf q}_1$ starts at $v_0$ 
and ${\bf q}_2$ is the image of ${\bf p}_i$ under the action of some element 
of ${\mathcal M}_{g,1}$.
We shall prove
that every closed path is a sum of paths conjugate to one of the paths 
${\bf p}_1$ -- ${\bf p}_6$ or their inverses. This will imply 
Theorem \ref{G presentation}.\par
\begin{figure}[ht!]
\cl{\relabelbox\small
\hglue-0.2in\epsfxsize1.05\hsize\epsffile{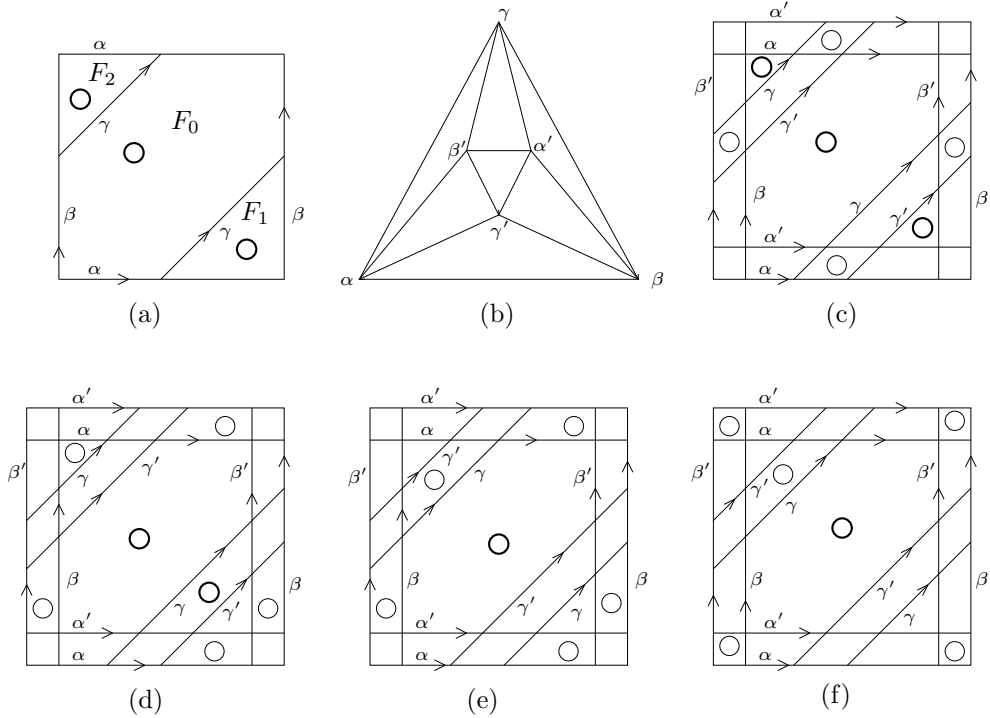}
\relabela <20pt,0pt> {A}{(a)}
\relabela <5pt,0pt> {B}{(b)}
\relabela <10pt,0pt> {C}{(c)}
\relabela <10pt,0pt> {D}{(d)}
\relabela <0pt,0pt> {E}{(e)}
\relabela <0pt,0pt> {F}{(f)}
\relabela <-5pt,0pt> {F1}{$F_0$}
\relabela <-4pt,2pt> {F2}{$F_2$}
\relabela <-5pt,0pt> {F3}{$F_1$}
\relabela <0pt,0pt> {a1}{$\ss \alpha'$}
\relabela <0pt,2pt> {a2}{$\ss \alpha'$}
\relabela <0pt,0pt> {a3}{$\ss \alpha'$}
\relabela <0pt,0pt> {a4}{$\ss \alpha$}
\relabela <0pt,0pt> {a5}{$\ss \alpha$}
\relabela <0pt,0pt> {a6}{$\ss \alpha$}
\relabela <0pt,0pt> {a7}{$\ss \alpha$}
\relabela <0pt,0pt> {a8}{$\ss \alpha$}
\relabela <0pt,0pt> {a9}{$\ss \alpha$}
\relabela <0pt,0pt> {a10}{$\ss \alpha$}
\relabela <-2pt,0pt> {a11}{$\ss \alpha$}
\relabela <0pt,0pt> {a12}{$\ss \alpha$}
\relabela <0pt,0pt> {a13}{$\ss \alpha$}
\relabela <0pt,0pt> {a14}{$\ss \alpha'$}
\relabela <0pt,0pt> {a15}{$\ss \alpha'$}
\relabela <0pt,0pt> {a16}{$\ss \alpha'$}
\relabela <0pt,0pt> {a17}{$\ss \alpha'$}
\relabela <0pt,0pt> {a18}{$\ss \alpha'$}
\relabela <0pt,0pt> {a19}{$\ss \alpha'$}
\relabela <0pt,0pt> {a20}{$\ss \alpha$}
\relabela <0pt,0pt> {b1}{$\ss \beta'$}
\relabela <0pt,0pt> {g1}{$\ss \gamma'$}
\relabela <0pt,0pt> {b2}{$\ss \beta'$}
\relabela <0pt,0pt> {g2}{$\ss \gamma'$}
\relabela <0pt,0pt> {b3}{$\ss \beta$}
\relabela <0pt,0pt> {g3}{$\ss \gamma$}
\relabela <-1pt,0pt> {b4}{$\ss \beta$}
\relabela <0pt,0pt> {g4}{$\ss \gamma$}
\relabela <-1pt,0pt> {b5}{$\ss \beta$}
\relabela <0pt,0pt> {b6}{$\ss \beta$}
\relabela <-1pt,0pt> {b7}{$\ss \beta$}
\relabela <0pt,0pt> {g7}{$\ss \gamma$}
\relabela <-1pt,0pt> {b8}{$\ss \beta$}
\relabela <0pt,0pt> {g8}{$\ss \gamma'$}
\relabela <-1pt,0pt> {b9}{$\ss \beta$}
\relabela <0pt,0pt> {g9}{$\ss \gamma$}
\relabela <0pt,0pt> {b10}{$\ss \beta$}
\relabela <-1pt,0pt> {g10}{$\ss \gamma'$}
\relabela <0pt,0pt> {b11}{$\ss \beta$}
\relabela <0pt,0pt> {g11}{$\ss \gamma$}
\relabela <-1pt,0pt> {b12}{$\ss \beta$}
\relabela <0pt,0pt> {g12}{$\ss \gamma'$}
\relabela <-1pt,0pt> {b13}{$\ss \beta'$}
\relabela <0pt,0pt> {g13}{$\ss \gamma$}
\relabela <-1pt,0pt> {b14}{$\ss \beta'$}
\relabela <0pt,0pt> {g14}{$\ss \gamma$}
\relabela <-1pt,0pt> {b15}{$\ss \beta'$}
\relabela <0pt,0pt> {g15}{$\ss \gamma'$}
\relabela <-1pt,0pt> {b16}{$\ss \beta'$}
\relabela <0pt,0pt> {g16}{$\ss \gamma$}
\relabela <-1pt,0pt> {b17}{$\ss \beta'$}
\relabela <0pt,0pt> {g17}{$\ss \gamma'$}
\relabela <-1pt,0pt> {b18}{$\ss \beta'$}
\relabela <-1pt,1pt> {g18}{$\ss \gamma'$}
\relabela <0pt,0pt> {b19}{$\ss \beta'$}
\relabela <0pt,0pt> {g19}{$\ss \gamma$}
\relabela <0pt,0pt> {b20}{$\ss \beta$}
\relabela <0pt,0pt> {g20}{$\ss \gamma'$}
\relabela <-1pt,-2pt> {g22}{$\ss \gamma$}
\relabela <-5pt,-4pt> {g21}{$\ss \gamma$}
\endrelabelbox}
\caption{Reduction for paths of type (C3)}
\label{squares}
\end{figure}

We shall start with Harer's reduction for paths of type (C3)
 (see \cite{Harer}). We fix curves 
 $\alpha_2,\dots,\alpha_g$
  and consider cut-systems containing one additional curve. Consider 
 a triangular path $(\lan\alpha\ran,\lan\beta\ran,\lan\gamma\ran)$.
 Give orientations to curves $\alpha,\beta,\gamma$ so that the algebraic
 intersection numbers satisfy  $(\alpha,\beta)=(\alpha,\gamma)=1$. Switching 
 $\beta$ and $\gamma$, if necessary, we may also assume $(\beta,\gamma)=-1$.
 We cut $S$ open along $\alpha_2,\dots,\alpha_g$ and get 
 a torus $S_{1,2g-1}$ with $2g-1$ holes. Each square on Figure \ref{squares}
  represents a part of
 the universal cover of a closed torus, punctured above holes 
 of $S_{1,2g-1}$. We show more than one preimage of some curves 
 in the universal cover. A fundamental region is a square bounded by
 $\alpha$ and $\beta$ with $2g-1$ holes, 
 one of them bounded by the boundary $\partial$ of $S$. 
 We may assume that $\gamma$ crosses
$\alpha$ and $\beta$ in two distinct points. Then $\gamma$ splits 
the square into
three regions: $F_1$ to the right of $\gamma$ after $\gamma$ crosses
$\alpha$, $F_2$ to the left of $\gamma$ after $\gamma$ crosses $\beta$,
and $F_0$ (see Figure \ref{squares} (a)). Reversing the orientations
  of all curves we can switch the regions $F_1$ and $F_2$. Let $l_i$ be the 
number of holes in $F_i$ for $i=0,1,2$. We want to prove that every
triangular
path is a sum of triangular 
paths with $l_1\leq 2$, $l_2=0$,
and with the hole $\partial$ not in the region $F_1$.
We shall prove it by induction on $l_1+l_2$.
A thin circle on Figure \ref{squares} denotes
a single hole and a thick circle denotes all remaining holes of the
region.

Suppose first that $l_1>1$.
Move curves $\alpha$, $\beta$, $\gamma$ off itself (on a closed torus) 
in such
a way that the region bounded 
by $\alpha$ and $\alpha^\prime$
contains only one hole, belonging to $F_1$, the region bounded by $\beta$ and
  $\beta^\prime$ contains another hole of $F_1$ and the
   region bounded by $\gamma$ and 
$\gamma^\prime$ contains both of these holes. Up to an isotopy (translating
holes and straightening curves) the situation looks like on
 Figure \ref{squares} (c).

Now consider Figure \ref{squares} (b) --- an octahedron projected onto one
  of its faces. All faces of the octahedron are \lq\lq triangles"
in $X$. In order to understand regions $F_0$, $F_1$, $F_2$ corresponding to
each face we translate the fundamental domain to a suitable square in 
the universal cover. 
For every face different from
$(\alpha,\beta,\gamma)$ the region $F_1$ has at least one hole less than
$l_1$ and the region $F_0$  has all of its original holes and possibly
some more. Clearly the boundary of the face $(\alpha,\beta,\gamma)$ is
a sum of conjugates of the boundaries of the other faces.  So by induction 
we may assume that $l_1\leq 1$.
By symmetry we may assume that $l_2\leq 1$ and $l_1>2$. We chose  new curves 
 $\alpha^\prime$, $\beta^\prime$ and $\gamma^\prime$ whose 
liftings are shown on Figure \ref{squares} (d). Consider again the octahedron
on Figure \ref{squares} (b). Now region $F_2$ is fixed for all faces
of the octahedron and region $F_0$  has all of its original holes and
some additional holes for all faces different from $(\alpha,\beta,\gamma)$.
So by induction we may assume that $l_1\leq 2$ and $l_2\leq 1$.
Suppose $l_1=2$ and $l_2=1$ and choose new curves 
 $\alpha^\prime$, $\beta^\prime$ and $\gamma^\prime$ as on 
 Figure \ref{squares} (e). Again consider the octahedron. Now for each 
 face different from $(\alpha,\beta,\gamma)$ at least one of the regions
 $F_1$ and $F_2$ looses at least one of its holes. Suppose now that
  each region has exactly one hole. Consider curves 
$\alpha^\prime$, $\beta^\prime$ and $\gamma^\prime$ on Figure
\ref{squares} (f). Again for each face of the corresponding octahedron
different from the face $(\alpha,\beta,\gamma)$ at least one of the regions
 $F_1$ and $F_2$ looses at least one of its holes. So we may assume that
 $l_2=0$ and $l_1\leq 2$.
 Finally it
  may happen that a hole in $F_1$ is bounded by $\partial$. We can
  isotop $\gamma$ (on the torus with holes) over $F_2$ to the other side
   of the intersection of $\alpha$ and $\beta$.
 Clearly $F_1$ becomes now $F_0$. We can repeat the
previous reduction and the hole $\partial$ will remain in $F_0$. Thus we are
left with four cases:\\
$l_1=l_2=0$,\\
$l_1=1$, $l_2=0$,\\ 
$l_1=2$, $l_2=0$ and the two holes in $F_1$ correspond to the same curve
$\alpha_i$,\\
$l_1=2$, $l_2=0$ and the two holes in $F_1$ correspond to different
curves $\alpha_i$ and $\alpha_j$.

 Clearly each case is uniquely determined up to
homeomorphism (of one square with holes onto another, preserving $\gamma$).
Paths ${\bf p}_i$, $i=1,2,3,4$ represent triangle paths of these 
four types with $\alpha$ corresponding to $\alpha_1$, $\beta$
corresponding to $\beta_1$ and $\gamma$ corresponding to $\gamma_i$.
Given any path of type (C3) we can map it onto a path which starts at $v_0$.
Applying an element of $H$ we may assume that the second vertex is 
$v_0^\prime$.
Then the path is of the type considered in the above reduction and it is
 a sum of conjugates of ${\bf p}_1$, ${\bf p}_2$, ${\bf p}_3$, 
${\bf p}_4$  and their inverses
(we may need to switch $\beta$ and $\gamma$ in the proof).\par

Consider  now the path ${\bf p}_5$. When we cut $S$ open along
 $\alpha_1,\alpha_2,
\dots,\alpha_g,\beta_1$ and $\gamma_5$ we get a disk with two
\lq\lq big" holes and $2g-4$ \lq\lq small" holes. Any other path $\bf q$
of type
(C4) produces a similar configuration. 
There exists $g\in G$ which takes ${\bf p}_5$ on $\bf q$. Therefore 
every path of type (C4) is a conjugate of ${\bf p}_5$.\par

Consider now ${\bf p}_6$ and another path $\bf q$ of type (C5). When we cut
$S$ open along the first four changing curves of ${\bf p}_6$: $\alpha_1,
\alpha_2,\beta_1,\epsilon_1$ and along all fixed curves of the
cut systems we get a disk with one \lq\lq big" hole and $2g-4$ \lq\lq
small" holes. When we do the same for the curves in $\bf q$ we get 
a similar configuration. 
We may 
map $\bf q$ onto a path 

\begin{figure}[ht!]\unitlength0.9pt\small  
\cl{\begin{picture}(380,175)
\put(80,20){$\lan\alpha_1,\gamma\ran$}  
\put(220,23){\vector(-1,0) {25}}
\put(280,20){$\lan\epsilon_1,\gamma\ran$}  
\put(95,50){\vector(0,1) {100}}
\put(315,50){\vector(0,1) {100}}
\put(80,170){$\lan\alpha_2,\alpha_1\ran$}
\put(120,40){\vector(1,1) {20}}
\put(285,40){\vector(-1,1) {20}}
\put(120,70){$\lan\alpha_1,\gamma_5\ran$}
\put(220,73){\vector(-1,0) {25}}
\put(240,70){$\lan\epsilon_1,\gamma_5\ran$}
\put(150,90){\vector(0,1) {20}}
\put(265,90){\vector(0,1) {20}}
\put(120,120){$\lan\alpha_1,\gamma_6\ran$}
\put(220,123){\vector(-1,0) {25}}
\put(240,120){$\lan\epsilon_1,\gamma_6\ran$}
\put(120,160){\vector(1,-1) {20}}
\put(285,160){\vector(-1,-1) {20}}
\put(155,173){\vector(1,0) {20}}
\put(185,170){$\lan\alpha_2,\beta_1\ran$}
\put(250,173){\vector(1,0) {20}}
\put(280,170){$\lan\epsilon_1,\beta_1\ran$}
\end{picture}}
\caption{Reduction for paths of type (C5)}
\label{pentagon}
\end{figure}
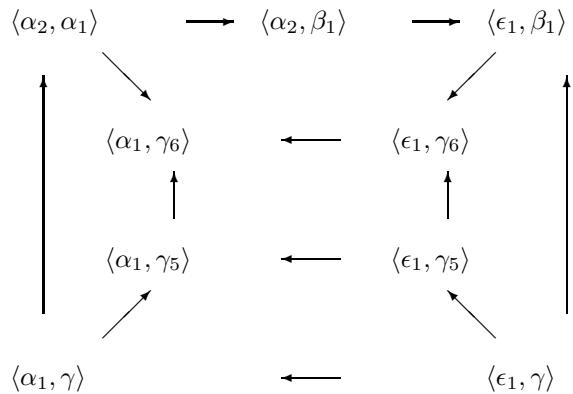

$(\lan\alpha_1,\alpha_2\ran\to\lan\beta_1,\alpha_2
\ran\to\lan\beta_1,\epsilon_1
\ran\to\lan\gamma,\epsilon_1\ran\to\lan\gamma,\alpha_1\ran\to
\lan\alpha_2,\alpha_1\ran)$\\
in which only the curve $\gamma$
is different from the corresponding curve $\gamma_6$ in ${\bf p}_6$ and 
all other curves are as in ${\bf p}_6$. 
Consider curve $\gamma_5=\beta_2$
on Figure \ref{gamma curves}. Curve $\gamma_5$ can be homotop onto the union
  of $\beta_1$, $\epsilon_1$ and a part of $\alpha_1$.
   Since 
  $\gamma$ intersects $\beta_1$ once and does not intersect 
  $\alpha_1$ nor $\epsilon_1$ it must intersect $\gamma_5$ once.
  Therefore we may form a subgraph of $X$ as on Figure \ref{pentagon}.

  Here bottom and middle square edge-paths are of type (C4).
  In the square edge-path on the left side (and on the right side)
 only one curve changes so the path
is a sum of triangles by Proposition \ref{genus1}.
   The top pentagon edge-path is equal to ${\bf p}_6$.
 The new edge-path coincides with the outside pentagon. It is a sum of
 conjugates of the other edge paths.\par

We can now complete the proof of Theorem \ref{G presentation}. Let 
$W\in H*\ints$ be such that $\eta(W)=1$. We want to prove that $W=1$
in $(H*\ints)/((P9),(P10),(P11))$. Modulo (P10) we can write $W$ as 
an $h$--product $g$ which represents a closed edge path $\bf p$. 
By Theorem \ref{Hatcher-Thurston} the edge-path $\bf p$ is a sum 
of conjugates of path of type (C3), (C4) and (C5). By the above 
discussion it is a sum of conjugates of the special paths
${\bf p}_1$ -- ${\bf p}_6$ and their inverses, modulo backtracking.
So, modulo backtracking, 
${\bf p}=\Pi{\bf q}_if_i({\bf p}_{j_i}^{\pm1}){\bf q}_i^{-1}$
for some $f_i\in {\mathcal M}_{g,1}$.
Each path ${\bf q}_if_i({\bf p}_{j_i}^{\pm1}){\bf q}_i^{-1}$ can be
represented by an $h$--product $g_i$. Since the path is closed  we
have $\eta(g_i)\in H$, so we can correct the $h$--product $g_i$ and
assume that $\eta(g_i)=1$. The product $\Pi g_i$ represents
the path $\Pi{\bf q}_if_i({\bf p}_{j_i}^{\pm1}){\bf q}_i^{-1}$, so,
by Claims 1 and 2, $W$ is equal to 
$\Pi g_i$ in $(H*\ints)/((P9),(P10))$. It suffices to prove that
$g_i=1$ in $(H*\ints)/((P9),(P10),(P11))$. 
Let $k=j_i$ and suppose that $g_i$ represents a path  
${\bf q}_if_i({\bf p}_k){\bf q}_i^{-1}$.
 Let $u_i$ be an 
$h$--product representing ${\bf q}_i$. Then  $u_i(v_0)=f_i(v_0)$ 
is the first
vertex of $f_i({\bf p}_k)$, hence $\eta(u_i)^{-1}f_i=h_i\in H$. 
Recall that ${\bf p}_k$ can be represented by an $h$--product $V_k$
which is equal to the identity in $(H*\ints)/((P9),(P10),(P11))$. 
The $h$--product $u_ih_iV_k$ represents
${\bf q}_if_i({\bf p}_k)$ and there exists an $h$--product $u_ih_iV_kw_i$
such that $\eta(u_ih_iV_iw_i)=1$ and $u_ih_iV_kw_i$ represents the edge-path 
${\bf q}_if_i({\bf p}_k){\bf q}_i^{-1}$. Clearly $u_ih_iw_i$ represents the
edge-path ${\bf q}_1{\bf q}_1^{-1}$ which is null-homotopic by backtracking.
By Claims 1 and 2 the $h$--product  $g_i$ is equal to $u_ih_iV_kw_i$ in
$(H*\ints)/((P9),(P10))$,  $u_ih_iV_kw_i$ is equal to
$u_ih_iw_i$ in $(H*\ints)/((P9),(P10),(P11))$ and $u_ih_iw_i=1$ in
$(H*\ints)/((P9),(P10))$.

The path inverse to ${\bf p}_k$ is represented by some other $h$--product 
$V_k^{\prime}$ but then $V_kV_k^{\prime}$ represents a path contractible
by back-tracking. Thus $V_kV_k^\prime=1$ in $(H*\ints)/((P9),(P10))$
and $V_k=1$ in $(H*\ints)/((P9),(P10),(P11))$ hence also
$V_k^{\prime}=1$ \\
in $(H*\ints)/((P9),(P10),(P11))$.\par
This concludes the proof of Theorem \ref{G presentation}.\endpf

\section {Reduction to a simple presentation}

Let us recall the following obvious direction of Tietze's Theorem.

\begin{lemma} \label{different presentations} Consider a presentation of 
$G$ with generators $g_1,\dots,g_k$ and relations $R_1,\dots,R_s$. If we 
add another relation $R_{s+1}$, which is valid in $G$, we get
another presentation of $G$. If we express some element $g_{k+1}$ of $G$
as a product $P$ of the generators $g_1,\dots,g_k$ we get a new 
presentation of $G$ with generators $g_1,\dots,g_{k+1}$ and with defining
 relations $R_1,\dots,R_s,g_{k+1}^{-1}P$. Conversely suppose that we can 
 express some
 generator,
say $g_k$, as a product $P$ of $g_1,\dots,g_{k-1}$. Let us replace
every appearence of $g_k$ in each relation $R_i$, $i=1,\dots,s$
by $P$, getting a new relation $P_i$. Then $G$ has a presentation 
with generators $g_1,\dots,g_{k-1}$ and with relations $P_1,\dots,P_s$.
\end{lemma}  

We start with the presentation of the mapping class group
$\cM_{g,1}$ established in Theorem \ref{G presentation}. The 
generators represent the mapping classes of corresponding homeomorphisms
of the surface $S=S_{g,1}$ represented on Figure \ref{general surface}.
We adjoin additional generators
\begin{equation}\label{extra generators} b_1=a_1^{-1}ra_1^{-1},\ \ 
b_2=(t_1a_1b_1)*d_{1,2},\ \ e_1=(rd_{1,2}a_2^{-1})*b_2\end{equation}
$\kern 1.85cm e_{i+1}=(t_it_{i+1})*e_i \ \ \hbox{for} \ \ i=1,\dots,g-2$.
\par
These generators also represent the corresponding twists in $\cM_{g,1}$.
We adjoin the relations (M1) -- (M3). Now, by Theorem 
\ref {G presentation}, generators $s,t_i,d_{i,j},r$ can be expressed 
in $\cM_{g,1}$ by the formulas from 
 Definition \ref{generators}. We substitute for each of these generators
 the corresponding product of 
 $b_2,b_1,a_1,e_1,a_2,\dots,a_{g-1},e_{g-1},a_g$
  in all relations (P1) -- (P11) and 
in (\ref{extra generators}). We check easily that the relations
(\ref{extra generators}) become trivial modulo the relations (M1) -- (M3)
(the last one will be proven in (\ref{A12}).)
We want to prove that
  all relations (P1) -- (P11) follow from
relations (M1) -- (M3). 

\begin{remark} {\rm We shall establish many auxiliary relations 
of increasing complexity which follow from the relations (M1) -- (M3).
 We shall explain some 
standard technique which one can use (some proofs will be left to the 
reader). From the braid relation $aba=bab$ one can derive several
other useful relations, like: $a*b=b^{-1}*a$, $a*(b^2)=b^{-1}*(a^2)$,
$(ab)*a=b$.
When we want to prove that $[a,b]=0$ we shall usually 
try to prove that $a*b=b$. The relation $aba=bab$ tells us that $a$ 
can \lq\lq jump" over $ba$ to the right becoming $b$. By consecutive 
jumping to the right we can prove that $a_1(e_1a_1a_2e_1e_2a_2)=
(e_1a_1a_2e_1e_2a_2)e_2$. We also get
$e_1a_1a_2e_1e_2a_2=e_1a_2e_2a_1e_1a_2$ and
$(b_1a_1e_1a_2)*b_2=(b_2^{-1}a_2^{-1}e_1^{-1}a_1^{-1})*b_1$.
We shall say that some relations follow by (J) -- jumping, if they follow
{\em easily} from (M1) by the above technique.}\end{remark}

We start the list of the auxiliary relations.
\begin{equation}\label{A1}
t_i*a_i=a_{i+1}, \ t_i*a_{i+1}=a_i, \ t_i*a_k=a_k 
\ \ \hbox{for} \ \ k\neq i,i+1,\end{equation}
$s*a_i=a_i \ \ \hbox{for} \ \ i=1,\dots,g \ \ \ {\rm by \ (J).}$\par\medskip
Let \ $w_0=a_ge_{g-1}a_{g-1}e_{g-2}\dots e_1a_1b_1$.
\begin{equation}\label{A2}
w_0^{-1}*b_2=d_{1,2}, \ w_0^{-1}*b_1=a_1, \ w_0^{-1}*a_i=e_i, 
\ w_0^{-1}*e_i=a_{i+1}, \
\end{equation}
$d_{1,2}*b_1=b_1^{-1}*d_{1,2}, \ d_{1,2}*e_2=e_2^{-1}*d_{1,2}, \ 
[d_{1,2},a_i]=1, \ [d_{1,2},e_j]=1, \ \hbox{for} \ j\neq 2,$ 
$\strut[d_{1,2},t_j]=1 \ \hbox{for} \ \ j\neq 2.$
\proof[Proof of \rm(5)] We have $w_0^{-1}*b_2=$
 (by (M1)) $(b_1^{-1}a_1^{-1}e_1^{-1}a_2^{-1})*b_2=
d_{1,2}$. Other results of
conjugation by $w_0$ follow by (J).
Now\nl
 $\strut d_{1,2}*b_1=(b_1^{-1}a_1^{-1}e_1^{-1}a_2^{-1}b_2
a_2e_1a_1b_1)*b_1=$ (by J)\nl
$\strut(b_1^{-1}a_1^{-1}e_1^{-1}a_2^{-1}b_1^{-1}a_1^{-1}e_1^{-1}a_2^{-1})
*b_2=$ (by jumping from left side to the right)\nl
 $\strut(b_1^{-1}b_1^{-1}a_1^{-1}e_1^{-1}a_2^{-1})*b_2 = b_1^{-1}*d_{1,2}$.
 \par
Other relations follow from (M1) by conjugation by $w_0^{-1}$.\qed
\begin{equation}\label{A3}
a_k*d_{i,j}=d_{i,j} \ \ \hbox{for all}\ \ i,j,k, \
{\rm by \ (J)\ and \ \ (\ref{A1}) \ and
 \ \ (\ref{A2}).}
 \end{equation}
\begin{equation}\label{ti,tj}
t_i*t_{i+1}=t_{i+1}^{-1}*t_i \ for\ i=1,2,\dots,g-2, \
\end{equation}
(by the calculations similar to the proof of  Lemma
  \ref{basic relations}  (iv)),\nl
$\strut[t_i,t_k]=1\ \hbox{for} \ |i-k|>1, \ [t_i,s]=1 \ \hbox{for} \  i>1, \ {\rm by\ (M1).}$
\par\noindent
Using relations (\ref{A2}) and (\ref{ti,tj}) we can write the elements
$d_{i,j}$ in a different way.
\begin{equation}\label{new dij}
d_{i,i+1}=(t_{i-1}t_it_{i-2}t_{i-1}\dots t_1t_2)*d_{1,2} \ \hbox{for} \ i>0,
\end{equation}
$ d_{-i-1,-i}=(t_{i-1}^{-1}t_{i}^{-1}t_{i-2}^{-1}t_{i-1}^{-1}\dots
 t_1^{-1}t_2^{-1})*d_{-2,-1} \ \hbox{fo}r \ i>0,$\nl
$\strut d_{i,i+1}=(t_{i-1}t_i)*d_{i-1,i}, \ \ t_k*d_{i,i+1}=d_{i,i+1}
\ \ \hbox{for}\ \ |k-i|\neq 1. $
\par\medskip
Let $w_1=a_2e_1a_1b_1^2a_1e_1a_2$.
\begin{equation}\label{A5}
(b_1a_1e_1)^4=b_2w_1b_2w_1^{-1}, \ w_1*b_2=w_1^{-1}*b_2, \ [w_1*b_2,b_2]=1
\end{equation}
\proof[Proof of \rm(9)] We have $b_2w_1b_2=$ (by (M2)) $(b_1a_1e_1a_2)^5=$
(as in the proof of\nl Lemma \ref{basic relations} (v))
$\strut(b_1a_1e_1)^4w_1=$ (by (J))
$w_1(b_1a_1e_1)^4$.\nl Also $\strut b_2(b_1a_1e_1)^4=(b_1a_1e_1)^4b_2$, by (M1).
\par\noindent
Therefore $w_1b_2w_1^{-1}=b_2^{-1}(b_1a_1e_1)^4=(b_1a_1e_1)^4b_2^{-1}=
w_1^{-1}b_2w_1$ \ commutes with $b_2$.\qed

\begin{equation}\label{A6}
st_1s=b_1a_1e_1a_2^2e_1a_1b_1t_1=t_1b_1a_1e_1a_2^2e_1a_1b_1
\ \ \hbox{hence}\ \ st_1st_1=t_1st_1s\end{equation}
\proof[Proof of \rm(10)] We have a sequence of transformations by (J).\nl
 $\strut st_1s=$\nl
$\strut b_1a_1a_1b_1e_1a_1a_2e_1b_1a_1a_1b_1=b_1a_1a_1e_1b_1a_1b_1a_2e_1a_1a_1b_1=$
\nl
$\strut b_1a_1a_1e_1a_1b_1a_2a_1e_1a_1a_1b_1=b_1a_1e_1a_1e_1b_1a_2e_1e_1a_1e_1b_1=$
\nl
$\strut b_1a_1e_1a_1b_1a_2e_1a_2e_1a_1b_1e_1=b_1a_1e_1a_2a_1b_1a_2e_1a_1b_1a_2e_1=$
\nl
$\strut b_1a_1e_1a_2a_2e_1a_1b_1t_1.$\par
The second equality follows immediately by symmetry.\endpf
\noindent
Let $w_2=e_2a_2e_1a_1^2e_1a_2e_2$.
\begin{equation}\label{A7} (a_1e_1a_2)^4=
t_1^2a_1^2a_2^2=d_{1,2}d_{-2,-1}, \ \ [d_{1,2},d_{-2,-1}]=1.
\end{equation}
$d_{-2,-1}=w_2*d_{1,2}=w_2^{-1}*d_{1,2}=(b_1a_1e_1a_2)*b_2 $

\proof[Proof of \rm(11)]
We have $d_{-2,-1}=(s^{-1}t_1^{-1}s^{-1})*d_{1,2}=$ (by (\ref{A6}))
\nl
$\strut ((b_1a_1e_1a_2a_2e_1a_1b_1)^{-1}t_1^{-1}
b_1^{-1}a_1^{-1}e_1^{-1}a_2^{-1})*b_2=$
(by (J))\nl
 $\strut ((b_1a_1e_1a_2a_2e_1a_1b_1)^{-1}b_1^{-1}a_1^{-1}e_1^{-1}a_2^{-1})*b_2=$
(by (\ref{A5})) \ 
$(b_1a_1e_1a_2)*b_2$.  \par \noindent Conjugating (\ref{A5}) by 
$(b_1^{-1}a_1^{-1}e_1^{-1}a_2^{-1})$ we get $(a_1e_1a_2)^4=d_{1,2}d_{-2,-1}$.
Also, by (M1),\par\noindent
 $t_1^2a_1^2a_2^2=(a_1e_1a_2)^4$. This proves the first 
relation. The second relation follows from it by (\ref{A2}). Conjugating
(\ref{A5}) by $w_0^{-1}$ we get, by (\ref{A2}),
$(a_1e_1a_2)^4=d_{1,2}w_2d_{1,2}w_2^{-1}$ \ and 
\ $w_2*d_{1,2}=w_2^{-1}*d_{1,2}$.
 Therefore, from the first relation,
$d_{-2,-1}=w_2*d_{1,2}=w_2^{-1}*d_{1,2}=(b_1a_1e_1a_2)*b_2 $.\endpf
\begin{definition}{\rm If $A$ is a product of the generators we 
denote by $A^\prime$ the element obtained from $A$ by replacing each
generator by its inverse. We call $A^\prime$ the element {\em symmetric to}
 $A$. }\end{definition}
 \remark\label{symmetry}{\rm Relations (M1) and (M2) are symmetric.
  They remain valid
 when we replace each generator by its inverse. Therefore every relation
 between some elements of $\cM_{g,1}$ (products of generators) which follows from
 (M1) and (M2) remains valid if we replace each element by the element
 symmetric to it.}

\begin{equation}\label{A8}
\hbox{For} \ i+ j\neq 0 \ \ d_{i,j} \ \ \hbox{is symmetric to} \ d_{-j,-i}^{-1}.
\end{equation}
\proof[Proof of \rm(12)] Element $d_{1,2}$ is symmetric to $d_{-2,-1}^{-1}$, by (\ref{A7}).
Also $t_1^\prime=t_1^{-1}$ and $s^\prime = s^{-1}$.
 We see immediately that 
$d_{i,j}^\prime=d_{-j,-i}^{-1}$ for $i>0$. If $i<0$ and $i+j>0$ then
\par\noindent
$d_{i,j}=(t_{-i-1}^{-1}\dots t_1^{-1}s^{-1}t_{j-1}\dots t_2)*d_{1,2}$\nl and
$\strut d_{-j,-i}^{-1}=(t_{j-1}^{-1}\dots t_1^{-1}s^{-1}t_{-i}
\dots t_2)*d_{1,2}^{-1}$.\par\noindent
 Jumping with the positive powers of $t_k$ to the 
left we get 
\par\noindent
$d_{-j,-i}^{-1}=(t_{-i-1}\dots t_1st_{j-1}^{-1}\dots t_2^{-1}s^{-1}t_1^{-1}s^{-1})*
d_{1,2}^{-1}=$\nl$\strut(t_{-i-1}\dots t_1st_{j-1}^{-1}\dots t_2^{-1})*d_{1,2}^\prime=
d_{i,j}^\prime$.\qed

\begin{equation}\label{A9}
d_{i+1,i+2}=(t_i^{-1}t_{i+1}^{-1})*d_{i,i+1}=
(t_it_{i+1})*d_{i,i+1}\ \ \hbox{for}\ \ i=1,\dots,g-2.
\end{equation}
\proof[Proof of \rm(13)]
 For $i=1$ we have $(t_2t_1^2t_2)*d_{1,2}=$\ (by \ref{A7})\nl 
$\strut(e_2a_2a_3e_2e_1a_1a_2e_1e_1a_1a_2e_1e_2a_2a_3
e_2e_2^{-1}a_2^{-1}e_1^{-1}a_1^{-2}
e_1^{-1}a_2^{-1}e_2^{-1})*d_{-2,-1}$\nl$=$ \  (by (J)) 
$\strut (e_2a_2a_3e_2e_1a_1a_2e_1e_1a_2e_2a_3a_1^{-1}e_1^{-1}a_2^{-1}e_2^{-1})
*d_{-2,-1}=$ (by (J)) \nl
 $\strut(e_2a_2a_3e_2e_1a_1a_2e_1e_1a_1^{-1}a_2e_1^{-1}e_2a_2^{-1}a_3e_2^{-1})
*d_{-2,-1}=$ (by (J)) \nl
 $\strut(e_2a_2a_3e_2e_1a_2e_1^{-1}a_1a_1e_1a_2e_1^{-1}e_2a_2^{-1}a_3e_2^{-1})
*d_{-2,-1}=$ (by (J))  \nl
 $\strut(a_3^{-1}e_2a_2e_1a_1^2e_1a_2e_2a_3)*d_{-2,-1}=$ (by (\ref{A3}) and
 (\ref{A7})) \  $a_3^{-1}*d_{1,2}=$ (by (\ref{A3})) \  $d_{1,2}$. \\
 So
 (\ref{A9})
  is true for $i=1$. We continue by induction.
 Conjugating relation (\ref{A9}) by
 $t_it_{i+1}t_{i+2}$ we get, by (\ref{new dij}) and (\ref{ti,tj}), 
 relation (\ref{A9})
 for index $i+1$. \qed

\begin{equation}\label{A10}
t_i^2=d_{i,i+1}d_{-i-1,-i}a_i^{-2}a_{i+1}^{-2}.
\end{equation}
\proof[Proof of \rm(14)] The relation is true for $i=1$, by (\ref{A7}). 
We proceed by induction.\\
$t_{i+1}^2=(t_i^{-1}t_{i+1}^{-1})*t_i^2=(t_i^{-1}t_{i+1}^{-1})*
(d_{i,i+1}d_{-i-1,-i}a_i^{-2}a_{i+1}^{-2})$\nl$=
$ (by (\ref{ti,tj}), (\ref{A9}),
 (\ref{new dij})
   and (\ref{A1}))
$\strut d_{i+1,i+2}d_{-i-2,-i-1}a_{i+1}^{-2}a_{i+2}^{-2}$.\qed

\begin{equation}\label{A11}
d_{-i,i} \ \ \hbox{is symmetric to} \ \ d_{-i,i}^{-1}.
\end{equation}
\proof[Proof of \rm(15)] The relation is true for $i=1$. The general case follows from
(\ref{A10}), (\ref{A1}), (\ref{A2}), (\ref{A8}) and the definitions.\qed\par
\begin{equation} \label{b_1,d22}
[b_1,d_{-2,2}]=1
\end{equation}
 \proof[Proof of \rm(16)] By the definition $d_{-2,2}=(t_1^{-1}d_{1,2})*(s^2a_1^4)=$ 
(by (\ref{A1}))\nl$\strut a_2^4((d_{1,2}t_1^{-1})*s^2)$.\nl Now
$\strut t_1^{-1}*s=t_1^{-1}st_1ss^{-1}=$ \ (by (\ref{A6})) \ 
$ b_1a_1e_1a_2^2e_1a_1^{-1}b_1^{-1}$.\par\noindent
 Taking squares we get $t_1^{-1}*s^2=b_1a_1e_1a_2^2e_1^2a_2^2e_1a_1^{-1}
 b_1^{-1}$ and\nl
 $\strut d_{-2,2}=a_2^4d_{1,2}b_1a_1e_1a_2^2e_1^2a_2^2e_1a_1^{-1}
 b_1^{-1}d_{1,2}^{-1}$. Now $b_1$ commutes with $d_{-2,2}$, by 
 (\ref{A1}), $\strut$(\ref{A2}) and (J).\qed\par
 \begin{equation}\label{A12}
(t_it_{i+1})*e_i=e_{i+1}\ 
\ \hbox{for}\ \ i=1,\dots,g-2.
\end{equation}
\proof[Proof of \rm(17)] We have, by (J),  $(t_it_{i+1})*e_i=$\nl
$\strut(e_ia_ia_{i+1}e_ie_{i+1}
a_{i+1}a_{i+2}e_{i+1})*e_i=$
$(e_ia_ia_{i+1}e_{i+1}e_i
a_{i+1})*e_i=e_{i+1}$.\qed\par
\begin{equation}\label{A13}
[b_1,d_{i,i+1}]=1 \ if \ i>1,
\end{equation}
$[e_k,d_{i,i+1}]=1 \ \ \hbox{if}\ \ |k-i|\neq 1, i>0,\ \ 
d_{i,i+1}*e_k=e_k^{-1}*d_{i,i+1} \ \ \hbox{if} \ \ |k-i|=1.$
\proof[Proof of \rm(18)]  By (\ref{A9}) and  (J)\nl 
 $\strut d_{2,3}=(e_1^{-1}a_2^{-1}e_2^{-1}a_3^{-1}a_1^{-1}e_1^{-1}a_2^{-1}
 e_2^{-1}b_1^{-1}a_1^{-1}e_1^{-1}a_2^{-1})*b_2=$\nl
 $\strut(e_1^{-1}a_2^{-1}e_2^{-1}a_3^{-1}a_1^{-1}b_1^{-1}e_1^{-1}a_1^{-1}a_2^{-1}
 e_1^{-1}e_2^{-1}a_2^{-1})*b_2$.\par Now $b_1$ commutes with
 $d_{2,3}$, by (J). For $i>2$ we have $b_1$ commutes with $d_{i,i+1}$
 by (M1) and  (\ref{new dij}). Conjugating by $a_1b_1t_1$
  we get, by (\ref {ti,tj}) and (J), $[e_1,d_{i,i+1}]=1$, \ for $i>2$.
 We also have $d_{1,2}*b_1=b_1^{-1}
 *d_{1,2}$, by (\ref{A2}). We conjugate this equality by
 $u=a_1b_1t_1t_2$
 and get   $d_{2,3}*e_1=e_1^{-1}*d_{2,3}$, by (\ref{new dij}) and  
 the first part of the proof. Conjugating relations (\ref{A2}) and 
 the above relations by suitable  products
 $t_it_{i+1}$ we get all remaining relations, by (\ref{A12}).\qed\par
\begin{equation}\label{A14}
[b_1,d_{1,2}sd_{1,2}]=1, \ \ \hbox{hence} \ \ [s,d_{1,2}sd_{1,2}]=1,
\end{equation}
$[e_j,d_{i,i+1}t_jd_{i,i+1}]=1, \ \ \hbox{hence} \ \
[t_j,d_{i,i+1}t_jd_{i,i+1}]=1 , \ \ \hbox{if} \ \ |i-j|=1.$
\proof[Proof of \rm(19)] By (\ref{A2}) we have $(d_{1,2}sd_{1,2})*b_1=
(d_{1,2}b_1a_1a_1b_1b_1^{-1})*d_{1,2}=(d_{1,2}b_1)*d_{1,2}=b_1$.
The other case is similar, but we use (\ref{A13}) instead of (\ref{A2}).\endpf
\par
\begin{remark}{\rm  Relation $uvuv=vuvu$ \ implies 
\ $(uv)*u=v^{-1}*u$ and $(u^{-1}v^{-1})*u=v*u$. Relations (\ref{A14}) will be 
often used in this form.}\end{remark}
Observe that relations (\ref{A1}) -- (\ref{A10}) imply in particular that
relations (P1) and relations (P3) -- (P7) follow from (M1) and (M2). We shall 
prove now that relations (P8) follow from (M1) and (M2).

\begin{definition}\label{proper move} {\rm We say that homeomorphism $t_k$ 
(respectively $t_k^{-1}d_{k+1,k+2}$ or $s$) moves a curve $\delta_{i,j}$ 
properly if
it takes it to some curve $\delta_{p,q}$.}\end{definition}
\begin{lemma}\label{t-conjugation} If $t_k$ (respectively $s$) moves curve
$\delta_{i,j}$ to some curve $\delta_{p,q}$ then $t_k*d_{i,j}=d_{p,q}$
(respectively $s*d_{i,j}=d_{p,q}$).
\end{lemma}
\begin{remark}{\rm Since the action of $t_k$ and $s$ on a curve 
$\delta_{i,j}$ is described by Lemma \ref{action of tk} and is easily 
determined, Lemma \ref{t-conjugation} helps us to understand the 
result of the conjugation. In fact the action corresponds exactly to 
relations (P8), so we have to prove that relations (P8) follow from 
(M1) and (M2).}\end{remark}\noindent
{\bf Proof of Lemma \ref{t-conjugation}}\qua We know that $[t_i,d_{i,i+1}]=1$,
by (\ref{new dij}), (\ref{ti,tj}), and (\ref{A2}). \par\noindent
If $i<0$ and $i+j=1$ then $t_{j-1}^{-1}*d_{i,j}=
(t_{j-1}^{-1}t_{j-2}^{-1}\dots t_1^{-1}s^{-1}t_{j-1}\dots t_2)*d_{1,2}=
d_{i-1,j-1}$.  \par\noindent
For $i>0$ or $i<0, \ i+j>0$ all other cases of conjugation by $t_k$
follow from (\ref{A2}), (\ref{ti,tj}) and the definitions. The other cases of 
$i\neq -j$ follow by symmetry. \\
Consider conjugation by $s$ for $i\neq -j$. Again it suffices to consider
$i>0$ or $i<0, \ i+j>0$. The other cases follow by symmetry.
We have $s^{-1}*d_{1,j}=d_{-1,j}$.\\
If $i>1$ then $d_{i,j}=$ (by (\ref{ti,tj})) \  
$(t_{i-1}\dots t_2t_{j-1}\dots t_3)
*d_{2,3}$ and $s*d_{i,j}=d_{i,j}$, by (\ref{ti,tj}), (\ref{A3})
 and (\ref{A13}).
 \par\noindent
  If
$i<1$ then $d_{i,j}=$ (by \ref{ti,tj}) \ 
 $(t_{-i-1}^{-1}\dots t_2^{-1}t_{j-1}\dots
t_3)*d_{-2,3}$. Also $s^{-1}*d_{-2,3}=$ \par\noindent
$(s^{-1}t_1^{-1}s^{-1}t_2)*d_{1,2}=
(s^{-1}t_1^{-1}s^{-1}t_1^{-1})*d_{2,3}=$
 (by (\ref{A6})) \  
$(t_1^{-1}s^{-1}t_1^{-1}s^{-1})*d_{2,3}= d_{-2,3}$
(by (\ref{A3}) and (\ref{A13})).
So $s^{-1}*d_{i,j}=d_{i,j}$, by (\ref{ti,tj}). \\
We now consider conjugation of $d_{-j,j}$. Clearly $s*d_{-1,1}=d_{-1,1}$,
by (M1). Also $s*d_{-2,2}=d_{-2,2}$, by (\ref{A3}) and (\ref{b_1,d22}). 
For $i>1$
we have $[s,t_i]=1$, by (M1), hence $s*d_{-j,j}=d_{-j,j}$ for all $j$,
by the first part of the proof.\\
 Consider conjugation by $t_k$.\\
  For $k>j$
we have $t_k*d_{-j,j}=d_{-j,j}$, by (\ref{ti,tj}) and the first part.\\
Curves $t_j(\delta_{-j,j})$ and $t_{j-1}(\delta_{-j,j})$ are not of 
the form $\delta_{p,q}$.\\
 Consider $k=j-2$ (the other cases follow by conjugation and
by the first part of the proof). We have \nl
$\strut d_{-j,j}=(d_{j-1,j}t_{j-1}^{-1}t_{j-2}^{-1}d_{j-2,j-1})*d_{2-j,j-2}=$
(by (\ref{new dij}) and (\ref{A9})) \nl 
$\strut (t_{j-2}^{-1}t_{j-1}^{-1}d_{j-2,j-1}t_{j-1}t_{j-2}t_{j-1}
^{-1}t_{j-2}^{-1}d_{j-2,j-1})*d_{2-j,j-2}=$ (by (\ref{ti,tj}) and 
(\ref{new dij}))  \nl
$\strut (t_{j-2}^{-1}t_{j-1}^{-1}t_{j-2}^{-1}d_{j-2,j-1}t_{j-1}
d_{j-2,j-1})*d_{2-j,j-2}$. Now, by (\ref{ti,tj})  and (\ref{A14}),
 \nl
$\strut t_{j-2}^{-1}*d_{-j,j}=
(t_{j-2}^{-1}t_{j-1}^{-1}t_{j-2}^{-1}d_{j-2,j-1}t_{j-1}
d_{j-2,j-1}t_{j-1}^{-1})*d_{2-j,j-2}=d_{-j,j}$
(by the previous case).\endpf

We now pass to the biggest task of this section: the relations (P2).

\begin{equation}\label{A15}
d_{i,j} \ \hbox{commutes with} \ d_{-1,1} \ \hbox{if} \  i,j\neq \pm1,
\end{equation}
$ d_{i,j} \ \hbox{commutes with} \ d_{k,k+1}$ if all indices are  
distinct.

\proof[Proof of \rm(20)]
We know by Lemma \ref{t-conjugation} that $d_{i,j}$ commutes with $s$,
hence also with $d_{-1,1}$. Consider the other cases. We assume first that 
$k>0$.
 Consider curves $\delta_{i,j}$ and $\delta_{k,k+1}$.
We want to move the curves properly to some standard position 
by application of products of
$s$ and $t_m$'s. This moves holes (see definition \ref{half-twist})
and corresponds to conjugation of $d_{i,j}$ and $d_{k,k+1}$,
by the previous lemma.\par\noindent
{\bf Case 1}\qua $i\neq -j$\qua Observe that for every
$k$  either $t_{k+1}t_{k}$ or $t_{k+1}^{-1}t_k^{-1}$ 
moves $\delta_{i,j}$ properly. Both products take 
$\delta_{k+1,k+2}$ onto $\delta_{k,k+1}$ and conjugation of 
$d_{k+1,k+2}$ by either product produces $d_{k,k+1}$.
 If either $|i|$ 
or $|j|$ is bigger than $k+1$ we may assume, moving 
 $\delta_{i,j}$ properly, and not moving $\delta
_{k,k+1}$, that either $|i|$ or $|j|$ is equal to $k+2$. Then 
either $t_kt_{k+1}$ or $t_k^{-1}t_{k+1}^{-1}$ moves $\delta_{i,j}$
properly and moves $\delta_{k,k+1}$ to $\delta_{k+1,k+2}$ and leaves at
most one index $|i|$ or $|j|$ bigger than $k+1$. Applying this procedure 
again, if necessary, we may assume $j<k$. If $j>0$ we can move $\delta_{i,j}$
properly, without moving $\delta_{k,k+1}$, and get a curve $\delta_{i,j}$
with $j<0$. Now applying consecutive products 
$t_{m+1}t_{m}$ or $t_{m+1}^{-1}t_m^{-1}$ we reach $k=1$. Further moves
will produce one of the following three cases:\par\medskip\noindent
{\bf Case 1a}\qua $d_{1,2}$ commutes with $d_{-2,-1}$.\qua True, by (\ref{A7}).
\par\medskip\noindent
{\bf Case 1b}\qua $d_{1,2}$ commutes with $d_{-3,-1}$ or $d_{-3,-2}$.\qua We can 
 conjugate $d_{-3,-2}$ by $t_1$ and get
  $d_{-3,-1}$.
Now $d_{-3,-1}=t_2^{-1}*d_{-2,-1}=$ (by (\ref{A7})) \  $(e_2^{-1}a_3^{-1}
a_2^{-1}e_2^{-1}w_2)*d_{1,2}=$ (by (J) and
 (\ref{A2}))
$(e_2^{-1}a_3^{-1}e_1a_1d_{1,2}^{-1}e_2^{-1}a_2^{-1}e_1^{-1})*a_1=$
(by (J) and (\ref{A2}))\nl
$\strut(e_2^{-1}d_{1,2}^{-1}a_3^{-1}e_2^{-1}e_1a_1a_2^{-1}
e_1^{-1})*a_1$.\par\noindent
 The last expression commutes with $d_{1,2}$, by (J) and
(\ref{A2}).
\par\medskip\noindent
{\bf Case 1c}\qua $d_{-2,-1}$ commutes with $d_{3,4}$.\qua The proof is 
rather long. We consider
relation (M3):
  $d_3=a_3^{-1}a_2^{-1}a_1^{-1}
d_{1,2}d_{1,3}d_{2,3}$, where\nl
 $\strut d_3=(b_1^{-1}b_2a_2e_1e_2a_2a_3e_2b_1^{-1}a_1^{-1}e_1^{-1}a_2^{-1})
*b_2=$
\nl$\strut(b_1^{-1}b_2a_2e_1e_2a_2a_3e_2b_2a_2e_1a_1)*b_1$.
It follows, by (J), that $d_3$ commutes with \ $a_2e_1e_2a_2a_3e_2$,
hence \ $d_3=(b_1^{-1}((a_2e_1e_2a_2a_3e_2)^{-1}*b_2)
b_1^{-1}a_1^{-1}e_1^{-1}a_2^{-1})*b_2$. We now conjugate relation
 (M3) by
$u=(a_4e_3a_3e_2a_2e_1a_1b_1)^{-1}a_3e_2a_2e_1a_1b_1$. When we write
$d_{1,2}=(b_1^{-1}a_1^{-1}e_1^{-1}a_2^{-1})*b_2$ we see that $d_{1,2}$
commutes with $u$, by (J). All other factors on the RHS commute with 
$u$ by (J), so we may replace $d_3$ in the relation (M3)
by $u*d_3=$\nl 
$\strut(a_4e_3a_3e_2a_2e_1a_1b_1)^{-1}*
((a_3e_2a_2e_1a_1((a_2e_1e_2a_2a_3e_2)^{-1}*b_2)
b_1^{-1}a_1^{-1}e_1^{-1}a_2^{-1})*b_2)$.\par\noindent 
 We conjugate each term by 
$(a_4e_3a_3e_2a_2e_1a_1b_1)^{-1}$ and get\nl
$\strut u*d_3=(e_3a_3e_2a_2e_1((e_2a_2a_3e_2e_3a_3)^{-1}*d_{1,2})
a_1^{-1}e_1^{-1}a_2^{-1}e_2^{-1})*d_{1,2}$. \par\noindent 
We now conjugate each term of 
relation (M3) by $t_2^{-1}t_3^{-1}=$\nl
$\strut e_2^{-1}a_3^{-1}e_3^{-1}a_4^{-1}
a_2^{-1}e_2^{-1}a_3^{-1}e_3^{-1}$.\par\noindent 
 The RHS becomes 
$a_4^{-1}a_3^{-1}a_1^{-1}
(t_2^{-1}*d_{1,2})
((t_2^{-1}t_3^{-1}t_2)*d_{1,2})d_{3,4}$. \par\noindent 
The LHS becomes\nl
$\strut(e_2^{-1}a_3^{-1}e_3^{-1}a_4^{-1})*
((e_1((e_2a_2a_3e_2e_3a_3)^{-1}*d_{1,2})
(a_1^{-1}e_1^{-1}a_2^{-1}e_2^{-1}))*d_{1,2})$. \par\noindent 
We conjugate each bracket.\nl
$\strut(e_2^{-1}a_3^{-1}e_3^{-1}a_4^{-1})*e_1=e_1$.\nl
$\strut(e_2^{-1}a_3^{-1}e_3^{-1}a_4^{-1}(e_2a_2a_3e_2e_3a_3)^{-1})*d_{1,2}
=$ (by (J) and (\ref{A2})) \ $(t_3^{-1}t_2^{-1})*d_{1,2}$.
\nl
$\strut(e_2^{-1}a_3^{-1}e_3^{-1}a_4^{-1}a_1^{-1}e_1^{-1}a_2^{-1}e_2^{-1})
*d_{1,2}=$
(by (J) and (\ref{A2})) \ $(a_1^{-1}e_1^{-1}t_2^{-1})*d_{1,2}$.\par\noindent 
We shall prove that all terms of the obtained equation commute with 
$d_{-2,-1}$,
except possibly for $d_{3,4}$. Therefore $d_{3,4}$ also commute with
$d_{-2,-1}$. Clearly $a_1,e_1,a_2,a_3,a_4,t_1,t_3$ commute with $d_{-2,-1}$
by (\ref{A2}) and symmetry. Also $d_{1,2}$ commutes with $d_{-2,-1}$ by 
 Case 1a and $d_{1,3}$ and $d_{2,3}$ commute with $d_{-2,-1}$ by 
 Case 1b. Therefore $t_2^{-1}*d_{1,2}=t_1*d_{2,3}$ commute with
$d_{-2,-1}$ and $(t_2^{-1}t_3^{-1}t_2)*d_{1,2}=$ (by (\ref{ti,tj}) and
(\ref{A2})) \ $(t_3t_2^{-1})*d_{1,2}$ \ commutes with $d_{-2,-1}$. It follows that 
$d_{3,4}$ commutes with $d_{-2,-1}$\par\medskip\noindent 
{\bf Case 2}\qua $i=-j$\qua
 We have to prove that $d_{k,k+1}$ commutes with
$d_{-j,j}$ if $j\neq k,k+1$. If $j<k$ then the result follows by Lemma
\ref{t-conjugation} and by Case 1. If $j>k+2$ we can properly move 
$\delta_{k,k+1}$ to $\delta_{j-2,j-1}$, without moving $\delta_{-j,j}$.
So we may assume $j=k+2$.
 We have \nl
 $\strut d_{-j,j}=(d_{k+1,k+2}t_{k+1}^{-1}
t_k^{-1}d_{k,k+1})*d_{-k,k}=$ (by (\ref{A3}) and (\ref{A10}))\nl
$\strut(d_{-k-2,-k-1}^{-1}t_{k+1}t_kd_{-k-1,-k}^{-1})*$
 $d_{-k,k}$. \par\noindent 
We conjugate $d_{k,k+1}$ by
$(d_{-k-2,-k-1}^{-1}t_{k+1}t_kd_{-k-1,-k}^{-1})^{-1}$
 and get, by Case 1 and (\ref{A9}),
$(d_{-1-k,-k}t_{k}^{-1}t_{k+1}^{-1}d_{-k-2,-k-1})*d_{k,k+1}=
d_{-1-k,-k}*d_{k+1,k+2}=
d_{k+1,k+2}$, which commutes with $d_{-k,k}$ by the first part of Case 2.
 
 Suppose now that $k<0$. Cases 1a and 1b follow by symmetry. 
 In the Case 1c we arrive by symmetry at the situation where we have 
 to prove that 
 $d_{1,2}$ commutes with $d_{-4,-3}$. Conjugating by $t_2t_1t_3t_2$ we
 get a pair $d_{3,4}$, $d_{-2,-1}$, which commutes by Case 1c.
 Now Case 2 follows by symmetry.
 \endpf
\begin{lemma} \label{td-conjugation}
 If $t_k^{-1}d_{k,k+1}$ takes
a curve $\delta_{i,j}$ to $\delta_{p,q}$ then $t_k^{-1}d_{k,k+1}*d_{i,j}=
d_{p,q}$.\end{lemma}
\proof We shall list the relevant cases. If $i$ or $j$ is equal to
$k$ it becomes $k+1$, if  $i$ or $j$ is equal to
$-k$ it becomes $-k-1$. In particular 
$t_k^{-1}d_{k,k+1}(\delta_{-k,k})=\delta_{-k-1,k+1}$.
Indices $k+1$ and $-k-1$ are forbidden, 
they do not move properly, except for $\delta_{k,k+1}$ 
and $\delta_{-k-1,-k}$
which are fixed by $t_k^{-1}d_{k,k+1}$. Other indices $i$, $j$ do not change. \\
We now pass to the proof of the Lemma. 
If $i=-j$ the result follows from the definitions and from 
(\ref{A15}) and Lemma \ref{t-conjugation}. 
Suppose $i\neq -j$. If $i$ and $j$ are different from $k$ 
 then $d_{i,j}$ commutes with $d_{k,k+1}$, by (\ref{A15}), 
and $t_k^{-1}$ moves $\delta_{i,j}$ properly, so we are done by Lemma
\ref{t-conjugation}. If $i$ and $j$ are different from $-k$ 
we can replace $t_k^{-1}d_{k,k+1}$ by 
$t_kd_{-k-1,-k}^{-1}a_k^2a_{k+1}^2$, using (\ref{A10}), and we are
done by a similar argument. \endpf\par
Lemmas \ref{t-conjugation} and \ref{td-conjugation} allow us to 
reduce relations
(P2) to relatively small number of cases. We can apply product of half-twists
$h_{k,k+1}$ and $h_{-k-1,-k}$ either in the same direction, conjugating by
$t_k$,  or in opposite
directions, conjugating by $t_k^{-1}d_{k,k+1}$,  and move properly curves 
corresponding to elements $d_{p,q}$
in the relations (P2) into a small number of standard configurations.

\begin{equation}\label{A16}
d_{i,j}\ \ \hbox{commutes with}\ \ d_{r,s}\ \ \hbox{if}\ 
\ r<s<i<j\ \ \hbox{or}\ \ i<r<s<j.
\end{equation}

\proof[Proof of \rm(21)]
Moving curves $\delta_{i,j}$ and $\delta_{r,s}$ properly we
can arrive at a situation  $s=r+1$ or $j=i+1$ or $-r=s=1$
or $-i=j=1$. In particular if 
$i<r<s<j$ and $r=-s$ then conjugating by 
$t_2d_{2,3}^{-1}\dots t_{s-1}d_{1-s,s-1}^{-1}$ we get a pair $d_{-1,1}$,
$d_{i,j}$. 
Then (\ref{A16}) follows from (\ref{A15}).\qed

\begin{equation}\label{A17}
d_{r,i}^{-1}*d_{i,j}=d_{r,j}*d_{i,j}\ \ \hbox{if}\ \ r<i<j.
\end{equation}

\proof[Proof of \rm(22)] 
Indices $(i,-i)$ move together (either remain fixed or
move to $(i-1,1-i)$ or to $(i+1,-i-1)$) when we conjugate by $t_k$ or 
$s$ or $t_k^{-1}d_{k,k+1}$. Moving curves properly we can arrive at one 
of the following four cases depending on the  pairs of opposite indices.
\par\medskip\noindent
{\bf Case 1}\qua There is no pair of opposite numbers among $r,i,j$.\qua We may
assume that $(r,i,j)=(1,2,3)$.\\
 $d_{1,2}^{-1}*d_{2,3}=$ (by (\ref{A9})) 
 $(d_{1,2}^{-1}t_1^{-1}t_2^{-1})*d_{1,2}=$ (by (\ref{A2}) and (\ref{A14})) \
$(t_1^{-1}t_2)*d_{1,2}$.\nl
$\strut d_{1,3}*d_{2,3}=(t_2d_{1,2}t_2^{-1}t_1t_2)*d_{1,2}=$ (by (\ref{A2})
 and (\ref{ti,tj}))
$(t_2t_1d_{1,2}t_2)*d_{1,2}=$ (by (\ref{A14})) $(t_2t_1t_2^{-1})*d_{1,2}
=$ (by (\ref{A2}) and (\ref{ti,tj})) \ $(t_1^{-1}t_2)*d_{1,2}$.\par\medskip\noindent
{\bf Case 2}\qua $i=-r$\qua We may assume that $(r,i,j)=(-1,1,2)$.\\
$d_{-1,1}^{-1}*d_{1,2}=(s^{-2}a_1^{-4})*d_{1,2}=$ (by (\ref{A14})) \ 
$(s^{-1}d_{1,2}s)*d_{1,2}=d_{-1,2}*d_{1,2}$.\par\medskip\noindent
{\bf Case 3}\qua $r=-j$\qua We may assume that $(r,i,j)=(-2,1,2)$. We have to prove
that $d_{-2,1}^{-1}*d_{1,2}=d_{-2,2}*d_{1,2}$. We conjugte by
$d_{1,2}^{-1}t_1$. The right hand side becomes $s^2*d_{1,2}$ and the left hand 
side becomes $(d_{1,2}^{-1}t_1t_1^{-1}s^{-1}d_{1,2}^{-1}st_1)*d_{1,2}=$
(by (\ref{A2}) and (\ref{A14})) 
$(sd_{1,2}^{-1}s^{-1}d_{1,2}^{-1})*d_{1,2}=$ (by (\ref{A14})) $s^2*d_{1,2}$.
\par\medskip\noindent
{\bf Case 4}\qua $i=-j$\qua We may assume that $(r,i,j)=(-2,-1,1)$.\\
$d_{-2,-1}^{-1}*d_{-1,1}=a_1^4(d_{-2,-1}^{-1}*s^2)=$ (by (\ref{A14}) and 
symmetry) $a_1^4((sd_{-2,-1}s^{-1})*s^2)=d_{-2,1}*d_{-1,1}$.\qed

\begin{equation}\label{A18}
d_{r,j}^{-1}*d_{r,i}=d_{i,j}*d_{r,i} \ \ \hbox{if}\ \ r<i<j.
\end{equation}

\proof[Proof of \rm(23)]
Let us apply symmetry to relation (\ref{A17}). We get\nl
$\strut d_{-i,-r}*d_{-j,-i}^{-1}=d_{-j,-r}^{-1}*d_{-j,-i}^{-1}$ if $-j<-i<-r$.
This is relation (\ref{A18}) after a suitable change of indices.\qed

\begin{equation}\label{A19}
[d_{i,j},d_{r,j}^{-1}*d_{r,s}]=1 \ \ \hbox{if}\ \ r<i<s<j.
\end{equation}

\proof[Proof of \rm(24)] 
Again we have to consider different cases depending on pairs of 
opposite indices.
For each of them we move curves properly to some standard position. If we
apply symmetry we get $[d_{-j,-i},d_{-j,-r}*d_{-s,-r}]=1$. Now conjugate by
$d_{-j,-r}^{-1}$ and get $[d_{-s,-r},d_{-j,-r}^{-1}*d_{-j,-i}]=1$ 
if $-j<-s<-i<-r$. This is again relation (\ref{A19}) with  different
pairs of opposite indices. So $i=-j$ is equivalent to $r=-s$ and
$s=-j$ is equivalent to $r=-i$. We are left with the following five cases.
\par\medskip\noindent
{\bf Case 1}\qua  There is no pair of opposite numbers among $r,i,s,j$.\qua We may
assume that $(r,i,s,j)=(1,2,3,4)$.
 Conjugate by $t_3^{-1}$.\nl
  $\strut t_3^{-1}*d_{2,4}=d_{2,3}$.\nl
$\strut(t_3^{-1}d_{1,4}^{-1})*d_{1,3}=(t_2d_{1,2}^{-1}t_2^{-1}t_3^{-1}t_2)
*d_{1,2}=$
(by (\ref{ti,tj})) $(t_2d_{1,2}^{-1}t_3t_2^{-1}t_3^{-1})*d_{1,2}=$ (by (\ref{A2}))
$(t_2t_3d_{1,2}^{-1}t_2^{-1})*d_{1,2}=$ (by (\ref{A14})) $(t_2t_3t_2)*
d_{1,2}=d_{1,4}$, and it commutes with $d_{2,3}$, by (\ref{A16}).
\par\medskip\noindent
{\bf Case 2}\qua $r=-i$\qua  We may assume that $(r,i,s,j)=(-1,1,2,3)$. We 
conjugate by $t_2^{-1}$.\nl
$\strut t_2^{-1}*d_{1,3}=d_{1,2}$.\nl
$\strut(t_2^{-1}d_{-1,3}^{-1})*d_{-1,2}=(t_2^{-1}s^{-1}t_2d_{1,2}^{-1}t_2^{-1})
*d_{1,2}=$ (by (\ref{ti,tj}) and (\ref{A14})) $(s^{-1}t_2)*d_{1,2}=d_{-1,3}$,
and it commutes with $d_{1,2}$, by (\ref{A16}).\par\noindent
{\bf Case 3}\qua $r=-s$\qua We may assume that $(r,i,s,j)=(-2,1,2,3)$. By 
(\ref{A18}) we have $d_{-2,3}^{-1}*d_{-2,2}=d_{2,3}*d_{-2,2}$, so we have to
prove that $d_{1,3}$ commutes with $d_{2,3}*d_{-2,2}$. We conjugate by 
$t_2^{-1}$ and get\nl
$\strut t_2^{-1}*d_{1,3}=d_{1,2}$.\nl
$\strut(t_2^{-1}d_{2,3})*d_{-2,2}=d_{-3,3}$, and it commutes with $d_{1,2}$, by
 (\ref{A16}).\par\medskip\noindent
{\bf Case 4}\qua $r=-j$ and $i\neq -s$\qua  
We may assume that $(r,i,s,j)=(-3,1,2,3)$. By 
(\ref{A18}) we have $d_{-3,3}^{-1}*d_{-3,2}=d_{2,3}*d_{-3,2}$. After
conjugation by $t_1^{-1}d_{2,3}^{-1}$  we have to
prove that $t_1^{-1}*d_{-3,2}=d_{-3,1}$ commutes with 
$(t_1^{-1}d_{2,3}^{-1})*d_{1,3}=(t_1^{-1}d_{2,3}^{-1}t_1^{-1})*d_{2,3}=$
(by (\ref{A14})) \ $d_{2,3}$. This is true by (\ref{A16}).
\par\medskip\noindent
{\bf Case 5}\qua $i=-s$\qua  We may assume that $(r,i,s,j)=(-2,-1,1,m)$, where 
$m=2$ or $m=3$. By 
(\ref{A18}) we have $d_{-2,m}^{-1}*d_{-2,1}=d_{1,m}*d_{-2,1}$. After 
conjugation by $s^{-1}d_{1,m}^{-1}$ we have to consider 
$s^{-1}*d_{-2,1}=d_{-2,-1}$ and $(s^{-1}d_{1,m}^{-1})*d_{-1,m}$.
For $m=2$ the last expression is equal $(s^{-1}d_{1,2}^{-1}s^{-1})*d_{1,2}=
d_{1,2}$, by (\ref{A14}). For $m=3$ we get 
$(s^{-1}d_{1,3}^{-1}s^{-1})*d_{1,3}=
(s^{-1}t_2d_{1,2}^{-1}t_2^{-1}t_2s^{-1})*d_{1,2}=$ (by (\ref{A14})) \
 $d_{1,3}$.
Both elements commute with $d_{-2,-1}$, by (\ref{A16}).\endpf
\par
This concludes the proof of the fact 
that relations (P2) follow from (M1) -- (M3).
\par

We now pass to the relations (P9) -- (P11).\\
Consider first relations (P9). Clearly $a_1$ commutes with all the elements
in (P9), by (\ref{A1}) and (\ref{A3}), so it suffices to prove that $b_1$ 
commutes with these elements. It commutes with $a_1^2s$, $t_1st_1$,
$a_2$, and $t_i$, for $i>1$, by (J). Also $b_1$ commutes with 
$d_{-2,2}$, by (\ref{b_1,d22}), and commutes with $d_{2,3}$ by 
(\ref{A13}). Finally it commutes with $d_{1,2}sd_{1,2}$, by 
(\ref{A14}), hence also commutes with
$d_{-1,1}d_{-1,2}d_{1,2}a_1^{-2}a_2^{-1}=
a_1^4s^2s^{-1}d_{1,2}sd_{1,2}a_1^{-2}a_2^{-1}=a_1^2sd_{1,2}sd_{1,2}a_2^{-1}$.
This proves relations (P9). Relation (P10) follows from the
definitions and (M1).\par
We now pass to relations (P11).

For $i=1,2,3,4$ we have $g_i=(k_ir)^3$. By 
(\ref{A3})  $k_i*a_1=a_1$  therefore it suffices to prove that
 $k_i*b_1=b_1^{-1}*k_i$. Then\nl
 $\strut g_i=(k_ia_1b_1a_1)^3=$ (by (\ref{A3}))
 $k_ia_1b_1a_1^2k_ib_1k_ia_1^2b_1a_1=k_ia_1sk_isa_1$
and this is exactly relation (P11) for $i=1,2,3,4$.

\par\medskip
\begin{equation} k_1*b_1=b_1^{-1}*k_1. \end{equation}
\proof[Proof of \rm(25)]
Since $k_1=a_1$ the result follows from (M1).\qed

\begin{equation} k_2*b_1=b_1^{-1}*k_2. \end{equation}
\proof[Proof of \rm(26)]
Since $k_2=d_{1,2}$ the result follows from (\ref{A2}).\qed

\begin{equation} k_3*b_1=b_1^{-1}*k_3. \end{equation}
\proof[Proof of \rm(27)]
We have \nl
$\strut k_3=a_1^{-1}a_2^{-2}d_{1,2}d_{-2,1}d_{-2,2}=
a_1^{-1}a_2^{-2}d_{1,2}t_1^{-1}s^{-1}d_{1,2}st_1t_1^{-1}
d_{1,2}s^2a_1^4d_{1,2}^{-1}t_1=$\nl
(by (\ref{A2}) and (J)) \ $\strut a_1^{-1}t_1^{-1}d_{1,2}s^{-1}d_{1,2}sd_{1,2}s^2
a_1^2d_{1,2}^{-1}t_1=$ (by (\ref{A14})) \nl
$\strut a_1^{-1}t_1^{-1}d_{1,2}sd_{1,2}sa_1^2t_1=$ (by (\ref{A2}) and (J))
$a_1^{-1}t_1^{-1}(d_{1,2}b_1a_1)^4t_1$.\par\noindent
Let $u=t_1b_1d_{1,2}a_1b_1$. It follows from (\ref{A2}) and (J) that
$u*a_1=d_{1,2}$, $u*e_1=b_1$, $u*a_2=a_1$, $u*d_{1,2}=a_2$. Conjugating 
(\ref{A7}) by $u$ we get $(d_{1,2}b_1a_1)^4=a_2(u*d_{-2,-1})$, hence
$k_3=$ (by (\ref{A1})) $(t_1^{-1}u)*d_{-2,-1}=(b_1d_{1,2}a_1b_1)*d_{-2,-1}$.
We want to prove that $b_1$ and $k_3$ are braided. It suffices to prove it 
for their common conjugates. We conjugate by $b_1^{-1}$ and get $b_1$ and
$(d_{1,2}a_1b_1)*d_{-2,-1}$. Now we conjugate by 
$a_1^{-1}b_1^{-1}d_{1,2}^{-1}$ and get, by (\ref{A2}), (\ref{A3})
and symmetry, \
$d_{1,2}$ and $b_1*d_{-2,-1}=d_{-2,-1}^{-1}*b_1$.  Now we conjugate
by $d_{-2,-1}$ and get, by (\ref{A15}), \ $d_{1,2}$ and $b_1$, 
which are braided.\qed

\begin{equation} k_4*b_1=b_1^{-1}*k_4 \  \end{equation}
\proof[Proof of \rm(28)]
By the definition and by (M3) we have $k_4=d_3=(b_2a_2e_1b_1^{-1})*d_{1,3}$.
Conjugating $k_4$ and $b_1$ by $t_2^{-1}b_1e_1^{-1}a_2^{-1}b_2^{-1}$
we get $d_{1,2}$ and $b_1$, which are braided, by (\ref{A2}).\qed
\par

\begin{equation} g_5=sa_1^2k_5sa_1^2k_5^{-1} \end{equation}
\proof[Proof of \rm(29)]
By the definition $g_5=(rk_5rk_5^{-1})^2$ where $k_5=a_2d_{1,2}^{-1}t_1$. 
We shall prove 
that $r$ commutes with $k_5rk_5^{-1}$. Then $g_5=r^2k_5r^2k_5^{-1}=
sa_1^2k_5sa_1^2k_5^{-1}$, as required.
\par
$k_5rk_5^{-1}=a_2d_{1,2}^{-1}e_1a_1a_2e_1a_1b_1a_1e_1^{-1}a_1^{-1}a_2^{-1}
e_1^{-1}d_{1,2}a_2^{-1}=$ (by (J) and (\ref{A2}))\nl
$\strut a_2^2d_{1,2}^{-1}e_1a_1b_1a_1^{-1}e_1^{-1}d_{1,2}=$ (by (J))\nl
$\strut a_2^2d_{1,2}^{-1}b_1^{-1}a_1^{-1}e_1a_1b_1d_{1,2}$.
 \par\noindent
The last expression commutes with $a_1$ and $b_1$, by (J) 
and (\ref{A2})).
\qed

\begin{equation} g_6=(sa_1^2t_1)^4 \end{equation} 
\proof[Proof of \rm(30)]
By the definition $g_6=(ra_1t_1)^5$. The required relation is proved by 
a rather long computation, using (J). Observe first that 
$sa_1^2=(b_1a_1)^3$, 
and that $ra_1=(b_1a_1)^2$. We also have

$\strut t_1(b_1a_1)^2t_1=e_1a_1a_2e_1b_1a_1b_1a_1e_1a_1a_2e_1=$\nl
$\strut e_1a_1b_1a_2e_1a_1b_1a_1e_1a_1a_2e_1=b_1e_1a_1b_1a_2e_1a_1b_1a_1e_1a_2e_1=$
\nl
$\strut b_1a_1e_1a_1b_1a_2e_1a_1b_1a_1e_1a_2=b_1a_1e_1a_1a_2e_1b_1a_1b_1a_1e_1a_2=$
\nl
$(b_1a_1)t_1(b_1a_1)^2e_1a_2$\\
and\par
$e_1a_2(b_1a_1)^2t_1=b_1e_1a_1a_2b_1a_1e_1a_1a_2e_1=
b_1e_1a_1a_2e_1b_1a_1e_1a_2e_1=$\nl
$\strut b_1e_1a_1a_2e_1b_1a_2a_1e_1a_2=
b_1a_1e_1a_1a_2e_1b_1a_1e_1a_2=
(b_1a_1)t_1(b_1a_1)e_1a_2=$\nl
$\strut b_1a_1e_1a_1b_1a_2e_1a_1e_1a_2=b_1a_1b_1e_1a_1b_1a_2e_1a_1a_2=
b_1a_1b_1a_1e_1a_1b_1a_2e_1a_1=$\nl
$\strut (b_1a_1)^2t_1(b_1a_1)$.\par\noindent
The required result follows from the above relations.\endpf\par
This concludes the proof of Theorem \ref{simple presentation}.

\section{Mapping class group of a closed surface}
We shall consider in this section the mapping class group $\cM_g$ of
a closed surface $S_{g,0}$ of genus $g>1$. 
We shall keep the notation from the previous 
section. In particular ${\mathcal M}_{g,1}$ is the 
mapping class group of $S=S_{g,1}$ and
$S_{g,0}$ is obtained from $S$ by capping the boundary $\partial$
of $S$ by a disk $D$ with a distinguished center $p$, 
and then forgetting $p$.
We have two exact sequences 
$$ 1\to\ints\mapright{\psi}
\cM_{g,1}\mapright{\phi}\cM_{g,0,1}\to 1$$
$$ 1\to \pi_1(S_{g,0,1},p)\mapright{\sigma} \cM_{g,0,1}
  \mapright{e_*} \cM_{g,0,0}\to 1 $$
In the first sequence 
the Dehn twist $\Delta=T_\partial$ belongs to the kernel of $\phi$.
We shall prove now that it generates the kernel. When we split the surface
$S_{g,1}$ open along the curves $\beta_1,\alpha_1,\epsilon_1,\alpha_2,
\dots,\epsilon_{g-1},\alpha_g$ we get an annulus $N$, and one boundary 
of $N$ is equal
 to $\partial$. If $h\in ker(\phi)$ than $h$ takes each curve 
 $\gamma\in\{\beta_1,\alpha_1,\epsilon_1,\alpha_2,
\dots,\epsilon_{g-1},\alpha_g\}$ onto a curve $h(\gamma)$ which is 
isotopic to $\gamma$ in $S_{g,0}$ by an isotopy fixed on $p$. Therefore
$\gamma$ and $h(\gamma)$ form 2--gons which are disjoint from $p$ and
hence from $D$. It follows that $h$ is isotopic in $S_{g,1}$ to a 
homeomorphism equal to the identity on all curves
$\beta_1,\alpha_1,\epsilon_1,\alpha_2,
\dots,\epsilon_{g-1},\alpha_g$. But then it is a homeomorphism of 
the annulus $N$ so it is isotopic to a power of $\Delta$.

\begin{figure}[ht!]
\relabelbox\small
\epsffile{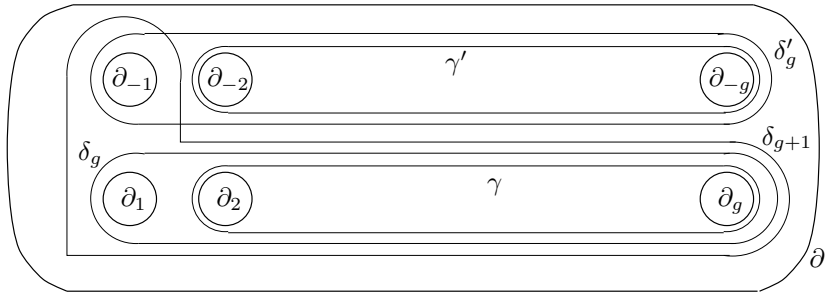}
\relabel {del}{$\partial$}
\relabela <-3pt,3pt>  {d1}{$\partial_1$}
\relabela <-3pt,3pt>  {d2}{$\partial_2$}
\relabela <-3pt,3pt>  {delg}{$\partial_g$}
\relabela <-6.5pt,3pt>  {dm1}{$\partial_{-1}$}
\relabela <-6.5pt,3pt> {dm2}{$\partial_{-2}$}
\relabela <-6.5pt,3pt> {dmg}{$\partial_{-g}$}
\relabela <-1pt,0pt> {dg}{$\delta_g$}
\relabela <0pt,-1pt> {dg1}{$\delta_{g+1}$}
\relabel {dgp}{$\delta_g'$}
\relabel {g1}{$\gamma$}
\relabel {ga}{$\gamma'$}
\endrelabelbox
\caption{Spin maps in the proof of Theorem 3}
\label{spin maps}
\end{figure}

The second sequence is described in \cite{Birman}, Theorem 4.3. The 
kernel of $e_*$ is generated by spin-maps
$T_{\gamma^\prime}T_\gamma^{-1}$, where $\gamma$ and $\gamma^\prime$
are simple, nonseparating curves which bound an annulus on $S_{g,0,1}$
 containing the distinguished point $p$.
The composition $e_*\phi$ is an epimorphism from the group 
$\cM_{g,1}$ onto the group
$\cM_g=\cM_{g,0,0}$ and its kernel is generated by $\Delta$ and 
by the spin maps 
 $T_{\gamma^\prime}T_\gamma^{-1}$, where  $\gamma$ and $\gamma^\prime$
are simple, nonseparating curves on $S$ which bound an annulus with a hole
bounded by $\partial$. Clearly all such annuli are equivalent by 
a homeomorphism of $S$, hence all spin maps are conjugate in 
${\mathcal M}_{g,1}$. 
It suffices
to consider one spin map $T_{\delta_g^\prime}T_{\delta_g}^{-1}$, where 
$\delta_g$ and $\delta_g^\prime$ are curves on Figure \ref{spin maps}.
$T_{\delta_g}$ is equal to the element $d_g$ in the relation (M4).

Let $w=b_1a_1e_1a_2\dots a_{g-1}e_{g-1}a_g^2e_{g-1}a_{g-1}\dots a_2e_1a_1b_1$. 
It is easy
to check, drawing pictures, that $w(\delta_g)=\delta_g^\prime$. Therefore,
by Lemma \ref{conjugation}, relation 
(M4) is equivalent,
 modulo 
relations in ${\mathcal M}_{g,1}$, to $T_{\delta_g}=T_{\delta_g^\prime}$. 
By the above 
argument $\cM_g$ has a presentation with relations (M1) -- (M4)
and relation $\Delta=1$. We have to prove that the last relation follows
from the others.\par
Let $\cM^\prime=({\mathcal M}_{g,1})/(M4)$.
Let $d_g$, $d_g^\prime$, $d_{g+1}$, $c$, $c^\prime$ be twists along curves
$\delta_g$, $\delta_g^\prime$, $\delta_{g+1}$, $\gamma$, $\gamma^\prime$
respectively, depicted on Figure \ref{spin maps}. Each element of 
${\mathcal M}_{g,1}$ represents an element in $\cM^\prime$ which we denote 
by the same symbol.

 From now on, till the end of this section, symbols denote elements in
 $\cM^\prime$. We want to prove that $\Delta=1$.

 All relations  from Lemma \ref{basic relations} are true in $\cM^\prime$.
  We have $d_g*b_1=b_1^{-1}*d_g$ and $d_g$ commutes with
every $a_i$ and $e_i$. By Lemma \ref{basic relations}, (iii) and 
 (v) we have: 

$\Delta=(a_ge_{g-1}a_{g-1}\dots e_1a_1b_1d_g)^{2g+2}=$\nl
$\strut(a_ge_{g-1}a_{g-1}\dots e_1a_1)^{2g}(b_1a_1\dots a_ga_g\dots a_1b_1)
(d_gb_1a_1\dots a_ga_g\dots a_1b_1d_g)$,\nl
$\strut d_gd_g^\prime=(a_ge_{g-1}a_{g-1}\dots e_1a_1)^{2g}=$ \nl
$\strut(a_ge_{g-1}a_{g-1}\dots e_2a_2)^{2g-2}
(e_1a_2\dots a_ga_g\dots a_2e_1)
(a_1e_1a_2\dots a_ga_g\dots a_2e_1a_1)$,\nl
$\strut(a_ge_{g-1}a_{g-1}\dots e_2a_2)^{2g-2}=cc^\prime$,\nl
$\strut(d_gb_1a_1)^4=cd_{g+1}$.

We also see that $c^\prime d_{g+1}^{-1}$ and $d_g^\prime d_g^{-1}$ are 
spin maps, hence $c^\prime=d_{g+1}$ and $d_g=d_g^\prime$.
Therefore\par
$(a_ge_{g-1}a_{g-1}\dots e_2a_2)^{2g-2}=(d_gb_1a_1)^4$,\nl
$\strut d_g^2=(a_ge_{g-1}a_{g-1}\dots e_2a_2)^{2g-2}
(e_1a_2\dots a_ga_g\dots a_2e_1)
(a_1e_1a_2\dots a_ga_g\dots a_2e_1a_1)$,\par hence\par
$(a_1e_1a_2\dots a_ga_g\dots a_2e_1a_1)=(e_1a_2\dots a_ga_g\dots a_2e_1)^{-1}
(d_gb_1a_1)^{-4}d_g^2$.\\
Further\par
$\Delta=(a_ge_{g-1}a_{g-1}\dots e_1a_1)^{2g}(b_1a_1\dots a_ga_g\dots a_1b_1)
(d_gb_1a_1\dots a_ga_g\dots a_1b_1d_g)=$ \nl
$\strut d_g^2b_1(a_1e_1a_2\dots a_ga_g\dots a_2e_1a_1)b_1
(d_gb_1a_1\dots a_ga_g\dots a_1b_1d_g)=$ \nl
$\strut d_g^2b_1(e_1a_2\dots a_ga_g\dots a_2e_1)^{-1}
(d_gb_1a_1)^{-4}d_g^2b_1(d_gb_1a_1\dots a_ga_g\dots a_1b_1d_g)$.
 \par\noindent
 Now $(d_gb_1a_1\dots a_ga_g\dots a_1b_1d_g)$ commutes with $d_g$, by (M4),
 and commutes with $b_1$ and $a_1$, by (J), hence\par 
$\Delta=$\nl
\cl{$\strut d_g^2b_1(e_1a_2\dots a_ga_g\dots a_2e_1)^{-1}(d_gb_1a_1)^{-1}
(d_gb_1a_1\dots a_ga_g\dots a_1b_1d_g)
(d_gb_1a_1)^{-3}
d_g^2b_1$}\nl
$\strut =d_g^2b_1a_1b_1d_g(a_1^{-1}b_1^{-1}d_g^{-1})^3d_g^2b_1=1$, by (J).  
\par
 This concludes the proof of Theorem \ref{presentation closed}.

\section{Equivalence of presentations} 

In this section we shall prove that the presentations of $\cM_{g,1}$ in
Theorems 1 and $1^\prime$ are equivalent. The relations 
(M1) coincide with the relations (A). It follows from relations (A) that
$b_2$ commutes with the left hand side of the relation (B). Thus (B) is 
equivalent to \par
$(b_2a_2e_1a_1b_1^2a_1e_1a_2b_2)
(a_2e_1a_1b_1^2a_1e_1a_2)^{-1}=(b_1a_1e_1)^4$.

\noindent
Multiplying by $(a_2e_1a_1b_1^2a_1e_1a_2)$ on the right
we get \par
$(b_2a_2e_1a_1b_1^2a_1e_1a_2b_2)=(b_1a_1e_1)^4(a_2e_1a_1b_1^2a_1e_1a_2)=
(b_1a_1e_1a_2)^5$, as in the proof of  Lemma
\ref{basic relations} (v), and we get a relation identical with (M2).

We now pass to relation (M3). We shall transform it using relations 
(M1) and (M2) and then we shall conjugate it by $w=a_3e_2a_2e_1a_1b_1$
to get the relation (C). Since (M1)=(A) and (M2) is equivalent to (B) 
in presence of (M1),
it will prove that (M3) is equivalent to (C) in presence of (M1) and (M2).
 It follows from (M1) and the definitions that
each factor on the right hand side of (M3) commutes with $a_1a_2a_3$,
therefore $d_3$ also commutes with $a_1a_2a_3$. 
 Recall the relations (\ref {A2}),
(\ref{ti,tj}) and (\ref{A14}) from section 4, which follow from
 the relations (M1) and (M2).\par
(\ref{A2}) \ \ $d_{1,2}t_1=t_1d_{1,2}$,\par
(\ref{ti,tj}) \ \ $t_1t_2t_1=t_2t_1t_2$,\par
(\ref{A14}) \ $t_2d_{1,2}t_2d_{1,2}=d_{1,2}t_2d_{1,2}t_2$.\par
\noindent We now have  

$d_{1,2}d_{1,3}d_{2,3}=
d_{1,2}t_2d_{1,2}t_2^{-1}t_1t_2d_{1,2}t_2^{-1}t_1^{-1}=$ (by \ref{ti,tj})
\nl
$\strut d_{1,2}t_2d_{1,2}t_1t_2t_1^{-1}d_{1,2}t_2^{-1}t_1^{-1}=$ (by \ref{A2})
$d_{1,2}t_2t_1d_{1,2}t_2d_{1,2}t_1^{-1}t_2^{-1}t_1^{-1}=$ (by \ref{ti,tj})
\nl
$\strut d_{1,2}t_2t_1d_{1,2}t_2d_{1,2}t_2^{-1}t_1^{-1}t_2^{-1}=$ (by \ref{A14})
$d_{1,2}t_2t_1t_2^{-1}d_{1,2}t_2d_{1,2}t_1^{-1}t_2^{-1}=$ (by \ref{ti,tj}
and \ref{A2}) \nl
$\strut t_1^{-1}d_{1,2}t_2t_1d_{1,2}t_2d_{1,2}t_1^{-1}t_2^{-1}=$ (by \ref{A2})
$t_1^{-1}d_{1,2}t_2d_{1,2}t_1t_2t_1^{-1}d_{1,2}t_2^{-1}=$ (by \ref{ti,tj}
and \ref{A14}) \nl
$\strut t_1^{-1}t_2^{-1}d_{1,2}t_2d_{1,2}t_1t_2d_{1,2}t_2^{-1}=$ (by \ref{A2}
and \ref{A14}) 
$t_1^{-1}t_2^{-1}d_{1,2}t_2t_1t_2^{-1}d_{1,2}t_2d_{1,2}$.

We now conjugate everything by $w$ and, using (M1), we  get \par
$w*a_1=b_1$, \ $w*e_1=a_1$, \ $w*a_2=e_1$, \ $w*e_2=a_2$, \ $w*a_3=e_2$,
\nl
$\strut w*t_1=a_1b_1e_1a_1=\tilde t_1$, \ $w*t_2=a_2e_2e_1a_2=\tilde t_2$,
 \ $w*d_{1,2}=b_2$.
\par\noindent
Therefore after conjugation by $w$  the right hand side of (M3) becomes 
the right hand side of (C).

We have shown in the proof of relations (\ref{A15}), Case 1c, 
using only relations (M1), that

$d_3=(b_1^{-1}((a_2e_1e_2a_2a_3e_2)^{-1}*b_2)
b_1^{-1}a_1^{-1}e_1^{-1}a_2^{-1})*b_2$.

When we conjugate the last expression by $w$
 we get exactly the expression for $\tilde d_3$ in 
Theorem $1^\prime$.\par
This proves the equivalence of 
the presentations in Theorems 1 and $1^\prime$.\par

In order to compare Theorems 3 and $3^\prime$ we need another set 
of generators. Let us call the curves $\beta_2,\beta_1,\alpha_1,
\epsilon_1,\alpha_2\dots,\epsilon_{g-1},\alpha_g$ --- {\em the generating 
curves}. Let $\beta_g$ be the curve shown on Figure \ref{general surface}
and let $\beta_2^\prime$ be a curve which intersects $\epsilon_{g-2}$
once and intersects $\beta_g$ once and is disjoint from the other
{\em generating curves}. Then the curves $\beta_2^\prime,
\alpha_g,\epsilon_{g-1},\alpha_{g-1},\dots,\epsilon_1,\alpha_1,\beta_1$
have the same intersection pattern as the {\em generating curves}
and the curve $\beta_g$ plays the same role with respect to these 
curves as the curve $\delta_g$ with respect to the {\em generating curves}.
Let $b_g$ and $b_2^\prime$ be twists along the curves $\beta_g$ and
$\beta_2^\prime$ respectively. Then, by Theorem 1, we have a new 
presentation of ${\mathcal M}_{g,1}$ with generators 
$b_2^\prime,a_g,e_{g-1},a_{g-1},\dots,a_1,b_1$ and with defining relations
(M1$^\prime)$, (M2$^\prime)$, (M3$^\prime)$ corresponding to 
(M1), (M2), (M3). It is a presentation of the same group and therefore, 
when we express $b_2^\prime$ in terms of the generators from Theorem 1,
 it is equivalent to 
the presentation ((M1), (M2), (M3)) and to the presentation
((A), (B), (C)). By Theorem 3 the group ${\mathcal M}_{g,0}$ 
has a presentation
with relations (M1$^\prime)$, (M2$^\prime)$, (M3$^\prime)$ and one more
relation
 
\noindent 
(M4$^\prime)$ \ \ $[a_ge_{g-1}a_{g-1}\dots e_1a_1b_1^2a_1e_1\dots 
a_{g-1}e_{g-1}a_g,b_g]=1.$

Here $b_g$ is some product of generators which represents the Dehn twist 
of $S_{g,1}$ along the curve $\beta_g$. All such products are equivalent 
modulo relations (M1$^\prime)$, (M2$^\prime)$, (M3$^\prime)$.
Relation (D) of Theorem $3^\prime$ has the same form with $b_g$ replaced by
$\tilde d_g$.
Therefore in order to check that the presentations
in Theorems 3 and $3^\prime$ are equivalent it suffices to prove that
the expression for $\tilde d_g$ in (D) also represents the Dehn twist 
with respect to the  curve $\beta_g$.
This task (of drawing very many pictures) is left to the reader.

\end{document}